%% file: goodnonsimple_archiv
\input amstex
\input generic_macros
\input papermacros_entirefunctions

\input diagrams

\def\diag{\text{{\rm diag}}\,}

\def\Q{\text{\bf Q}}
\def\Inf{\text{Inf\,}}
\def\Ann{\text{Ann\,}}
\def\paren #1{\/{\rm(}#1\/{\rm)}}

\def\cvx{\text{cvx\,}}
\long\def\Ex #1.#2\par{\vskip 2pt \noindent {\it #1} #2}

\font\sm=cmr5 at 5pt
\def\Gd{\subset {} \hskip -1.05em  \raise .8pt \hbox{ \text{{\sm G}}}\,\ }
\def\Ideal{\triangleleft}
\let\ideal=\Ideal
\def\Asp{\text{Aspan\,}}

\let\hat=\widehat

\def\Aff{\text{Aff\,}}
\NoBlackBoxes

\def\oneone{1.1}
\def\onetwo{1.2}
\def\onethr{1.3}
\def\onefou{1.4}
\def\onefiv{1.5}
\def\onesix{1.6}
\def\onesev{1.7}

\def\twoone{2.1}
\def\twotwo{2.2}
\def\twothr{2.3}
\def\twofou{2.4}
\def\twofiv{2.5}

\def\threig{4.8}

\def\fouone{3.1}
\def\foutwo{3.2}
\def\fouthr{3.3}
\def\foufou{3.4}
\def\foufiv{3.5}
\def\fousix{3.6}
\def\fousev{3.7}

\def\fivone{5.1}
\def\fivtwo{5.2}
\def\fivthr{5.3}
\def\fivfou{5.4}
\def\fivfiv{5.5}
\def\fivsix{5.3}
\def\fivsev{5.6}
\def\fiveig{5.8}
\def\fivnin{5.4}
\def\fivten{5.5}

\def\sixone{6.1}
\def\sixtwo{6.2}
\def\sixthr{6.2}
\def\sixfou{6.3}
\def\sixfiv{6.4}
\def\sixsix{6.5}
\def\sixsev{6.6}
\def\sixeig{6.7}
\def\sixnin{6.8}
\def\sixten{6.9}
\def\sixele{6.10}

\let\iso=\cong

\def\Int{\text{int}\,}

\def\tripnorm #1xxx{\left\|\hglue-.2ex\left|#1\right|\hglue-.2ex\right\|}

\def\diag{\text{diag}\,}

\Title Good measures for non-simple dimension groups

\Abstract Akin's notion of good measure, introduced to classify measures on Cantor sets has been  translated to dimension groups and corresponding traces by Bezuglyi and the author, but emphasizing the simple (minimal dynamical system) case. Here we deal with non-simple (non-minimal) dimension groups. In particular, goodness of tensor products of large classes of non-good traces (measures) is established. We also determine the pure faithful traces on the dimension groups associated to xerox type actions on AF C*-algebras; the criteria turn out to involve  algebraic geometry and number theory.
{\par}We also deal with a coproduct of dimension groups, wherein, despite expectations, goodness of direct sums is nontrivial. In addition, we verify a conjecture of [BeH] concerning good subsets of Choquet simplices, in the finite-dimensional case.

\noindent  {\it David Handelman}%
\plainfootnote{$^1$}{\hglue -.5em Supported in part by a Discovery grant from NSERC.}

\SecT Introduction \& definitions

Akin [Ak1, Ak2, ADMY, \dots] introduced and studied the notion of good measures, in connection with the
classification of (probability) measures on Cantor sets up to
homeomorphism. With the development in [Pu, HPS, GPS, etc],  of classification and
construction of minimal actions \wrt strong orbit and orbit equivalence
via Vershik maps and ordered Grothendieck groups of AF C*-algebras, this
and related properties were translated into the language of (traces on)
dimension groups (a class of partially ordered abelian groups) in [BeH].
In particular, the characterizations therein of goodness of traces on
simple dimension groups provided relatively easy constructions of good and
non-good measures on minimal systems. For more details, see the discussion
in the introduction to [BeH].

Recent work (e.g., [FO, P]) has extended  Vershik action(s) to non-minimal
systems, and correspondingly to non-simple dimension groups. Here we give
computable criteria for goodness in the general (approximately divisible)
case, and then use the criteria to  give a surprising result that  tensor
products of (some) non-good traces are good; this applies to  the ugly
traces of [BeH]. We also completely determine the pure faithful traces on
fixed point algebras under xerox actions of tori: these include Pascal's
triangle and variations corresponding to spatially and temporally
homogeneous random walks with finite support on the lattice $\Z^d$.

From [H1,
Theorem III.3], the pure faithful  traces correspond to points $r = (r_i)$ in
the strictly positive orthant of $\R^d$; those that are good are precisely
the ones that satisfy two number-theoretic conditions, which in the case
$d=1$ reduce to (i) no other algebraic conjugate of $r = r_1$ is positive
and (i) if the leading and terminal coefficients of the polynomial
implementing the random walk are $a_0$ and $a_k$, then there exists $s$
\st $a_0^s/r$ and $a_k^s r$ are both algebraic integers.

We also deal with a strict form of direct sum of dimension groups,
determining when the corresponding sum of traces is good; there are some
surprises here, as the direct sum can be good without either one being
good (in fact, we find for each $m$, a collection of simple dimension
groups with traces, $(G_i,\tau_i)$ \st that for any strict direct sum of
$m$ or fewer distinct summands, $\oplus_{i \in S} G_i$, the sum of the
traces is not good, but for any direct sum of more than $m$ direct
summands, the sum is good.

We then consider good sets of traces. The first problem is the definition;
it should be consistent with the current definition in the simple case,
and in the singleton case, and we discuss various possibilities; finally,
we settle on one. We show that for the class of dimension groups
considered above (arising from random walks on $\Z^d$), with any
reasonable definition, the notion is surprisingly restrictive, and even
order-unit goodness turns out to be sensitive to the Newton polyhedra of
the polynomials (unlike the case for single traces).

There are two appendices. The first characterizes order unit good traces on simplicial dimension groups, and the resulting characterization suggests that there are no effective  for goodness involving order unit goodness when there are discrete traces, in contrast to the approximately divisible situation discussed in the rest of this aritcle. The second appendix verifies, in the finite-dimensional trace space case, a conjecture made in [BeH, section 7] concerning the structure of good subsets relative to a simplex.

\noindent{\it Definitions.} A partially ordered abelian group $G$ with positive cone $G^+$ is {\it unperforated\/} if whenever $ n$ is a positive integer and $g \in G$, then $ng \in G^+$ entails $g \in G^+$. An {\it order unit\/} for $G$ is an element $u \in G^+$ \st for all $g \in G$, there exists a positive integer $K$ \st $-K u \leq g \leq Ku$. A {\it trace\/} (formerly, {\it state\/})
is a nonzero positive group homomorphism $\Arrow \tau; G.\R$; if $\tau(u) =1$ and $u$ is an order unit, we say $\tau$ is normalized (\wrt $u$). The trace $\tau$ is {\it faithful\/} if $\ker \tau \cap G^+ = \brcs{0}$ (this is {\it much\/} weaker than being one to one, and corresponds to faithful measure).

When $(G,u)$ is a partially ordered abelian group with order unit, we may form $S(G,u)$, the compact convex set of normalized traces, equipped with the weak (or point-open) topology. We denote by $\Aff S(G,u)$ the Banach space of continuous convex-linear (affine) real-valued functions on $S(G,u)$. There is a natural representation  $G \to \Aff S(G,u)$, given by $g \mapsto \hat g$, where $\hat g(\tau) = \tau(g)$.

If $(G,u)$ is an unperforated ordered abelian group, we say $G$ is {\it approximately divisible\/} if its range in $\Aff S(G,u)$ is norm-dense; for dimension groups with order unit, this is equivalent to $\tau(G)$ being dense in $\R$ for all pure traces $\tau$, or equivalently, for all order units $g \in G$, there exist order units $a, b$ of $G$ \st $g = 2a + 3b$ (and there are many other equivalent formulations).

When $I$ is a subgroup (typically an order ideal) of a partially ordered abelian group $G$, we say {\it $I$ has its own order unit $w$\/} or {\it $w$ is a relative order unit of $I$\/} if $w \in I$ is an order unit of $I$ \wrt the relative ordering inherited from $G$. This is to emphasize the fact that $w$ is {\it not\/} an order unit for $G$, merely for $I$.

If $G$ is an unperforated ordered abelian group, we say $G$ is {\it nearly divisible\/} if for every order ideal $(I,w)$ which has its own order unit, $(I,w)$  is approximately divisible; an equivalent form that does without the order ideals is that  for all $g \in G^+$,
there exists $a,b \in G^+$ \st $g = 2a + 3b$ and $g \leq ka , kb$ for some positive integer $k$).

For example, if $G = H \otimes U$ where $H$ is a partially ordered unperforated abelian group and $U$ is a noncyclic subgroup of the rationals, $\Q$,
then $G$ is nearly divisible, and it is approximately divisible if it has an order unit. We will see plenty of nearly divisible examples that are not of this type in later sections.

A trace on $G$ is {\it discrete\/} if its image $\tau(G)$ is a cyclic (that is, discrete) subgroup of $\R$. An alternative characterization of approximately divisible, for dimension groups, is that $(G,u)$ admit no discrete traces; for nearly divisible, the characterization is that no nonzero order ideal with order unit admits a discrete trace.

For general relevant results on partially ordered abelian groups, especially dimension groups, see [G].

An {\it interval\/} in a partially ordered group $G$, is a subset of the form
$[0,b]:= \Set{g \in G}{0 \leq g \leq b}$ for some $b \in G^+$.

Following [BeH], and based on Akin's notion for measures on Cantor sets, a trace $\Arrow \tau;G.\R$ is {\it good\/} (as a trace of $G$) if for all $b \in G^+$, $\tau([0,b]) = [0,\tau(b)]$, that is, if $a' \in G$ and $0 \leq \tau(a') \leq \tau(b)$, there exists $a \in [0,b]$ \st $a - a' \in \ker \tau$. If $(G,u)$ is a partially ordered abelian group with order unit, we say $\tau$ is {\it order unit good\/} if  in the definition of good, we restrict $b$ to be an order unit.

\SecT 1 Characterization of goodness

Order unit goodness is relatively easy to characterize when $(G,u)$ is approximately divisible [BeH, Proposition 1.7]: $\tau$ is order unit good iff the image of $\ker \tau$ in $\Aff S(G,u)$ is dense in $\tau^{\vdash}:= \Set{h \in \Aff S(G,u)}{h(\tau) = 0}$ (the latter is closed and codimension one subspace of $\Aff S(G,u))$. This makes examples and non-examples relatively easy to construct. There is a corresponding characterization for goodness, which we shall simplify a bit, and used to actually do something.

\Lem Proposition \oneone. Suppose $(G,u)$ is a  dimension group with order unit. Let $\tau$ be a faithful trace of
$G$. Then $\tau$ is good iff for all nonzero order ideals with order unit $(I,w)$,
both $\tau(I) = \tau(G)$ and $\tau|I$ is order unit good.

\Rmk Necessity  is shown in [BeH, Proposition 4.2]; although the statement hypothesizes that  $\tau$ be  pure, this is not used in the proof; also shown there
was that if $\tau$ is good, then $\tau|I$ is good (as a trace on the order ideal $I$), and
this implies  (in the case that $I$ is approximately divisible) that $\tau|I$ is   order unit good,
just from the definitions.

\Rmk It is always possible to reduce to the case that $\tau$ be faithful, by factoring out the maximal order ideal $J$ contained in $\ker \tau$ [BeH, Lemma 4.4]. In this case, the criteria apply to  $G/J$ (replacing $G$). This would make the statement somewhat more complicated.

\Pf Proof of necessity is given in [BeH; Proposition 4.2], requiring neither purity of $\tau$ nor approximate divisibility.

Conversely, suppose $a \in G$, $b \in G^+$ and $0 < \tau(a) < \tau(b)$. Form the order ideal $I $ generated
by $b$, that is, $I = \Set{c \in G}{\exists N \in \N \text{ \st } -N b \leq g \leq Nb}$. Then $I$ is an order
ideal with its own order unit, $b$. Since $\tau(I) = \tau(G)$, there exists $a_1 \in I$ \st $\tau(a_1) = \tau(a)$.
Now order unit goodness of $\tau|I$ yields $a' \in I$ \st $\tau(a') = \tau(a_1) = \tau(a)$ and $0 \leq a' \leq b$, verifying
goodness of $\tau$.
\qed

\comment
Next, it suffices to show that if $I$ and $J$ are approximately divisible
order ideals in a dimension group and $\tau$ is a trace on $I+J$ (itself
an order ideal in $G$), and the latter is approximately divisible (this
may be a consequence of other stuff), and the restrictions $\tau|I$ and
$\tau|J$ are order unit good, then $\tau|(I+J)$ is order unit good. We
need something like, every order unit $a$ of $I+J$ is a sum of $u+v$ where
$u $ is an order unit for $I$ and $v$ is an order unit of $J$. [Then the
usual proof works.]
\endcomment

Let $G$ be a dimension group, and let $I$ and $J$ be order ideals thereof. Then $H:=
I+J$ (the set of sums of elements in $I$ and $J$) and $I \cap J$ are both
order ideals. Most of the following are variations on [BeH, Lemma 1.3]. As in [BeH], and element $v$ of $G^+$ is $\tau$-good or $\tau$-order unit good if $\tau([0,v]) = [0,\tau(v)]$.

\Lem Lemma \onetwo. Suppose $G$ is a dimension group, and $I$ and $J$ each have
(relative) order units, $w,y$ respectively. Then
\item{(a)} $I +J$ is an order ideal of $G$ with a (relative) order unit.
\item{(b)} Let $\tau$ be a trace on $G$ \st $\ker \tau \cap G^+ =
\brcs{0}$ and $\tau(I) \cap \tau(J)$ is dense in $\R$. If $\tau|I$ and
$\tau|J$ are good (as traces on $I$ and $J$ respectively), then $\tau$ is
good.
\item{(c)} If $I+J$ is approximately divisible, then every order unit $b$ of $I+J$
can be written in the form $b = u + v$ where $u,v$ are relative order
units for $I$, $J$ respectively.
\item{(d)} If $v$ is $\tau$-order unit good (\wrt $I$) and $w$ is $\tau$-order unit good (\wrt $J$), and $\tau(I ) \cap \tau(J)$ is dense in $\R$, then $v + w $ is $\tau$-order unit good \wrt $I+J$.
\item{(e)}Suppose each of $I$, $J$ and
$I+J$ are approximately divisible, and $\tau$ is a trace on $I+J$ \st each
of $\tau|I$ and $\tau|J$ is order unit good, and $\tau(I) \cap \tau(J)$ is
dense in $\R$. Then $\tau$ is order unit good as a trace of $I+J$.

\Rmk Part (c) can fail if approximate divisibility is dropped; for example, take $G = \Z^3$ with the usual simplicial ordering, let $I$ be the order ideal generated by $(1,1,0)$ and let $J$ be the order ideal generated by $(0,1,1)$; then $I + J = G$ and the order unit $(1,1,1)$ cannot be realized as a sum of relative order units from $I$ and $J$ respectively.

\Pf (a) That $I+J$ is an order ideal is ancient, e.g., [G]. If $w$ and $y$
are respective order units for $I$ and $J$, then $z := w+y$ is an order
unit for $I+J$. To see this, let $f \in (I+J)^+$; for dimension groups,
$(I+J)^+= I^+ + J^+$, hence we can find $e \in I^+$ and $g \in J^+$ \st $f
= e+g$. Since there exist positive integers $k$, $k'$ \st $e \leq kw$ and
$g \leq k'v$, we have $f \leq k'' z$ where $k'' = \max\brcs{k,k'}$.

\noindent (b) Select $b \in G^+$ and $a \in G$ \st $\tau (a) < \tau (b)$. We may
write $b = i + j$ where $i \in I^+$ and $j \in J^+$ (see [G]). Then
$\tau(i), \tau(j) > 0$. We may write $\tau(a) =r + s$ where $r \in
\tau(I)$ and $s \in \tau(J)$.

Assume $\tau(a) \geq \tau(i)$. By density of $\tau(I) \cap \tau(J)$, given $0 <
\epsilon< \min\brcs{\tau(i), \tau(b) - \tau(a)}$, there exists $\delta \in
\tau(I) \cap \tau( J)$ \st $\tau(i) - \epsilon < r + \delta < \tau(i)$.
Then $s - \delta = \tau(a) -r - \delta $ satisfies
$$
\tau(a)- \tau(i)+ \epsilon > s-\delta > \tau(a) - \tau(i) > 0
$$
Hence we can write $\tau (a) = (r + \delta) + (s-\delta)$, where the
parenthesized terms are respectively in the intervals $(0,\tau(i))$ and
$(0,\tau(a) - \tau(i) + \epsilon)$. However, $\epsilon < \tau(b) - \tau
(a)$ entails $\tau(a) - \tau(i) + \epsilon < \tau(b) - \tau(i) =
\tau(j)$. Since $\pm \delta \in \tau(I \cap J)$, we may thus find $a_1 \in
I$ and $a_2 \in J$ \st $0 < \tau(a_1) < \tau(i)$ and $0 < \tau(a_2) <
\tau(j)$. Since each of $\tau|I$ and $\tau|J$ is good, there exist $c_1
\in [0,i]$ (the interval in $I$) and $c_2 \in [0,j]$ \st $\tau(c_1) < \tau
(i) $ and $\tau(c_2) < \tau(j)$. Hence we have $c:= c_1 + c_2 \in [0,b]$
and $\tau(c) = \tau (c_1) + \tau(c_2) < \tau(i) + \tau (j) = \tau (b)$,
verifying goodness in this case.

Reversing the roles of $i$ and $j$, the same conclusion results if
$\tau(a) \geq \tau(j)$, so we are reduced to the case that $\tau(a) <
\min\brcs{\tau(i), \tau(j)}$. If $\tau(a) = 0$, there is nothing to do
(except set $c = 0$). Otherwise, choose $
0 < \epsilon < \tau(a)/2$ find real $ \delta \in \tau(I \cap J)$ \st
$\tau(a)/2 - \epsilon < \delta + r < \tau(a)/2$, and consider $\tau(a) =
(r + \delta) + (s-\delta)$; then $r + \delta \in (0, \tau(a)/2) \subset
(0, \tau(i))$, so $s- \delta \in (\tau(a)/2, \tau(a)) \subset
(0,\tau(j))$. Now we can proceed as in the previous paragraph.

\noindent (c) Now let $b$ be an order unit of $I+J$. By approximate divisibility of
$I+J$, the range of $I+J$ in $\Aff S(I+J,b)$ is dense; hence given
$\epsilon > 0$, we may find $b_0 \in I+J$ \st $(1/2 - \epsilon)\pmb 1 <
\hat b_0 < \pmb 1/2$ (where $\ \hat{ }\ $ refers only to the representation on
$S(I+J,b)$, that is, $\hat b = \pmb 1$). Let $\epsilon < 1/8$, so that
$\hat b_0 \gg 0$ and thus $b_0$ is an order unit of $I+J$, and moreover,
$2b_0 \leq b$, and $b-b_0$ is also an order unit for $I+J$.

Now consider the set $\Cal S:=\Set{c \in I^+}{c \leq b_0}$. This is
directed, as if $c,c' \in \Cal S$, then we have $c,c' \leq b_0, c+c'$;
interpolating, we obtain $c''$ \st $c,c' \leq c'' \leq b_0, c+c'$; as
$c+c' \in I$, it follows that $c'' \in I$, so $c'' \in \Cal S$. As there
exists $k$ \st $w \leq kb_0$, we can write $w = \sum_{i=1}^k w_i$ where
$w_i \in I^+$ and each $w_i \leq b_0$. Then $w_i \in \Cal S$, so there
exists $u_0 \in I^+$ \st $w_i \leq u_0 \leq b_0$ for all $i$. Since $\sum w_i
= w$ is an order unit for $I$, $ku_0$ is an order unit for $I$, and thus $u_0$
is too. Hence there exists an order unit $u_0$ of $I$ \st $u_0 \leq
b_0$.

Since $b-b_0$ is also an order unit for $I+J$, applying the same process
to $J$ instead of $I$ yields an order unit $v_0$ of $J$ \st $v_0 \leq b -
b_0$. Thus $u_0 + v_0 \leq b_0 + (b-b_0) = b$. The element $b-(u_0 + v_0)$
is in the positive cone of $I+J$, so can be written $b - (u_0 + v_0) = c +
d$ where $c \in I^+$ and $d \in J^+$.
This yields $b = (u_0 + c) + (v_0 + d)$; setting $u = u_0 + c$, we see
that $u \in I^+$ and is larger than an order unit for $I$, so is itself an
order unit for $I$; similarly $v = v_0 + d$ is an order unit for $J$.

\noindent (d) \& (e) Select an order unit $b$ for $I+J$, and $a \in I+J$ \st $0 < \tau(a) <
\tau(b)$. By (c), we may write $b = u + v$ where $u$ and $v$ are order units for $I$ and
$J$ respectively. We can write $a = r + s$ where $r \in I$ and $s\in J$,
and set $t = \tau(u)$ (as $\tau|I$ is order unit good, it does not vanish
identically, hence $ t > 0$), so that $\tau(v) = \tau(b ) - t$, which is
again positive. Now proceed as in the proof of (b).
\qed

\comment
Other possibilities: $I \cap J = 0$ (the intersection is always an order
ideal); in this case, $(I+J)/(I \cap J) \iso I \oplus J$ with the direct sum
(coordinatewise) ordering; then the question boils down to $\sigma \oplus
\tau$ on $I \oplus J$ being good if $\sigma$ and $\tau$ are on the
corresponding dimension groups. Hasn't this been done? [Note that $\tau(I
\cap J) = 0$ implies $I\cap J = 0$ because of our assumption that $\tau$
kills no positive elements and the fact that $I \cap J$ is an order
ideal.]

Next, $\tau(I \cap J) $ is a discrete subgroup. We can eliminate this
possibility by saying that all order ideals with order unit are
approximately divisible (so if $I \cap J \neq 0$, and $a$ is a nonzero
positive element therein, $\tau(\langle a \rangle)$ is a dense subgroup of
$\R$), but it may not be necessary. Since $\tau|I$ is good, its
restriction to any order ideal is good; but a good trace with discrete
value group either has a nonzero positive element in its kernel (which is
ruled out here) or the intersection is just a simplicial group (which
forces the trace to be given by $(1,1,\dots,1, 0,0,\dots)$, hence just by
$(1,1,1,\dots,1)$.

Seems to boil down to the easier condition $\tau(I) \cap \tau(J)$ is
dense/discrete/zero. With density, everything is fine (same argument);
with zero, there is a problem: the sum is not good! Also, really need
density (even discreteness).
\endcomment

The density requirement on  $\tau(I) \cap \tau(J)$ is essential.

\Lem Lemma \onethr. Suppose that $u$ and $v$ are elements of $G^+$, and let $\tau$
be a trace \st each is $\tau$-order unit good on the order ideals they
generate, $I(u)$ and $I(v)$ respectively.
\item{(a)} If $u+v$ is $\tau$-order unit
good on $I(u) + I(v) = I(u+v)$ and $\tau(I(u)) + \tau(I(v))$ is dense in
$\R$, then $\tau(I(u))\cap \tau(I(v)) \neq \brcs{0}$;
\item{(b)} if additionally,
both $\tau(I(u))$ and $\tau(I(v))$ are dense subgroups of $\R$, then so is
$\tau(I(u))\cap \tau(I(v))$.

\Pf Suppose the intersection consists of just $0$. We may find positive
real numbers $s \in \tau(I(u))$ and $t \in \tau(I(v))$ \st $s > \tau(u)$,
$t > \tau(v)$, and $0 < r:= s-t < \tau(u+v)$ (since the value group is
dense). By order unit goodness, there exists
$a $ \st $0 \leq a \leq u+v$ and $\tau(a) = r$. Riesz decomposition
entails $a = a_1 + a_2$ where $0 \leq a_1 \leq u$ and $0 \leq a_2 \leq v$.
Set $s' = \tau(a_1) \geq 0$ and $t' = \tau(a_2) \geq 0$. Then $s-t = s' +
t'$, so $s-s' = t+t'$. The intersection consisting of $0$ forces $s = s'$
and $t = -t'$; the latter forces $t= t' = 0$, a contradiction.

Now suppose the intersection is nonzero and not dense. Then it is cyclic,
so there exists $x \in \R$, which we may assume positive, \st $\tau(I(u))
\cap \tau(I(v)) = x\Z$. We may find $0 < s,t < x$ with $s \in \tau(I(u))$ and
$t \in \tau(I(v))$ \st $0 < r := s-t $. Find $a \leq u+v$ as above with $r
= \tau(a)$, similarly decompose $a = a_1 + a_2$, and define $s'$, $t'$ as
in the preceding paragraph.
We deduce $s - s'= t+ t'$; hence there exists an integer $m$ \st $s - s' =
mx = t+ t'$; as $t,t' \geq 0$, we have $m \geq 0$, but as $s <x$, we have
$m < 1$; hence $m = 0$. This forces $t = t' = 0$, again a contradiction.
\qed

\Lem Corollary \onefou. Let $G$ be a nearly divisible dimension group with a faithful trace $\tau$. Suppose that $I$ and $J$ are order ideals with their own order units \st each of $\tau|I$, $\tau|J$, and $\tau|(I+J)$ is order unit good. Then $\tau(I) \cap \tau(J)$ is a dense subgroup of $\R$.

\Pf Since $\tau$ is faithful, $\tau|I$ and $\tau|J$ are nonzero, and since every trace on an order ideal with order unit is nondiscrete (as the order ideals are approximately divisible by definition), it follows that $\tau(I)$ and $\tau(J) $ are dense. Now Lemma \onethr(b) applies. \qed

Let $(G,u)$ be a dimension group. Let $\Cal J$ be a collection of nonzero
order ideals each with their own order unit, \st every order ideal of $G$
with order unit can be expressed as a sum of order ideals from $\Cal J$
(such a sum can always be made finite, as the order ideal has an order
unit); then we say $\Cal J$ is a {\it generating set of order ideals of
$G$.}

The criteria in Proposition \onetwo\ for goodness can be reduced to that on a generating
set of order ideals. This will make the computations of section 4 much simpler.

\Lem Lemma \onefiv. Let $(G,u)$ be a nearly divisible dimension group, let $\Cal
J $ be a generating set of order ideals of $G$, and let $\tau$ be a
faithful trace of $G$. Sufficient for $\tau$ to be a good trace of $G$ is
that it satisfy
\item{(i)} for all $J \in \Cal J$, $\tau(J) = \tau(G)$ and
\item{(ii)} for all $J \in \Cal J$, $\tau|I$ is an order unit good trace
of $I$.

\Pf We can express a nonzero order ideal $I$ with order unit as $I = \sum
J_{\alpha}$ for some $ J_{\alpha} \in \Cal J$. Thus $\tau(I) = \sum
\tau(J_{\alpha}) = \tau(G)$.

Since $I$ has an order unit, the sum can be made finite; now we apply
induction (on the number of summands) to \onetwo(d); this verifies the second property in Proposition \oneone.
\qed

\comment
[This will be useful for $R_P$, where $\Cal J$ consists of $fR_P$ as $f$
varies over products of formal monomials.]
\endcomment

Verifying the various criteria for goodness and related properties is much simpler when the partially ordered abelian group is an ordered ring having $1$ as an order unit.

\Lem Lemma \onesix. Let $(R,1)$ be a (commutative) partially ordered commutative
ring with $1$ as order unit. If $R$ is approximately divisible, then it is
nearly divisible.

\Pf Approximate divisibility implies the existence of order units $u$ and
$v$ \st $1 = 2 u + 3v$; for any $r \in R^+\setminus \brcs{0}$, we thus
have $r = 2(ru) + 3(rv)$. From $1 \leq ku, kv $ for some positive integer
$k$, we deduce $r \leq k(ru), k(rv)$, verifying the definition of nearly
divisible.
\qed

The following is implicit in the proof of [BeH, Corollary 7.12].

\Lem Lemma \onesev. Let $(R,1)$ be a partially ordered (commutative) unperforated
ring with $1$ as order unit, that is approximately divisible. Let $\tau$
be a faithful pure trace. Then $\tau$ is order unit good iff for all
$\sigma \in \partial_e S(R,1) \setminus \brcs{\tau}$, $\sigma(\ker
\tau) \neq \brcs{0}$.

\Pf Since $1$ is an order unit of the partially ordered ring, $X:= \partial_e S(R,1)$ is
compact and consists precisely of the normalized multiplicative traces of
$R$; moreover, $\Aff S(R,1) = C(X,\R)$ with the affine representation
re-interpreted as $\tilde g(\phi) = \phi(g)$ for $\phi \in X$ (note the
use of $\ \tilde{\ }$ rather than $\ \hat{\ }$, to distinguish them). By
approximate divisibility, the image of $R$ is dense in $C(X,\R)$. If $A$
is any ideal of $R$, then its closure in $C(X,\R)$ is an ideal therein,
hence of the form
$\Ann(Y): = \Set{f \in C(X,\R)}{f|Y \equiv 0}$ for a unique compact
subset $Y$ of $X$.

Since $\tau$ is pure, it is multiplicative, and therefore $ \ker \tau$ is
an ideal of $R$ [{\it not\/} an order ideal, unless $\ker \tau = 0$, as
$\ker \tau \cap R^+ = \brcs{0}$ is the definition of faithfulness]. The
closure of the image of $\ker \tau $in $C(X,\R)$ can thus be written in
the form $\Ann (Y)$ for some compact subset $Y$.

If $\tau$ is order unit good, then $\Ann (Y)$ is $\Ann(\brcs{\tau})$
(corresponding to $\tau^{\perp}$ in $\Aff S(R,1)$), from which it follows
that $Y = \brcs{\tau}$. Hence if $\sigma \in X \setminus \brcs{\tau}$, there exists continuous $\Arrow f; X.[0,1]$ \st $f(\tau) = 0$ but
$f(\sigma) = 1$; then $f \in \Ann(\brcs{\tau})$, hence there exist $g_n \in
R$ \st $g_n \in \ker \tau$ and $\tilde g_n \to f$ uniformly. Applying
$\sigma$, there exists $n$ \st $\sigma(g_n) \neq 0$, so that $\sigma(\ker
\tau) \neq \brcs{0}$.

Conversely, suppose for every $\sigma \in X\setminus \brcs{\tau}$,
$\sigma(\ker \tau) \neq \brcs{0}$. Then $\sigma \not\in Y$; hence $Y =
\brcs{\tau}$, so that the closure of the image of $\ker \tau$ is
codimension one in $C(X,\R)$, hence equal to $\tau^{\perp} $ in $\Aff
S(G,u)$. Thus $\tau$ is order unit good.
\qed

\SecT 2 Tensor products

If $G$ and $H$ are partially ordered abelian groups, we may form the
tensor product (as $\Z$-modules) $G \otimes_{\Z} H$ (usually, we delete
the subscripted $\Z$); it is equipped with a cone which makes it into a
partially ordered group, $\Set{\sum g_i \otimes h_i}{g_i \in G^+ \text{
and } h_i \in H^+}$ [GH2, Proposition 2.1]. If both are dimension groups, then
so is $G\otimes H$, and if $u,v$ are respectively order units for $G,H$,
then $u\otimes v$ is an order unit for $G\otimes H$. If $\sigma$, $\tau$
are respective (normalized) traces on $(G,u)$ and $(H,v)$, then $\sigma
\otimes \tau$ (defined in the obvious way) is a (normalized) trace of $(G
\otimes H,u \otimes v)$.

A special case occurs when we  form the divisible hull of a dimension
group, $G \otimes \Q$, the rational vector space that $G$ generates. Then
$\tau$ extends to a trace $G \otimes \Q$ in the obvious way, denoted $\tau
\otimes 1_{\Q}$. In general, $\tau$ being order unit good or good implies
the corresponding property for $\tau \otimes 1_{\Q}$, but the converse
fails practically generically. As a special case, we [BeH] defined a trace
$\tau$ to be {\it ugly\/} if $\tau \otimes 1_{\Q}$ is good and $\ker \tau$
has {\it discrete\/} image in (the Banach space) $\Aff S(G,u)$. Ugly
traces exist in profusion.

In Akin's original context of measures on Cantor sets, he showed that
(what amounts to) the tensor product of good traces is good; in the
context of simple dimension groups or more generally for approximately
divisible dimension groups, the tensor product of order unit good traces
was shown to be order unit good. Here, we show a somewhat surprising
result for order unit goodness: if $(G,u)$ and $(H,v)$
are approximately divisible, and both $\sigma \otimes 1_{\Q}$ and $\tau
\otimes 1_{\Q}$ are order unit good on their respective groups, then
$\sigma \otimes \tau$ is order unit good (as a trace on $G \otimes H$).
This means that the tensor product has a stronger property (in general)
than its constituents. In particular, the tensor product of ugly traces is
at least order unit good.

Using the criterion of Proposition \oneone, we then obtain a corresponding
criterion for goodness of the tensor product ($G$ and $H$ are nearly
divisible, $\sigma \otimes 1_{\Q}$ and $\tau \otimes 1_{\Q}$ are good, and
a condition that guarantees the value groups on the order ideals is the
same as the full value group).

\Lem Proposition \twoone. Let $(G,u)$ and $(H,v)$ be approximately divisible dimension
groups with traces $\sigma$ and $\tau$ respectively. If each of $\sigma
\otimes 1_{\Q}$ and $\tau \otimes 1_{\Q}$ on $G \otimes \Q$ and $H \otimes \Q$ respectively is
order unit good, then the trace on $(G\otimes H, u\otimes v)$ given by
$\sigma \otimes \tau$, is order unit good.

If we only require that $ \sigma \otimes \tau \otimes 1_{\Q}$ (a trace on
$G \otimes H \otimes \Q$) be order unit good (in place of each of $\sigma
\otimes 1_{\Q}$ and $\tau \otimes 1_{\Q}$ being good), the conclusion may
still be true. In any event, I know of no counter-examples.

We require a number of elementary results about tensor products. Here the
tensors will be over one of the rings $\Z$, $\Q$, or $\R$; {\it
torsion-free   \ \paren{module}\/} means torsion-free abelian group when the
underlying ring is $\Z$, otherwise is just means vector space over the
relevant field.

\Lem Lemma \twotwo. Let $A$ and $B$ be torsion-free modules, and $A' \subset A$,
$B' \subset B$ submodules \st $A/A'$ and $B/B'$ are torsion-free.
\item{(a)} The kernel of the map $A \otimes B \to A \otimes B/B'$ is  $A
\otimes B'$.
\item{(b)} The kernel of the map $A \otimes B \to A/A' \otimes B/B'$ is $A
\otimes B' + A' \otimes B$.

\Pf (a) One inclusion is obvious. Because the quotient is torsion-free, $A
\otimes B/B'$ is torsion-free. We have an induced map $(A \otimes B)/(A
\otimes B') \to A\otimes B/B'$. If $z $ is in the kernel, find a nonzero
integer $n$ \st $nz$ has a representative in $A \otimes B$ of least length
(as $n$ varies over nonzero integers), say $nz = \sum a_i \otimes b_i + (A
\otimes B')$. Then $\brcs{a_i}$ is rationally linearly independent, hence
the image, $n\overline {z} $, yields, $0 = \sum a_i \otimes (b_i + B')$.
Since $B/B'$ is torsion-free, this easily implies all $b_i + B' = 0$
(tensor with $\Q$ if necessary, so we are working over a field, then use a
basis for $B'\Q$, extended to $B\Q$). [This proof works for all fields.]

(b) First, $A \otimes B/(A \otimes B')$ is naturally isomorphic to $A
\otimes B/B'$ by (a). Then another application of (a) with the order
reversed yields a natural isomorphism
$(A \otimes B/B')/(A' \otimes B/B') \iso A/A' \otimes B/B'$. Then the
kernel of the first map is $A \otimes B'$, and of the second is $A'
\otimes B/B'$, which pulls back to $A \otimes B' + A' \otimes B$.
\qed

\Pf (of Proposition \twoone) We will show that that the closure of the image of $\ker
\sigma \otimes \tau$ in $\Aff S(G \otimes H, u\otimes v)$ is $(\sigma
\otimes \tau)^{\vdash}$; by [BeH, Proposition 1.7], $\sigma \otimes \tau$ is order unit
good.

First, we identify $\Aff S(G,u)\otimes_{\R} \Aff S(H,v)$ with a subspace
of $\Aff S(G \otimes H, u\otimes v)$ in the obvious way. Standard results
(e.g., pure traces are pure tensors) yields that it is a dense subspace.

We note that $(\ker \sigma) \otimes H + G \otimes (\ker \tau) \subseteq
\ker \sigma \otimes \tau$. It easily follows that the closure of the image
of $(\ker \sigma) \otimes H$ contains everything in $y \otimes \Aff
S(H,v)$ (real tensors) where $y$ varies over the image of $\ker \sigma$
(in $\sigma^{\vdash} \subset \Aff S(G,u)$). For $y$ fixed, $y \otimes \Aff
S(H,v) $ is a real vector space, and this means that we can rewrite it as $y\R \otimes \Aff
S(H,v)$ (just approximate real multiples of $\hat v$ by elements of $\hat
H$, and transfer through the tensor product). Taking finite sums, we see
that the closure of the image of $\ker \sigma \otimes H$ includes the
closure of $\Im (\ker \sigma )\Q \otimes \Aff S(H,v)$.

Now $\sigma \otimes 1_{\Q}$ being order unit good implies $\ker \sigma
\otimes \Q $ has dense image in $\sigma^{\vdash}$ (in $\Aff S(G,u)$). If
$e$ is an element of $G \otimes \Q$, there exists a nonzero integer $m$
\st $me \in G$. If in addition, $\sigma \otimes 1_{\Q}(e) = 0$, then
$\sigma (me) = 0$; thus $\ker \sigma \otimes 1_{\Q} \subseteq (\ker \sigma
) \Q$ (the reverse inclusion is trivial, but never needed).

Thus the closure of the image of $(\ker \sigma) \otimes H $ contains $\Im
(\ker \sigma )\Q \otimes \Aff S(H,v)$, which in turn contains  the closure
of $\Im (\ker \sigma)\Q \otimes \Aff S(H,v)$, and thus includes
$\sigma^{\vdash} \otimes \Aff S(H,v)$.

Similarly, the closure of the image of $G \otimes \ker \tau$ contains
$\Aff S(G,u) \otimes \tau^{\vdash}$. Set $A = \Aff S(G,u)$, $A' =
\sigma^{\vdash}$, $B = \Aff S(H,v)$, and $B' = \tau^{\vdash}$; then each is
a Banach space, and $A/A'$ and $B/B'$ are both one-dimensional, and the
closure of the image of $\ker \sigma \otimes \tau$ contains $A' \otimes B
+ A \otimes B'$.

By (b) above, $A \otimes B/(A' \otimes B + A \otimes B')$ is
one-dimensional. Let $W = A' \otimes B + A \otimes B'$ and $Z = \Aff
S(G,u)\otimes \Aff S(H,v)$, so that $W$ is a codimension one subspace of
$Z$. It is an easy exercise to show that when we complete $Z$ to $\Aff S(G
\otimes H, u\otimes v)$, the closure, $\overline W$, is of at most
codimension one. (This is a general Banach space result; if $\overline W
\neq \overline Z$, then $W = \overline W \cap Z$ as $W$ is codimension one
in $\Z$; choose $z \in Z \setminus W$; the functional sending $z \mapsto 1$
and $W \mapsto 0$ is continuous (essentially the closed graph theorem), hence
extends to a bounded linear functional $p$ on $\overline W$; we may write
arbitary $y \in \overline{Z}$ as $\lim y_n$; then $y_n = p(y_n)z + (y_n -
p(y_n)z)$, and thus by continuity, $y = p(y)z +(y - p(y)z)$, and $y -
p(y)z$ is in $\overline W$; hence $z + \overline W = \overline Z$.

In particular, the closure of the image of $\ker \sigma \otimes \tau$ in
$\Aff S(G\otimes H, u \otimes v)$ is codimension one. As it is contained
in $(\sigma \otimes \tau)^{\vdash}$, which is proper, it follows that the
image of $\ker \sigma \otimes \tau$ is dense in $(\sigma \otimes
\tau)^{\vdash}$.
\qed

This explains  a phenomenon exemplified in [BeH, Example 9]. Let $G$
be a critical dimension group of rank $k+1$ (that is, a free rank $k+1$ abelian group densely embedded in $\R^d$, and equipped with the strict ordering therefrom [H4]). Then we say $G$ is {\it
basic\/} (as a critical group) if it is {\it order-isomorphic\/} to the
subgroup of $\R^k$ spanned by $\brcs{e_i; \sum \alpha_j e_j}$ where
$\brcs{e_i}$ is the standard basis and $\brcs{1, \alpha_1, \dots,
\alpha_k}$ is linearly independent over the rationals (this guarantees
density of the subgroup). Every critical group is {\it topologically
isomorphic\/} to a group of the latter form.

For basic critical groups, every pure trace is ugly, as is immediate from
the definitions.
Hence if $G_i$ are  basic critical groups (and there is more than one), their tensor
product (a simple dimension group) $\otimes G_i$ has all of its pure
traces good. In [BeH, Example 9],  an example was given of a basic critical group of
rank three, for which all pure traces on $G \otimes G$ are good. We also
asked whether the pure traces on $G \otimes G \otimes G$ are good, and now
we know that the answer is yes.

It is possible that among critical groups, basic ones are characterized
by all pure traces being ugly. There are lots of critical groups for which
all or some are bad, hence not ugly [BeH, section 2].

Now suppose that $(G,u)$ and $(H,v)$ are nearly divisible, and $\sigma$, $\tau$ are normalized traces
on $G$, $H$ respectively \st $\sigma \otimes 1_{\Q}$ and $\tau \otimes 1_{\Q} $ are both good. We expect
to obtain that $\sigma \otimes \tau$ is a good trace on $G \otimes H$.

\Lem Lemma \twothr. Let $(G,u)$ and $(H,v)$ be dimension groups with order unit.
\item{(a)} Then $G \otimes H$ is approximately divisible iff at least one
of $G$ or $H$ is;
\item{(b)} $G \otimes H$ is nearly divisible iff at least one of $G$ or
$H$ is.

\Pf (a) Suppose $G$ is approximately divisible. Every pure trace of
$(G\otimes H, u \otimes v)$ is of the form $\sigma \otimes \tau$ [GH2, Lemma 4.1], where
$\sigma$, $\tau$ are pure traces of $G$, $H$ respectively. Then $(\sigma
\otimes \tau) (G \otimes H) $ is $\sigma(G)\cdot \tau (H)$ (the set of sums
of terms of the form $\sigma(g)\cdot \tau(h)$); as $\sigma(G) $ is dense,
obviously so is $\sigma(G)\cdot \tau (H)$, so that $G \otimes H$ has no
discrete traces, and is thus approximately divisible. Obviously the same
argument applies if $H$ is approximately divisible.

If neither $G$ nor $H$ is approximately divisible, then there exists a
discrete trace $\sigma$ of $G$ and a discrete trace $\tau$ of $H$; as
these are normalized (at $u$, $v$ respectively), $\sigma(G) = (1/n) \Z$
and $\tau(H) = (1/m)\Z$ for some positive integers $m$ and $n$; then
$(\sigma \otimes \tau) (G \otimes H) = (1/mn)\Z$, which is discrete. Hence
$G \otimes H$ admits a discrete trace, thus is not approximately
divisible.

\noindent (b) Select $a = \sum g_i \otimes h_i \in (G\otimes H)^+$; from the defintion of the ordering on the tensor product,
we can assume each of $g_i$ and $h_i$ are positive in their respective groups. By definition, we can write $g_i = 2 a_i + 3 b_i$
where $0 \leq g_i \leq k a_i, kb_i$ for some positive integer $k$; since the sum is finite, we can take the same integer $k$ for
all $i$. Set $c_1 = \sum a_i \otimes h_i$ and $c_2 = \sum b_i \otimes h_i$. Then $a= 2c_1 + 3 c_2$; moreover, $\sum g_i \otimes h_i \leq k \sum a_i \otimes h_i$, that is, $a \leq k c_1$, and similarly $a \leq k c_2$.

If neither $G$ nor $H$ is nearly divisible, there exist an order ideal of
$G$ with its own order unit, $(I,w)$ together with a discrete trace (of
$I$) $\phi$, and an order ideal of $H$ with its own order unit, $(J, y)$
and a discrete trace on it, $\psi$. Then $\phi \otimes \psi$ is a discrete
trace (as above) of $I \otimes J$; this being an order ideal of $G \otimes
H$, the latter is not nearly divisible.
\qed

\Lem Lemma \twofou. Let $G$ and $H$ be nearly divisible, having  faithful traces $\sigma$ and $\tau$ respectively \st   $\sigma\otimes 1_{\Q}$ and $\tau\otimes 1_{\Q}$ are good as traces on $G \otimes \Q$, $H \otimes \Q$ respectively.
\item{(a)}Let $(I,w)$ be an order ideal of $G$ with its own order unit, and let $(J,y)$ be an order ideal of $H$ with its own order unit. Then $(\sigma \otimes \tau)|(I \otimes J)$ is order unit good.
\item{(b)} Suppose for each order ideal $I$ of $G$, $\sigma(I) = \sigma (G)$, and similarly, for each order ideal $J$ of $H$, we have $\tau(J) = \tau (H)$. Then for every nonzero order ideal $L$ of $G \otimes H$, we have $(\sigma \otimes \tau)(L) = (\sigma \otimes \tau)(G \otimes H)$
\item{(c)} Suppose the hypotheses of (b) apply. Let $(L, e)$ be an arbitary order ideal of $G \otimes H$ with its own order unit. Then $(\sigma \otimes \tau)|L$ is order unit good.

\Pf (a) Each of the restrictions of $\sigma\otimes 1_{\Q}$ and $\tau\otimes 1_{\Q}$ to $I \otimes \Q$ and $J \otimes \Q$
respectively is good, hence is order unit good, and thus $(\sigma \otimes \tau)|(I \otimes J)$ is an order unit good trace of $I \otimes J$.

\noindent (b) First, if $L = I \otimes J$ (where $I$ and $J$ are nonzero order ideals in $G$ and $H$ respectively), then $(\sigma \otimes \tau)(I \otimes J)$ is the subgroup of $\R$ generated by all terms of the form $\sigma (a)\cdot \tau(b)$, where $a \in I$ and $b \in J$, and $(\sigma \otimes \tau)(G \otimes H)$ has the same form, except $a$ and $b$ are allowed to vary over $G$ and $H$ respectively. Since for all $a \in G$, there exists $a' \in I$ \st $\sigma(a') = \sigma(a)$, and similarly for $\tau$, the two groups are equal.

If $e \in L^+$, then by definition of the tensor product ordering, we can write $e = \sum g_i \otimes h_i$. For an element $x $ in the positive cone of a dimension group, let $I(x)$ be the order ideal it generates; then it is easy to check (since sums of order ideals are again order ideals in a dimension group) that $L = I(e) = \sum I(g_i ) \otimes I(h_i)$; in particular, $L $ contains a tensor product of order ideals, so the previous paragraph applies.

\noindent (c) Every $e \in (G\otimes H)^+$ can be written in the form $e = \sum g_i \otimes h_i$ with $g_i \in G^+$ and $h_i \in H^+$. By (a), the restriction of $\sigma \otimes \tau$ to each of $I(g_i) \otimes I(h_i)$ is order unit good. Since $\sigma \otimes\tau(L) = (\sigma \otimes \tau)(G \otimes H)$, for any nonzero order ideal $L$ of $ G \otimes H$, we may apply \onetwo(e) (the intersection of the value groups is dense), so the restriction of $\sigma \otimes \tau$ to $L$ is order unit good.
\qed

\Lem Proposition \twofiv. Suppose that $(G,u,\sigma)$ and $(H, v, \tau)$ are nearly divisible dimension groups with faithful trace having the following properties:
\item{(i)} for all nonzero order ideals $I$ ($J$) of $G$ ($H$), $\sigma(I) = \sigma (G)$ ($\tau(J) = \tau(H)$);
\item{(ii)} each of $\sigma \otimes 1_{\Q}$ and $\tau \otimes 1_{\Q}$ is good on $G \otimes \Q$, $H \otimes \Q$ respectively.{\par}
\noindent Then $\sigma \otimes \tau$ is a good trace of $G \otimes H$.

\Pf Follows from \twothr, \twofou, and \oneone.
\qed

\SecT 3 Examples from xerox actions of tori on UHF algebras

We characterize the good faithful pure traces on the dimension groups arising
from xerox product type actions of tori on UHF C*-algebras. It turns out that
there is a surprising number-theoretic component.

Form the Laurent polynomial ring in $d$ variables over the integers, $\Z[x_i^{\pm1}]$, and let $\Z[x_i^{\pm1}]^+$
denote the set of those with only nonnegative coefficients. As in [H1, H2],
we adopt monomial notation, that is, for $w \in \Z^d$, define $x^w = x_1^{w(1)}\cdot x_2^{w(2)}\cdot \dots \cdot x_d^{w(d)}$.
For any $f \in \Z[x_i^{\pm1}]$, we denote the coefficient of $x^w$ in $f$ by $(f,x^w)$ (inner product notation, which is consistent
with the origins of the work), and we set $\Log f := \Set{w \in \Z^d}{(f,x^w) \neq 0}$.
Let $P = \sum a_w x^w \in \Z[x_i^{\pm1}]^+$ (where $a_w \in \Z^+$), and form the ring
$R_P = \Z[\brcs{x^w/P}_{ w\in \Log P}]$; equipped with the partial ordering generated added
and multiplicatively by  $\Set{x^w/P}{w \in \Log P}$, this is a dimension group and an ordered ring
with $1$ as order unit, and many more properties. We may also form $\Z[x_i^{\pm1}, 1/P]$ (a subring of
the field of fractions of the Laurent polynomial ring. It also has a partial ordering given
by $\Set{f/P^k}{\exists N \text{ \st } P^N f \text{ has no negative coefficients}}$. The restriction of this
to $R_P$ yields the original ordering.
\def\Ad{\text{Ad}\,}

This arose from the following construction. Let $n = P(1,1,1,\dots,1)$, and form $\Cal A = \otimes M_n \C$ (the
UHF C*-algebra). The Laurent polynomial $P$ is the character of an $n$-dimensional representation of the torus
$\T^d$, say given by $z \mapsto \diag(z^w)$ (one for each $w$ that appears in $P$, with repetitions as indicated by the
multiplicities, that is, the coefficients. This yields a map $\Arrow \pi; \T^d . M_n\C$ with nonzero entries along the diagonal.
Form $\Arrow \phi:=\otimes \Ad \pi;\T^d. \text{Aut }\Cal A$, and the corresponding fixed point subrings, $\Cal A^{\phi(\T^d)}$, and
$\Cal A \times_{\phi} \T^d$, the latter the C*-crossed product. Then $(\text{K}_0 (\Cal A^{\phi(\T^d)}),[1])$ is naturally ordered ring isomorphic to
$R_P$ and $\text{K}_0 (\Cal A\times_{\phi} \T^d)$ similarly isomorphic to the ordered ring $\Z[x_i^{\pm1}, 1/P]$. This will play a role in
what follows.

Renault [R] determined the positive cone and analyzed (inter alea) the structure of $R_P$ when $P = 1+x$. That was in 1980; people are still obliviously reproving his and other results (concerning Pascal's triangle Bratteli diagrams) 30+ years later!

We normally assume that $P$ is projectively faithful, that is, $\Log P - \Log P$ generates (as an abelian group) the standard
copy of $\Z^d$ in $\R^d$ (we can reduce to this case anyway). This has the effect that whenever $v \in \Log P^k \cap \Int \cvx \Log P^k$ for some positive
integer $k$,
$x^v /P^k$ belongs to $R_P$ and $R_P[(x^v/P^k)^{-1}] = \Z[x_i^{\pm1}, 1/P]$, i.e., the larger ring is obtained by inverting $x^v/P^k$.

We call an element of the form $x^w/P$ with $w \in \Log P$ a {\it formal monomial\/} in $R_P$. (It can happen that $x^w/P \in R_P$ even if $w \not\in \Log P$---e.g., if $w + \Log P^k \subseteq \Log P^{k+1}$ for some $k$. These aren't significant in what follows.)

In addition to the obvious facts about $R_P$ (it is a commutative,
finitely generatedhence noetheriandomain), the following results are
known [H1, H2]:
\item{} $R_P = \Set {g/P^k}{g \in \Z[x], \ \Log g \subset \Log P^k}$,
$R_P$ is a partially ordered ring with $1$ as an order unit,
and it is a dimension group [H1, section I];
\item{} all sums and finite intersections of order ideals are order
ideals are order ideals (this is true for all dimension groups) [G];
\item{} products of order ideals are order ideals (this is not generally
true for commutative partially ordered domains having $1$ as an order unit
and being dimension groups) [H1];
\item{} every order ideal is an order ideal (true in every partially  ordered commutative ring in which $1$ is an order unit) [H1, Proposition I.2];
\item{} if $f$ is a formal monomial, then $fR_P$ (the ideal generated by
$f$) is an order ideal [H2; Proposition II.2A];
\item{} every order ideal is the finite sum of ideals, $\sum f_i R_P$
where $f_i$ are formal monomials, and all such sums are order ideals [H2, p\,19];
\item{} if $f$ is a formal monomial and $a \in R_P$, then $f a \in R_P^+$
implies $a \in R_P^+$ (follows from the definitions); the conclusion is also true if we replace {\it formal monomial\/} by
{\it order unit,} a result that is very special for $R_P$ [H2; Proposition II.5];
\item{} the pure traces are exactly the multiplicative ones (true for any
partially ordered ring with $1$ as an order unit); the pure faithful
traces  are exactly those of the form
$\tau_r(g/P^k) = g(r)/P^k(r)$ where $r = (r_i)$ is a strictly positive
$d$-tuple in $\R^d$, and these extend in the obvious way to positive
homomorphisms $\Arrow \tau_r; \Z[x_i^{\pm1};1/P].\R$ (warning: although the ring
$\Z[x_i^{\pm1};1/P]$ is partially ordered, $1$ is not an order unit for it) [H1, Theorem III.3];
\item{} the weighted moment map/Legendre transform corresponding to $P$ implements a
homeomorphism $\partial_e S(R_P,1) \to \cvx \Log P$ (the latter is the
{\it Newton polytope\/} of $P$) sending the faithful pure traces onto the
interior; unexpectedly, the set of {\it pure\/} traces admits a type of
convex structure; in particular, the faces correspond to traces that
factor through quotients in a particularly nice way [H2, Theorem IV.1];
\item{} In general, $R_P$ is not a pure polynomial ring; only rarely does it have
unique factorization [H2, Appendix~ A, Theorem A.8A].\vskip 4pt

\noindent Now let us consider the following  property of a faithful pure
trace $\tau \equiv \tau_r$:
\item{(1)} for every nonzero order ideal $I$, $\tau_r (I) = \tau_r(R_P)$.

\noindent By
Proposition \oneone, this is one of the two necessary conditions for $\tau_r$ to be a good trace.

Here $r = (r_i) \in (\R^d)^{++}$ as described above. First we note that $\brcs{fR_P}$ (as $f$ varies over all products of formal monomials) is a generating set of order ideals with order unit (they are given as ring ideals, but in fact are order ideals by the properties above, and every order ideal is a finite sum of these).
Necessary and
sufficient for (1) to hold is simply that it hold for all ideals of
the form $I_w= (x^w/P)R_P$ (where $w \in \Log P$, a finite set). To see
this, note that $\tau_r(I_w) = (r^w/P(r)) \tau(R_P)$, hence $\tau_r (I_w)
= \tau_r (R_P)$ iff $P(r)/r^w \in \tau_r(R_P)$; thus if this holds for all
$w\in \Log P$, then each of $P(r)/r^w$ belong to $\tau(R_P)$, hence all
their products do; this means that for every formal monomial $f$,
$1/\tau_r(f)$ belongs to $\tau_r (R_P)$, hence $\tau_r (fR_P) =
\tau(R_P)$.

The upshot of this is that $\tau_r$ satisfies (1) if and only if for all
$w\in \Log P$, $P(r)/r^w \in \tau_r (R_P)$. The latter is simply
$\Z[r^w/P(r)]_{w \in \Log P}$. So we deduce

\Lem Lemma \fouone. For $r \in (\R^d)^{++}$, $\tau_r$ satisfies (1) iff for all $v
\in \Log P$,
$P(r)/r^v \in \Z[x^w/P]_{w \in \Log P}$.

This is a fairly drastic condition, even when $d = 1$ and $P = 1+x$ or $2
+ 3x$.

For $r \in (\R^d)^{++}$ and $P \in \Z[r_i]^+$, let $R_r = \Z[\brcs{r^w/P(r)}_{w \in \Log P}]$; this is exactly $\tau_r (R_P)$, and is a finitely generated unital subring of $\R$. The next lemma says that $r$ satisfies (1) iff when we extend $\tau_r$ all the way up to $\Z[x_1^{\pm1}, \dots, x_d^{\pm1}, P^{-1}]$, the image of $\tau_r$ does not increase---something we should have expected, in terms of the original definition.

\Lem Lemma \foutwo. Let $r  = (r_i)\in (\R^d)^{++}$ and $P \in \Z[r_i]^+$ be projectively faithful. Then $r$ satisfies (1) iff $R_r = \Z[r_i^{\pm 1}; P(r)^{-1}]$.

\Pf We may construct $R_P$ by beginning with $\Z[x_i^{\pm1}]$ (the Laurent polynomial ring) instead of $\Z[x_i]$; this is in fact how it was originally constructed in [H1, H2]. By replacing $P$ by $x^v P^t$ for some $v \in \Z^d$ and positive integer $t$ (this has no effect on $R_P$, up to order isomorphism), we can arrange that $\pmb 0$ is in the interior of $\cvx \Log P$ and in $\Log P$. Then $1/P \in R_P$ and we may invert $1/P$, creating
$R_P[P] = \Z[x_i^{\pm 1}; P^{-1}]$ [H2]. Let $I = (1/P)R_P$; this is an order ideal ([H2, p\,19]), and $\Z[x_i^{\pm 1}; P^{-1}] = \cup_{j \in \Z^+} P^{j}R_P$.

If $r$ satisfies (1) \wrt $P$, then applying it to $I$, we obtain $\tau_r(I) = \tau_r (1/P)\tau_r (R) = (1/P(r))\tau_r(R) = (1/P(r))R_r$; by hypothesis, this is $R_r$, so that $P(r) \in R_r$. Thus $\tau_r (P^j R_P) = P^j(r)R_r \subset R_r$. Taking the union, we obtain $\tau_r (\Z[x_i^{\pm 1}; P^{-1}]) \subseteq R_r$, and the reverse inclusion is trivial.

Conversely, suppose $R_r = \tau_r (\Z[x_i^{\pm 1}; P^{-1}])$. Then $\tau_r (x_i^{\pm 1}) = r_i^{\pm 1}$ and $\tau_r(P^{\pm 1}) = P^{\pm1}(r)$ belong to $R_r$ and are invertible therein. Hence if $f$ is any formal monomial, $\tau_r (f)$ is a product of terms of the form $r^w/P(r)$, hence is invertible in $R_r$. Thus if $I$ is an order ideal, it contains a formal monomial, and $\tau_r (I)$ contains an invertible element in $R_r$, and so $\tau_r(I) = R_r = \tau_r (R_P)$. Thus $r$ satisfies (1).
\qed

In other words, (1) holds iff the range of evaluation at $r$ on $R_P$ is the same as the range of the evaluation on the much larger ring $\Z[x_i^{\pm 1}, 1/P]$.

Now we consider what (1) means in the special case of $d= 1$.

Let $A$ be a unital subring of $\C$, the complexes. A complex number $r$ is {\it integral over $A$} (or {\it $r$ is an $A$-algebraic integer\/}) if it satisfies a monic polynomial with coefficients from $A$; equivalently, $r \in A[r^{-1}]$. The number $r$ is an {$A$-algebraic unit\/} if it satisfies a monic polynomial with coefficients from $A$ whose constant term is invertible in $A$; equivalently, $A[r] = A[r^{-1}]$. If $A = \Z$, we just write {\it integral\/} (adjective) or {\it algebraic integer\/} (noun). If $A = \Q$, these notions coincide, and we just say $r$ is {\it algebraic.} The {\it degree\/} of an integral or algebraic element is the degree of its minimal polynomial (over $A$).

\Lem Lemma \fouthr. Let $P$ be a projectively faithful element of $\Z[x]^+$ with smallest and largest degree coefficients $a_0$ and $a_k$ respectively. If $r \in \R^{++}$ satisfies (1) \wrt $P$, then there exist nonnegative integers $s$ and $t$ \st $a_0^s/r$ and $a_k^t r$ are integral.

\Pf Write $P = a_0 + \sum_{0< i < k} a_i x^i + a_k x^k$ where $a_i$ are
nonnegative integers (some can be zero, but we still need $\gcd \(
\Set{i}{a_i \neq 0} \cup \brcs{k}\) = 1$). From $P(r) \in
\Z[\brcs{r^j/P(r)}_{j\in \Log P}]$, we deduce an equation of the form
$P(r)^{m+1} = p(r)$ where $p \in \Z[x]$ and $\deg p \leq \deg P^m = km$.
The leading term of this expression is $a_k^{m+1} r^{(m+1)k}$, and so $r$ satisfies
a monic polynomial with coefficients from $A = \Z[a_k^{-1}]$. It follows that $a_k^t r$ is integral for all sufficiently large $s$.

Replacing $P$ by its {\it reversal\/} (also called {\it reciprocal\/}) $\tilde P$ (defined by $\tilde P (x) = P(x^{-1})x^k$), and redoing the process yields the other form, that $a_0^s/r$ is integral.
\qed

The following is true  if we weaken the hypotheses on $P$ to be projectively faithful (instead of requiring all the intermediate coefficients to be strictly positive). The modifications to the proof will muddy an already-complicated but elementary argument; so we just outline it afterwards. We can replace $P$ by any power of itself, without changing anything, so the no gaps  condition is just that the second largest and second smallest terms have nonzero coefficients.

\Lem Proposition \foufou. Let $r \in \R^{++}$ and $P \in \Z[x]^+$ be $\sum_{i=0}^{k} a_i x^i$ where all $a_i \neq 0$. Let $a_0$ and $a_k$ be the coefficient of the least and greatest degree terms in $P$. Let $R_r = \Z[\brcs{r^i/P}_{i \in \Log P}]$. Then
the following are equivalent
\item{(i)} $r$ satisfies (1) \wrt $P$
\item{(ii)} there exist nonnegative integers $s$ and $t$ \st both $a_k^s r$ and $a_0^t /r$ are algebraic integers
\item{(iii)} $R_r = \Z[r^{\pm 1}, P(r)^{\pm 1}]$
\item{(iv)} for all $j \in \Log P$, $P(r)/r^j \in R_r$.

\Pf (ii) implies (iv). Without loss of generality, we may assume $P = a_0 + \sum_{0 < i < k } a_i x^i + a_k x^k$.

If $c$ is an algebraic integer, then $\Z[c]$ is free on the $\Z$-basis $\brcs{1, c, c^2, \dots, c^{e-1}}$ where $e$ is the degree of $c$ (this is an alternative definition of integrality); in particular, for every positive integer $u$, we can write $c^u = \sum_{i=0}^{e-1} b_i c^i$, in other words, there exists a polynomial $p \in \Z[x]$ of degree at most $e-1$ \st $c^u = p(c)$.

Apply this to $c = a_k^s r$; for each positive integer $u$, we can write $(a_k^s r)^u = p_u (a_k^s r) = q_u (r)$ where $\deg q_u \leq e-1$. Multiplying this by $r^{u(s-1)}$, we obtain $(a_k r)^{us} = r^{u(s-1)}q_u$; setting $Q_u = x^{u(s-1)}q_u$, we have $(a_k r)^{us} = Q_u (r)$ where $Q_u \in \Z[x]$ and $\deg Q_u = u(s-1) + \deg q_u \leq u (s-1) + e-1$.
Hence (multiplying by an additional $r^j$), we have for every $j = 0,1,2, \dots$, $Q_{u,j}\in \Z[x]$ \st $\deg Q_{u,j} = u(s-1) + j$ and $(a_k r)^{us + j} = Q_{u,j}(r)$. We will subsequently choose $u$ to be fairly large.

Now let $N$ be a (large) positive integer, and consider the $k$ leading coefficients of $P^N$, that is, the coefficients of the terms $x^{kN}, x^{kN-1}, x^{kN-2}, \dots, x^{kN-k+1}$. They are respectively divisible by $a_k^N, a_k^{N-1}, \dots, a_k^{N-k+1}$ (as is trivially easy to see). Hence we may find integers $b_i$ (with $b_0 = 1$) \st
$$
P^N - \sum_{i=0}^{k-1} (a_k x)^{N-i} x^{N(k-1)}b_i:= G
$$
is a polynomial of degree at most $Nk - k$. Assume (as we may) that $N - k= us$ for some integer $u$. Replace each $(a_k x^{N-i})$ by $Q_{u, k-i}$; this has no effect on the value at $r$. Setting $H = \sum_{i=0}^{k-1} b_i Q_{u,k-i}x^{N(k-1)}$, we have
$P^N(r) = (G+H)(r)$. Then
$$\eqalign{
\deg (G + H) & \leq \max \brcs{\deg G, \deg H} \cr
& \leq \max\brcs{Nk-k, \max_i\brcs{\deg Q_{u,k-i} + Nk-N}}\cr
& \leq \max\brcs{Nk-k, u(s-1) + e-1 + Nk-N } = \max\brcs{Nk-k, Nk-N + e-1 + N-k - u}\cr
& \leq\max \brcs{Nk-k, Nk - k - u + e-1}.\cr
}$$
We can choose $u \geq e-1$ at the outset, and so guarantee that $\deg (G + H) \leq Nk-k$.
Thus $P(r) = (G+H)(r)/P^{N-1}(r)$. For every $0 \leq i \leq k$, $r^i/P(r) \in R_r$, and since$\deg (G+H) \leq Nk-k = \deg P^{N-1}$, we obtain $P(r) \in R_r$.

Now form the reversal of $P$, given by $\tilde P (x) = P(x^{-1})x^k$; this reverses the roles of $a_k$ and $a_0$, and the same process (using $a_0^t/r$ being integral) yields after translating back, $P(r)/r^k \in R_r$. From $P(r) \in R_r$, we obtain $r^i = (r^i/P(r))\cdot P(r) \in R_r$ for $i \in \Log P$, and thus for all $i \geq 0$. Since $P(r)/r^k \in R_r$, we deduce $r^{-k} \in R_r$, hence $r^{-j} \in R_r$ for all $j \geq 0$; hence $P(r)/r^j \in R_r$.

Now (i) implies (ii) was done in the previous lemma, and the equivalence of (i), (iii), and (iv) follows from the general results preceding this.
\qed

To prove the result when $P$ is projectively faithful, we can still write $P = a_0 + \sum_{1 \leq i \leq k-1} a_i x^i + a_k x^k$, only this time $\gcd{\Set{i}{a_i \neq 0}} = 1$ (equivalent after translation to projective faithfulness). Then it is elementary, and presumably well-known, that there exists $M$ \st for all $N$, $(P^N,x^i) \neq 0$ if $M < i < kN-M$. Now in the construction above, make sure that when the multiplications by powers of $r$ take place, that the exponent lands in the interval where all the coefficients are guaranteed nonzero (we are of course free to take arbitrary large powers of $P$).

A strange consequence is that when the hypotheses on $P$ are satisfied, the set of $r$ \st $\tau_r$ satisfies (1) is closed under multiplication; this follows immediately from (ii), but not obviously from any of the other equivalent properties.

This does not appear to extend to more than one variable. For example, if $P = 2 + 3x + 5y$, and we restrict to $r = (m,n)$ with positive integer coordinates, it is tedious but routine to see that $\tau_r$ satisfies (1) \wrt $P$ iff
for all primes $p$ and $q$,
$$
p|m \implies p|(2 + 5n) \qquad \text{and} \qquad q|n \implies q|(2 + 3m).
$$
For example, $(7,1), (3,11), (2^i,2^j)$ (where both $i,j> 0$) satisfy these conditions, but $(14,2)$ does not. Of course, there may be another, more appropriate, notion of multiplication \wrt which the set is closed.

Another general property concerns approximate divisibility. Let $K = \cvx
\Log P$; this is a compact convex polytope. Let $e \in K$ be an extreme
point (we do not use the usual term, {\it vertex,} because this might be
confused with lattice point); then $v \in \Log P$, and there is a pure
trace associated with $v$, $\sigma^v$, given by
$\sigma^v (g/P^k) = (g,x^{kv})/(P,x^v)^k$ (this can also be obtained as
the limit along a path of $\tau_r$, via l'Hpital's rule, as in [H1, section III, especially just before III.3]).

Since every order ideal of $R_P$ is of the form $\sum f_i R_P$ (finite
sum), if we assume that $R_P$ is approximately divisible, then $R_P$ is
nearly divisible. Thus every order ideal has its own order unit and is
approximately divisible. If $\tau$ is faithful, then $\tau(I \cap J) \neq
0$ (no finite intersections of order ideals can be zero since they are
also ideals in a domain), and $I \cap J$ is itself approximately
divisible, hence $\tau(I\cap J) $ is dense in $\R$. Thus for any faithful
trace that is order unit good for $R_P$, its restriction to any nonzero
order ideal is also order unit good.

Thus we have the following.

\Lem Lemma \foufiv. The ordered ring $R_P$ is approximately divisible iff for all
extreme points $v$ of $K = \cvx \Log P$, $(P,x^v) > 1$.

\Lem Lemma \fousix. Let $P = \sum \lambda_w x^w \in \Z[x^{\pm1}_i]^+$ with $(P,x^v) > 1$
for all extreme points of $K= \cvx \Log P$.
\item{(a)} Then $R_P$ is nearly divisible
\item{(b)} If $\tau$ is a faithful trace that is order unit good for
$R_P$, then its restriction to any nonzero ideal is order unit good for
that ideal.

If we replace $R_P$ by $S_P:=R_P \otimes \Q = \Q[x^w/P]$, then it is
divisible, which is of course stronger than nearly divisible, so that (a)
holds automatically (without the hypothesis on the coefficients at extreme
points), and (b) also holds by the same arguments.

\Lem Proposition \fousev. Let $r = (r_i) \in (\R^d)^{++}$, and let $P \in
\Z[x^{\pm1}_i]^+$ be projectively faithful.
\item{(a)} the pure trace $\tau_r$ on $R_P$ is good iff
\itemitem{(i)} $\tau_r$ is order unit good for $R_P$ and
\itemitem{(ii)} for all $v \in \Log P$, $P(r)/r^v \in \Z[r^w/P(r)]_{w \in
\Log P}$.
\item{(b)} the pure trace $\tau_r$ on $S_P$ is good iff
\itemitem{(i)} $\tau_r$ is order unit good for $R_P$.

\Rmk Note the absence of (ii) from (b), and the appearance of $R_P$ in
(bi). It is known (along the same lines as in [BeH, Proposition 5.10]), that if $\tau_r$ is
order unit good (for either coefficient ring), then each $r_i$ is
algebraic. Since $\Q[r_1, \dots, r_d]$ is thus a field, (ii) is redundant
in (b).

\Pf  We  show that if $\tau_r$ is order
unit good (which means that the closure of the image of $\ker \tau_r$ in
$\Aff S(R,1)$ is exactly $\tau_r^{\perp} = \Set{h \in \Aff
S(R,1)}{h(\tau_r) = 0}$), then its restriction to any order ideal is also
order unit good. It suffices to do this for $I = fR_P$ where $f$ is a formal monomial.

The map $R_P \to
fR_P$ given by $r \mapsto fr$ is an order-isomorphism of $R_P$ modules
(this of course uses the the fact that $fr \geq 0$ in $R_P$ entails $r
\geq 0$). Using $f$ as an order unit for $I$, the map on traces $\tau
\mapsto \tau/\tau(f)$ (restricted to those $\tau$ \st $\tau(f) \neq 0$
sends $\tau_r \to \tau_r/ \tau_r(f) = \tau'$, and $\ker \tau'= \ker \tau_r
\cap fR_P = f\cdot \ker \tau_r$ (since $f(r) \neq 0$). The map between
$R_P$ modules induces an affine homeomorphism between $S(R_P,1)$ and
$S(I,f)$, sending $\tau_r$ to $\tau'$, and it easily follows that $\tau'$
is order unit good. But $\tau'$ is just the normalization of $\tau|I$,
hence the latter is order unit good.

The rest follows from the preceding results.
\qed

In one variable, we can show that $\tau_r$ is order unit good iff none of
the algebraic conjugates of $r$ (except itself) are positive real. In
more than one variable, the situation is far more complicated, and there
is no decisive theorem (yet).

\Ex Example. Let $d = 1$ and $P = 1+x$; then we can rewrite $R_P = \Z[1/P,
x/P] = \Z[1-X,X]$ where $X = x/(1+x) $, and the positive cone translates
to $\langle X,1-X\rangle$. This goes back to Renault [R]. The translation
however, obscures some of the features, as we will see.

First, $R_P$ has two discrete pure traces, $\tau_0 = \sigma^0$ and
$\tau_{\infty} = \sigma^1$ ($0$ and $1$ are the extreme points of the
convex set $\cvx \Log P = [0,1]$), so is not approximately divisible.
However, it is interesting to calculate the condition that $\tau_r (I) =
\tau_r(R_P)$ for all nonzero order ideals.

By \fousev\ above, this amounts to $1+r, 1+1/r \in \Z[1/(1+r), r/(1+r)]$; as
$r/(1+r) = 1 - 1/(1+r)$, the condition (1) is equivalent to $1 + r^{\pm 1}
\in \Z[1/(1+r)]$. Now for a real number $s$, the condition $s \in \Z[1/s]$
is equivalent to $s$ be an algebraic integer (that is, satisfies a monic
integer polynomial). Hence we infer that if (1) holds for $\tau_r$, then
$r$ has to be an algebraic unit (that is, not only is its minimal
polynomial over the integers monic, but the constant term must be $\pm1$
as well). Conversely, if $r$ is an algebraic unit, then the desired
membership property holds.

We conclude that $\tau_r$ satisfies (1) iff $r$ is an algebraic unit.

In particular, if $r$ is an integer, then $\tau_r$ satisfies (1) iff $r =
1$ (we are restricting ourselves to actual traces, hence excluding
negative values for $r$).

The translation, $X = x/(1+x)$ converts $r$ to $r/(1+r)$; then of course
$\tau(X) $ is a fractional linear transformation of an algebraic unit, but
this characterization is not as pleasant as the pre-translation version.
\qed

Let $V \subset \C^d$. For $A$ a subring of $\C$, define $I_A(V)$ to be
the ideal in the polynomial ring $A[x_1,\dots, x_d]$ consisting of
polynomials that vanish at all points of $V$. Given an ideal $I$ of
$A[x_1, \dots, x_d]$, define $Z_A(I) $ to be the common zero set (in
$\C^d$) of all elements of $I$. The {\it variety generated by $V$ over
$A$} is simply $Z_A I_A (V)$. If $A = \Z$, we drop the subscript.

We say $r = (r_i) \in (\R^d)^{++}$ is {\it really isolated\/} if
$ZI(\brcs{r}) \cap (\R^{d})^{++} = \brcs{r}$. For example, if $d = 1$,
then $r$ is really isolated if $r$ is algebraic and all algebraic
conjugates of $r$ other than $r$ itself are {\it not\/} positive real. In
general, $r$ is really isolated means that the slice of the variety
generated by $r$ (or more simply, the Zariski closure of $\brcs{r}$) by
the positive orthant contains only $r$.

The argument in [BeH, 5.10] shows that if $r$ is really isolated (or more generally, $\brcs{r}$ is an isolated point in $(\R^d)^{++}\cap ZI(\brcs{r})$, then
all of its coordinates are algebraic (there is an assumption in [op\,cit]
concerning interior points which is automatic here). We remind the reader that we have assumed that $P$ is projectively faithful, which implies in particular, that its Newton polytope contains a $d$-ball.

The condition that $r $ be really isolated appears in the examples in [BeH, Examples 5 \& 10], for which the relevant dimension groups are remotely related to the ones appearing here.

\Lem Proposition \threig. Suppose $R_P$ is approximately divisible, and $\tau$ is
a pure faithful trace. Then
\item{(a)} $\tau$ is an order unit good trace of $R_P$ iff $\tau =
\tau_r$ where $r \in (\R^d)^{++}$ is really isolated.
\item{(b)} $\tau_r$ is a good trace of $R_P$ iff $r$ is really isolated
and for all $v \in \Log P$, $P(r)/r^v \in \Z[\brcs{r^w/P(r)}_{w \in \Log
P}]$.
\item{(c)} $\tau_r$ is a good trace of $R_P \otimes \Q$ iff $r$ is really
isolated.

\Pf Every pure faithful trace of $R_P$ is of the form $\tau_r$ for
(unique) $r $ in the positive orthant.

If $r$ is not really isolated, then there exists $r' \in ( \R^d)^{++}$ \st
every polynomial that vanishes at $r$ also vanishes at $r'$. Suppose $a:=
g/P^k \in R_P$; we may assume $\Log g \subseteq \Log P^k$. If $\tau_r (a)
= 0$, then $g(r) = 0$, hence $g(r') = 0$, whence $\tau_{r'}(a) = 0$; thus
with $\sigma = \tau_{r'}$, we have $\sigma \in \partial_e S(R,1) \setminus
\brcs{\tau_r}$ \st $\sigma|\ker \tau_r \equiv 0$. Hence $\tau_r$ is not
order unit good. The same of course applies with $R_P \otimes \Q$ in place
of $R_P$.

Conversely, suppose that $r$ is really isolated, but there exists $ \sigma
\in \partial_e S(R,1) \setminus \brcs{\tau_r}$ \st $\sigma|\ker \tau_r =
0$. Then $\sigma$ cannot be faithful (as otherwise, $\sigma= \tau_{r'}$
for some $r' \in (\R^d)^{++}$, and $r' \in ZI(\brcs{r})$). Consider $S =
R_P \otimes \Q$, and let $T_r$, $\Sigma$ be the extension to $S$ of
$\tau_r$ and $\Sigma$ (both extend, since the ranges are torsion-free
abelian groups). Then $T_r (S) = \Q[r^w/P(r)]$, which is a field (since
the coordinates are algebraic, so are all the $r^w/P(r)$). Then $\ker T_r
$ is a field, so $\ker T_r$ is a maximal ideal.
Also, $\ker T_r \cap R_P = \ker \tau_r$ and $\ker \Sigma \cap R_P = \sigma$.
If $\ker \tau_r \subseteq \ker \sigma$, then $\ker T_r \subset \ker
\Sigma$, but maximality of $\ker T_r$ implies $\ker T_r = \ker \Sigma$,
and thus $\ker \tau_r = \ker \sigma$. However, since $\sigma$ is not
faithful, $\ker \sigma $ contains a positive nonzero element of $R_P$,
whereas $\ker \tau_r$ does not, a contradiction.

Hence if $r$ is really isolated, then $\sigma \in \partial_e
S(R_P,1)\setminus \brcs{\tau_r}$ implies $\sigma(\ker \tau_r) \neq 0$, and
by Lemma \onesev\ above, this implies $\tau_r$ is order unit good. The same of course
applies to $T_r$ as a trace on $S_P$. This yields (a), and contributes to
(c).

Part (b) now follows from preceding results in this section.

Part (c) comes from $\Q[r^w/P(r)]$ being a field (which in turn arises
because the coordinates of $r$ are algebraic), so that condition (1) is
automatic.
\qed

A particular consequence is that the set of good pure faithful traces of
$S_P= R_P \otimes \Q$ is the same for all choices (with $d$ fixed) of
faithfully projective $P \in \Z[x_i]^+$ (or $P \in \Q[x_i]^+$), whereas
for $R_P$, there is dependence on $P$.

When $d=1$, the conditions for $\tau_r$ to be good are precisely that no distinct algebraic conjugate of $r$ be positive and the integrality condition, (ii), of Proposition \foufou.

\Ex Example. Let $d=1$ and $P = 2 + 3x$. By Proposition \foufou, the positive real number $r$ satisfies (1) iff there exists $s$ \st both $2^s/ r$ and $3^s r$ are integral. Let $K = \Q(r)$, and $\Z_K$ the ring of integers in $K$. The fractional ideal $r\Z_K$ factors as $\prod \Cal P_i/\prod \Cal Q_j$ (where $\Cal P_i$ and $\Cal Q_j$ are prime ideals in $\Z_K$, and we allow repetitions; the products might also be over the empty set). The intersections $\Cal P_i \cap \Z$ and $\Cal Q_j \cap \Z$ determine primes in $\Z$, denoted respectively $p_i$ and $q_j$. Then (1) is equivalent to $p_i = 2$ and $q_j = 3$ for all $i$ and $j$.

Hence $\tau_r$ is good for $R_P$ iff no non-identity algebraic conjugate is positive and the prime factorization of the fractional ideal $r\Z_K$ consists of primes sitting over $2$ in the numerator and over $3$ in the denominator.
\qed

In this section, we have restricted ourselves to pure {\it faithful\/} traces; this is a technical convenience. By the comment after Proposition \oneone, we can factor out the maximal order ideal contained in the kernel of a trace, and in the case that the dimension group is $R_P$, these correspond to quotients corresponding to faces of the Newton polytope ([H1, section VII]). This amounts to a reduction to a lower dimensional lattice and vector space, that is, a polynomial in fewer variables.

There are related naturally occurring classes of dimension groups whose pure traces can be similarly analyzed. For example, for the matrix-valued random walks appearing in [H5], in  non-degenerate cases, the pure faithful traces are similarly parameterized by the positive orthant (the non-faithful traces are generically terrible, but can be analyzed in reasonable cases). An amusing example appears in [P], where very specific local limit asymptotics were used to derive the one-parameter family (indexed by the unit interval) of pure traces. In fact, that random walk can be represented as $M= \(\smallmatrix 1 +x & x \\ 1 & 0 \\ \endmatrix \)$, and in this very simple case, via [H5], we can write down the pure traces parameterized by $[0, \infty]$ (the endpoints corresponding to the two non-faithful pure traces) via the large eigenvalue function. Alternatively, it is elementary that $(1+x)\hat M^{-1}$ is an order unit in $E_b (G_M)$, so on setting $P = 1+x$, we can view $M/P$ as a matrix with entries in $R_P$ without changing the pure trace space. This yields a parameterization of the pure traces by those of $R_P$ (again via the large eigenvalue function, an algebraic function), which are indexed by the unit interval.

\comment
However, there are some results (which don't intersect the approximately
divisible case) when $d = 1$. For example, suppose $P\in \Z[x]^+$ is
projectively faithful (meaning, the $\gcd$ of its exponents is $1$)
monic, its constant term is $1$ (so $R_P$ has exactly two pure discrete
traces). Let $r$ be a positive real number. Then (1) is satisfied for $r
\in \R^{++}$ iff $r$ is an algebraic unit (that is, $r$ satisfies a monic
integer polynomial with constant coefficient $\pm 1$, or equivalently,
$\Z[r] = \Z[1/r]$).

To see this write $P = 1 + \sum_{0< i < k} a_i x^i + x^k$ where $a_i$ are
nonnegative integers (some can be zero, but we still need $\gcd \(
\Set{i}{a_i \neq 0} \cup \brcs{k}\) = 1$). From $P(r) \in
\Z[\brcs{r^j/P(r)}{j\in \Log P}]$, we deduce an equation of the form
$P(r)^{m+1} = p(r)$ where $p \in \Z[x]$ and $\deg p \leq \deg P^m = km$.
The leading term of this expression is $r^{(m+1)k}$, and so $r$ satisfies
a monic polynomial with integer coefficients. Owing to the $\gcd$
condition, $\Log P \setminus \brcs{0,k}$ is not empty, so there exists $i$
\st $a_i \neq 0$; letting $i'$ be the smallest such, we have $P(r)/r^{i'}
\in [\Z[\brcs{r^j/P(r)}{j\in \Log P}]$, and it easily follows (from the
constant coefficient of $P$ being $1$) that $1/r^{i'}$ is a polynomial in
$r$, hence $r^{i'}$ is invertible in $\Z[r]$, and thus $r'$ is invertible
in $\Z[r]$, so $r'$ is an algebraic unit.

Conversely, suppose $r$ is an algebraic unit. Then $P(r)/r^j$ are all
integral (for $j \in \Log P$), hence $P(r)/r^i \in \Z[r^i/P(r)] $ for
each $i$, and thus $P(r)/r^i \in\Z[\brcs{r^j/P(r)}{j\in \Log P}]$ for all
$i \in \Log P$.

We can get necessary conditions rather easily. If $P = a_0 + \sum a_i x^i
+ a_k$, then $r$ satisfies (1) implies the prime factorization of $N(r)$
is $\prod p^{e(p)}/\prod q^{f(q)}$ where $p $ divides $a_0$ and $q$
divides $a_k$. (If $\gcd{a_0,a_k} \neq 1$, the set of these products
contains an infinite group; otherwise, it is merely a semigroup.)

For $G = \Q[X]$ with $G^+ = \Set{f\in \Q[X]}{f|(0,1)> 0}$; this is the
Pascal's triangle dimension group (see Renault [Rxxx]) tensored with $\Q$.
In this case, the $\tau(I) = \tau(R)$ condition has no effect, as we will
see.

On the other hand, let $P = 2 + 3x \in \Z[x]$, set $R = \Z[x/P, 1/P]$ (a
subring of $S=\Z[x,1/P]$) with $R^+$ generated additively and
multiplicatively by $\brcs{1/P, x/P}$. It can also be rewritten as $\Set
{f/P^k}{\deg f \leq 2k}$. This is a very very special case of the
constructions in [Hxxx, Hxxx]; it is free as an abelian group, but is
nearly divisible, as we will see. The previous example was the divisible
hull of the ring we would have obtained from $P = 1+x$ (set $X =
x/(1+x)$).
The pure traces are given by $\tau_r$ for $r \in [0,\infty]$, given by
$\tau_r (f/P^k) = f(r)/P(r)^k$, together with
$$
\tau_{\infty} (\frac f{P^k}): = \lim_{t\o \infty} \frac{f(t)}{P^k(t)} =
\frac{(f,x^{2k})}{(P^k,x^{2k})},
$$
where the inner product notation, $(g,x^r)$, indicates the coefficient of
$x^r$ in the polynomial $g$ (so the denominator in the displayed
expression is $3^{k}$).

The pure traces are multiplicative, and the only two with a positive
element in their kernel are $\tau_0$ and $\tau_{\infty}$, and in both
cases, their kernels are actually order ideals ([H,xxx]), so they are good
by default. In all other cases, $\tau_r$ extends uniquely and obvious to a
multiplicative trace on the ordered ring $S = \Z[x,1/P]$; unfortunately,
this has no order unit, so we have to be careful in dealing with it.

The condition $\tau(I) = \tau(G)$ for all order ideals, a necessary
condition for $\tau$ to be good, turns out to be amazingly restrictive.
For example, if $r$ is rational, then $\tau_r$ is good iff $r \in
\brcs{2^i3^{-j}}_{i,j\geq 0}$, and the constraints if $r$ is merely
algebraic are quite strict.

Let $I_0$ and $I_{\infty}$ be the ideals $(1/P)R$ and $(x/P)R$; these turn
out to be order ideals (it is rarely true that the order ideal generated
by a positive element $a$ is $aR$, the principle ideal it generates, but
it happens here for these two choices of $a$ [H,xxx]), and in fact
(because $P$ has only two terms), $\ker \tau_i = I_i$ for $i \in
\brcs{0,\infty}$, as is easy to see.

Thus for $r \in (0,\infty)$, $\tau_r(I_0) = (1/P(r))\tau(R)$, $\tau_r(I_1)
= (r/P(r)\tau(R))$, and $\tau(R) = \Z[1/P(r),r/P(r)] = \Z [1/(2+3r),
r/(2+3r)] \subset \R$. It is tempting to conjecture that the last is just
$\Z[1/(2+3r)]$ in general, but it isn't as very easy evaluations show.
Then both $\tau_r(I_i) = \tau(R)$ iff $P(r)/r$ and $P(r)$ belong to
$\tau(R)$ (it is enough to get $1 \in \tau_r (I_i)$. Thus we have,
$\tau_r(I_i) = \tau(R)$ iff both $2/r$ and $3r$ belong to $ \Z [1/(2+3r),
r/(2+3r)] $.

For example, if $r = a/b$ with $a,b$ relatively prime positive integers
(that is, $r \in \Q^{++}$, then $\tau_r(R) = \Z[1/(2+3a/b), a/(2b+3a)] =
\Z[b/(2a+3b), a/(2a + 3b)]$; since $\gcd (a,b) = 1$, $\tau(R) = \Z[1/(2a
+3b)]$. Necessary and sufficient that $3a/b$ belong to this is that $b =
3^j$ for some $j \geq0$ (using relative primeness of $(a,b)$), and $2/r =
2b/3a$ belongs iff $a = 2^i$ for some $i \geq 0$. Hence if $\tau_{r}$ is
good and $r$ is rational, then $r \in \brcs{2^i3^{-j}}_{i,j\geq 0}$.
(Rather oddly, the set of rational $r$ \st $\tau_r$ is good is closed
under multiplication.)

Suppose $3r, 2/r \in \tau_r(R)$; then $N(r) \in \brcs{ 2^i 3^j}_{i,j \in
\Z} $ (this is not the full converse, because it allows negative
exponents for $2$ and positive exponents for $3$; but it is a first step.

To see this, invert $6$, forming $A = \Z[1/6]$; then the condition boils
down to $r, 1/r \in A[1/(2 + 3r)]$ (since $r/(2+3r) = (1/3)(1-2/(2 +
3r))$). The first says there exist $a_i \in A$ \st $r =\sum_{i=0}^k
a_i/(2+3r)^i$ with $a_k \neq 0$; this yields $r(2+3r)^{k} = \sum_{0< i\leq
k} a_i (2+3r)^{k-i}$, yielding $3^k r^{k+1}$ is a polynomial in $r$ of
lesser degree with coefficients in $A$; since $3$ is invertible, this says
$r$ is an algebraic integer over $A$; similarly, $1/r$ is an algebraic
integer over $A$. Hence $N(r) \in A* = \brcs{\pm 2^i 3^j}$.

(In general, if $K$ is a finite-dimensional extension field of $\Q$, then
the norm $N$ has values in $\Q$; however, if $r$ is an algebraic integer
over a subring $A$ of $\A$, then $N(r) \in A$; if in addition, $r$ is an
invertible algebraic integer---that is, its inverse is an algebraic
integer over $A$, $N(r)$ is invertible in $A$, as $N$ is multiplicative.)

Now we wish to exclude the positive exponents of $3$ and the negative
exponents of $2$. Form $B = \Z[1/3]$; then $r/(2+3r) = (1/3)(1 -
2/(2+3r))$, so the conditions become, $r, 2/r$ belong to $B[1/(2+3r)]$; as
before, $r$ is an algebraic integer over $B$, so $N(r) \in B$. From $2/r
\in B[1/(2+3r)]$, we deduce $(2/r)(2+3r)^k = \sum a_i (2+ 3r)^{k-i}$, so
that $2^{k+1}/r$ is a polynomial in $r$ with coefficients from $B$. Thus
$2^{k+1}/r$ is integral over $B$, so that $N(2^{k+1}/r) \in B$. Since
$N(2^{k+1})$ is a power of $2$, we deduce that $N(r)$ divides a power of
$2$ in $B$. Hence $N(r) \in \brcs{\pm 2^i3^{j}}_{i \geq 0, j \in \Z}$.

Form $C = \Z[1/2] $; this time we rewrite $1/(2+3r)$ in terms of
$r/(2+3r)$, and interchange the roles of $r$ and $1/r$. Set $R = 1/r$;
then $1/(2 + 3r) = (1/2)(1 - 3/(2R+3))$, so $C[1/(2 + 3r), r/(2+3r)] =
C[1/(2R + 3)]$, and the conditions become
$3/R, R/2 \in C[1/(2R + 3)]$; as $2$ is invertible, the latter condition
is $R \in C[1/(2R+3)]$. Now we can interchange the roles of $2$ and $3$
from the previous paragraph, and deduce that $N(R) $ divides a power of
$3$ in $C$, so $N(R) \in \brcs{3^{-j}2^i}_{j \leq 0, i \in \Z}$. Hence
$N(r) = 1/N(
R)$ belongs to $\brcs{2^i3^{-j}}_{i \in \Z, j \geq 0}$. Combining this
with the result of the previous paragraph, we have $N(r) \in \brcs{2^i
3^{-j}}_{i, j \geq 0}$.

\endcomment


\SecT 5 Direct sums and goodness

For (noncyclic) simple dimension groups, there is a notion of direct sum (corresponding to coproduct; see [BeH, Appendix 2] for a discussion). This actually extends to nearly divisible dimension groups. Let $G$ and $H$ be nearly divisible. Form the group direct sum $G \oplus H$, and impose on it the ordering given by the positive cone
$$
\Set{(g,h)}{g \in G^+ \setminus \brcs{0} \text{ and }h \in H^+ \setminus \brcs{0} } \cup \brcs{(0,0)}.
$$
When $G$ and $H$ are simple (and noncyclic) dimension groups, the resulting {\it strict direct sum\/} $G \oplus_s H$ ({\it s\/} for strict) is also a simple dimension group. It is each to check that when both $G$ and $H$ are nearly divisible, then so is $G \oplus_s H$ (trivial), and when both are additionally dimension groups, so is the strict direct sum. We suppress the subscript $s$.

If $K = G \oplus H$, and $\sigma$ and $\tau$ are traces on $G$ and $H$ respectively, we consider the possibility that $\phi := \sigma \oplus \tau$ (defined by $(g,h) \mapsto \sigma(g) + \tau(h)$) be good or order unit good.
If we put the usual direct sum ordering on $K$ and consider the question of characterizing when $\phi$  is good, the answer is not exciting. However, if we put the strict direct sum ordering on $K$, then  order unit goodness is interesting. Iteration of this process yields some weird examples.

\Lem Lemma \fivone. [A consequence of the method of proof of [BeH, Proposition 1.7]] Suppose $(K,w)$ is
an approximately divisible dimension group with order unit, and $\phi$ is
an order unit good trace. Then whenever $a \in G$, $b \in G^{++}$ and $0 <
\phi(a) < \phi(b)$, for all $\epsilon > 0$, there exists $a' \in [0,b]$
\st $\phi(a') = a$ and $\| \hat a' - \hat b \sigma(a)/\sigma(b) \| <
\epsilon$.

\Pf Approximate divisibility implies density of $G$ in $\Aff S(G,u)$. Set
$j = \sigma(b)\hat b/\sigma(a)$, so that $j(\sigma) = \sigma (a)$ and $\inf j
= \sigma(a)\sigma(b)^{-1}\inf \hat b$. There exists $g_n \in G$ \st
$\hat g_n \to j$ uniformly. If for infinitely many $n$, $g_n(\sigma) =
\sigma(a)$, we are done (taking large enough $n$). Otherwise, select
$\sigma(a)(\sigma(b)2)^{-1}\inf b > \epsilon > 0$ and $\| \hat g_n - j\| <
\epsilon$, then  $|\sigma(g_n) - \sigma(a )| < \epsilon$ provided $n$ is
sufficiently large; if $\sigma(g_n) > \sigma(a)$, set $c_n = g_n - a$.
There exists an order unit $z_n$ \st $0 < \sigma(c_n) \pmb 1 < \hat z_n <
2\epsilon$. By order unit goodness, there exists $v_n \ll z_n$ \st
$\sigma(c_n) = \sigma(v_n)$, and of course, $\|v_n \| \leq \| \hat z_n\| <
2\epsilon$. Then $g_n - v_n$ has image within $3\epsilon$ of $j$, and it
is easy to check that $g_n - v_n$ is strictly positive, hence is an order
unit.

If instead, $\sigma(a) > \sigma(c_n)$ for infinitely many $n$, we obtain a
corresponding $c_n = g_n - a$ and $v_n \ll z_n$, and this time, $g_n +
v_n$ has all the right properties. In both cases, by taking $n$
sufficiently large, we make the error terms go to zero, hence obtain the
$a'$ as one of  $g_n \pm v_n$.\qed

In the following, the function $\psi$ will not be a group homomorphism (just a function, and usually a weird one, if it exists at all).

\Lem Lemma \fivtwo. Suppose $G$ and $H$ are nearly divisible dimension groups, each with order unit, and
respective trace $\sigma$ and $\tau$. Let $K = G\oplus H$ with the strict ordering, and
and suppose that the trace on $K$, $\phi := \sigma \oplus \tau$ is order
unit good. Then provided the following condition holds, $\sigma$ is order
unit good as a trace on $G$:
\item{} there exists a function $\Arrow \psi; \tau^{-1}(\sigma(G) \cap
\tau(H).\sigma^{-1}(\sigma(G) \cap \tau(H))$ that is pseudo-norm
continuous with the additional property that $\sigma\psi = \tau$.

\Rmk As we will see below, without the weird extra condition, the result
fails.

\Pf Select an order unit $b$ in $G$, and $a$ in $G$ \st $0 < \sigma(a) <
\sigma(b)$. As $H$ is approximately divisible, there is a sequence of
order units in $H$, $(h_n)$ \st $h_n \to 0$ (\wrt the pseudo-norm topology
on $H$; equivalently, as functions on $S(H,v)$, $\hat h_n$ converges
uniformly to zero). There also exists $\delta $ in $G$ \st $\sigma(b-a)/4
< \hat \delta < \min\brcs{\sigma(b - a)/2, \inf_{\theta \in S(G,u)} \theta(b)/2}$
uniformly on $S(G,u)$. Then $B_n:= (b- \delta,h_n)$ are order units of $G \oplus H$, and
$\phi(a,0) < \phi(B_n) = \sigma(b) - \sigma (\delta) + \tau (h_n)$.

Since $\phi$ is order unit good and each $B_n$ is an order unit, there
exist $(a_n, z_n)$ \st $0 \ll a_n \ll b- \delta$ and $0 \ll  z_n \ll h_n$
with $\phi((a_n, z_n)) = \sigma(a)$, and by the previous lemma,
$\inf_{S(G,u)} \hat a_n$ is bounded below (as $n \to \infty$); in
particular, $\| z_n\|_H \to 0$ and $\sigma(a_n) + \tau(z_n) = \sigma(a)$.
Thus $z_n \in  \tau^{-1}(\sigma(G) \cap \tau(H))$, so we may consider the
sequence $\psi(z_n) \in \sigma^{-1}(\sigma(G) \cap \tau(H)$. Since $\psi$
is pseudo-norm continuous, $\hat \psi(z_n) \to 0$ uniformly on $S(G,u)$.

Consider $a_n + \psi (z_n)$; its value at  $\sigma$ is $\sigma (a_n) +
\sigma(\psi(z_n)) = \sigma(a_n ) + \tau(z_n) = \sigma(a)$. If we choose
$n$ sufficiently large that $\| \widehat {\psi (z_n)}\| < \inf \delta$, then
$a_n + \psi(z_n) \ll b - \delta + \delta = b$. In addition,  we can also
choose $n$ sufficiently large that $\inf \widehat {\psi(z_n) }> - \inf_{S(G,u)}
\hat a_n$, by the uniform boundedness below of the $a_n$ (there is no
guarantee that $\psi(z_n)$ is positive). Then $a_n + \psi(z_n)$ is an
order unit in the interval $[0,b]$ and we are done.
\qed

\comment

\Lem Corollary \fivthr. Let $G$ be a nearly divisible dimension group with trace $\tau$. If $K = G \oplus G$ is equipped with the strict ordering and $\phi = \tau \oplus \tau$
is order unit good, then $\tau$ is order unit good.

This contrasts with the
situation with ordinary direct sum ordering, even with simplicial groups
(but order unit goodness), where we can even find, for each $k$, $(\Z^2, t)$
which is not order unit good,  $t \oplus t \oplus \dots \oplus t$ is order
unit good if and only if the number of direct summands is $k$ or more
(take $t = (1,k+1)$.

There exist simple $G, u, \sigma$ and $H,v, \tau$ \st $\phi = \sigma \oplus
\tau$ is good on $K = G \oplus H$ with the strict ordering, but $\sigma$
is not good (and we can arrange that either $\tau$ is good or not good).

\Lem Lemma \fivfou. Suppose that $(G,u,\sigma)$ and $(H,v,\tau)$ are dimension
groups with order unit and (not necessarily normalized) designated trace,
and let $g$ and $h$ be order units of $G$ and $H$ respectively,  \st the
following hold.
\item{(i)} $g$ is $\sigma$-good and $h$ is $\tau$-good;
\item{(ii)} $\sigma(G) \cap \tau (H)$ is dense in $\R$.
{\par} \noindent Then $(g,h)$ is $\phi$-good, where $\Arrow \phi; K:= G \oplus H.
\R$ is the trace on the strict direct sum of $G$ and $H$ given by
$\phi((a,b)) = \sigma(a) + \tau(b)$.

\Pf An immediate consequence of Lemma \onetwo(d).

\qed
\endcomment
\comment
Suppose $0<\phi((a,b)) < \sigma(g) + \tau (h)$. We proceed in several
steps.
\noindent {(i)} We can assume $\sigma (a) < \sigma(g)$ and $\tau(b) <
\tau(h)$. If $\sigma(a) > \sigma(g)$, then $\tau(b) = \phi((a,b)) -
\sigma(a) < \phi((a,b)) - \sigma(g) < \tau(h)$. Select $\delta > 0 $ \st $\sigma(a)
- \delta < \sigma (g)$ and $\tau(b) + \delta < \tau(h)$ (to see that such
$\delta$ exists, we need only note that $\sigma (g) - \sigma(a) < \tau(h)
- \tau (b)$). There exists a neighbourhood of $\delta$ that satisfies the
same inequalities. By density, there exist elements $e \in G$ and $f \in
H$ \st $\sigma(e) = \tau(f)$  and they belong to this interval. Then set
$a' = a-e$ and $b' = e+f$; we have $\phi((a',b')) = \phi(a,b)$, but now
$\sigma(a') <\sigma(g)$ and $\tau(b') < \tau (h)$.

If instead, $\tau(b) > \tau(h)$, just reverse the roles of $a$ and $b$.

\noindent {(ii)}. We can assume that $0 < \sigma (a) < \sigma(g) $ and $0
< \tau(b) < \tau (h)$. Using case (i), suppose that $\sigma(a) \leq 0$ and
$\tau(b) < \tau(h)$; then necessarily $\tau(b) =\phi((a,b)) - \sigma(a)
\geq \phi((a,b))$. We want to find $\delta> 0$ \st $\sigma(g) > \sigma(a)
+ \delta > 0$, $0 < \tau(b) - \delta < \tau(h)$; the existence of such a
$\delta$ (and its corresponding neighbourhood) is equivalent to the four
inequalities
$$
\sigma(g) - \sigma(a), \tau(b) > -\sigma (a), \tau(b) - \tau (h);
$$
each of the four is easy to check (the most interesting being $\sigma(g) -
\sigma(a), \ >  \tau(b) - \tau (h)$, which is equivalent to $\phi((g,h)) >
\phi((a,b))$). As in the previous argument, we can find $e \in G$ and $f
\in H$ \st $\sigma(e) = \tau(f)$ and lies in the interval. Then we replace
$(a,b)$ by $(a' = a+e, b' = b-f)$.

Having case (ii), by order unit goodness of each of $\sigma$ and $\tau$,
there exists $a'' $ with $0 \ll a'' \ll g$ and $b''$ with $0 \ll b'' \ll
h$ with $\sigma(a'') = \sigma(a)$ and $\tau(b'') = \tau(b)$; then
$(a'',b'')$ lies in the interval $[0, (g,h)]$ in $K$ and its value at
$\phi$ is $\sigma(a) + \tau(b) = \phi((a,b))$, and we are done.
\endcomment
\comment
\Lem Corollary \fivfiv. Suppose that $(G,u,\sigma)$ and $(H,v,\tau)$ are nearly divisible dimension
groups with order unit and (not necessarily normalized) designated trace,
and let $g$ and $h$ be order units of $G$ and $H$ respectively,  \st the
following hold.
\item{(i)} $\sigma$ and $
\tau $ are order unit good;
\item{(ii)} $\sigma(G) \cap \tau (H)$ is dense in $\R$.
{\par} \parindent=0em \noindent Then  $\phi = \sigma \oplus \tau$
is an order unit good trace on $K$ \wrt  either the strict or the usual direct sum ordering on $K = G\oplus H$.

In the presence of near divisibility, (ii) is necessary, by \onethr. However, (i) is not
(neither, or only one of them can be good), as simple examples show (later). There is an alternative proof of the following which yields more, but requires more definitions.
\endcomment

One advantage of not requiring normalization of $\sigma$ and $\tau$ is that we can replace them by any positive scalar multiples, in testing for order unit goodness of $\lambda \sigma \oplus \mu \tau$; the first hypotheses are unchanged, but the second translates to density of $(\lambda \sigma(G) )\cap (\mu \tau(G))$  in $\R$. In the following, we cannot apply earlier results directly, since $G \oplus 0$ is not an order ideal of $G \oplus H$ (strict ordering).

\Lem Lemma \fivsix. Suppose that $\sigma$ is a trace on $G$, $\tau $ is a trace on
$H$, and $\sigma \oplus \tau = \phi$ is order unit good for $K = G \oplus
H$ with the strict ordering, and moreover assume that each of $G$ and $H$
is approximately divisible. Then $\sigma(G) \cap \tau(H)$ is dense in
$\R$.

\Pf We use the characterization of order unit good traces on approximately
divisible dimension groups, namely $\ker \phi$ has dense image in
$\phi^{\vdash}$ [BeH, Proposition 1.7].

Suppose the intersection is not dense; then exists real $\delta \geq 0$
\st $\sigma(G) \cap \tau(H) = \delta \Z$. We have that $\ker \phi$ has
dense range in $\Aff S(K,(u,v)) = \Aff S(G,u) \times \Aff S(H,v)$. But
$\ker \phi = \Set{(g,h) \in G \oplus H}{\sigma(g) = - \tau(h)}$.

If $\delta = 0$, then $\ker \phi = \ker \sigma \oplus \ker \tau$ (since
$\sigma(g) = -\tau(h)$ implies $\sigma(g) \in \tau(H)\cap \sigma(H)$,
hence is zero). The image of $\ker \phi$ is then contained in
$\sigma^{\vdash} \times \tau^{\vdash}$, which is closed and  of
codimension two in $\Aff S(K,(u,v)) $, and so $\ker \phi$ cannot be
dense in $\phi^{\vdash}$ (which as codimension one), hence $\phi$ cannot
be order unit good.

If $\delta \neq 0$, select $g$ and $h$ in $G$ and $H$ respectively \st
$\sigma(g) = \delta = \tau(h)$. Then it is easy to see that  $\ker \phi =
(\ker \sigma \oplus \ker \tau) + (g,-h) \Z$, and then its range is
contained in $(\sigma^{\vdash} \times \tau^{\vdash}) + (\hat g, -\hat
h)\Z$. However, the latter is closed (easy to see), and so the image of
$\ker \phi$ is contained in a proper closed subspace (with a discrete
direct summand) of $\phi^{\vdash}$, hence in this case as well, $\phi $ is
not order unit good.\qed

Now we want to determine when $\sigma \oplus \tau$ is good or order unit good. Let $\Arrow \pi_G;  G\oplus H.G$ and $\Arrow \pi_H;  G\oplus H.H$ be the obvious projection maps. Unlike the inclusions $G,H \to G \oplus H$, these {\it are\/} order-preserving.
First, consider $\Arrow
\sigma \circ \pi_G; \ker \phi.\sigma(G) \cap \tau(H) \subseteq \R$. The
kernel is exactly $\ker \sigma \oplus \ker \tau$; we also note that
$\sigma$ extends to a map $\Arrow \Sigma; \phi^{\vdash} . \R$ (sending
$(j,l)$ to $j(\sigma)$), the kernel of which is $\sigma^{\vdash} \times
\tau^{\vdash}$. With the identification of $\Aff S(K,(u,v))$  with $\Aff S(G,u) \times \Aff S(H,v)$, we have the following diagram.

$$
\diagram
0 &\rTo^{} & \ker \sigma \oplus \ker\tau & \rTo^{} &
\ker \phi & \rTo^{\sigma\circ\pi_G} & \sigma(G) \cap \tau(H) &\rTo &0\\
&&\dTo^{}& & \dTo^{} &&\dTo^{} &&& \\
& & \overline{\widehat{\ker \sigma}} \times \overline{\widehat{\ker \tau}}& \rTo^{} &
\overline{\widehat{\ker \phi}} & \rTo^{\Sigma} & \R  &&\\
&&\dTo^{}& & \dTo^{} &&\dTo^{} &&& \\
0&\rTo^{} & \sigma^{\vdash} \times \tau^{\vdash}& \rTo^{} &
\phi^{\vdash} & \rTo^{\Sigma} & \R  &\rTo &0\\
\enddiagram
$$

The long horizontal overlines indicate closure, as subgroups of the affine function vector spaces; of course, there is no requirement that any of the three overlined groups be real vector spaces (they are norm-complete subgroups).
The  two leftmost  top vertical arrows are just induced by the affine representations; the right one is the inclusion, compatible with $\Sigma$ restricted to the image of $\ker \phi$. The two leftmost bottom vertical arrows are the obvious inclusions. The $\Sigma$ in the middle row is an abuse of notation; it represents the restriction of $\Sigma$ to $\overline{\widehat{\ker \phi}}$, the closed subgroup of $\phi^{\vdash}$, but the notation is already rather complex.

The middle row need not be exact at either end (for example, if $\ker \phi$ has dense image in $\phi^{\vdash}$ but one or both of $\ker \sigma$ or $\ker \tau$ does not have dense image in $\sigma^{\vdash}$ or respectively $\tau^{\vdash}$, then it is not left exact; if $\sigma(G) \cap \tau(H)$ is discrete, then the middle line is not right exact).

If $\ker \phi$ has dense image in $\phi^{\vdash}$, then $\sigma(G) \cap \tau(H)$ is a dense subgroup of $\R$: we simply note that density of the image of $\ker \phi$ in $\phi^{\vdash}$, the latter being a closed and therefore a norm-complete subspace of $\Aff S(K,(u,v))$, entails that for every bounded linear functional that is not zero on $\phi^{\vdash}$, its restriction to a dense subgroup must be not zero and have dense range in the reals.
\qed

It also leads to a straightforward proof that if $\ker \sigma$ and $\ker \tau$ have dense images in $\sigma^{\vdash}$ and $\tau^{\vdash}$ respectively, and if $\sigma(G) \cap \tau(H)$ is a dense subgroup of $\R$, then $\ker \phi$  has dense image in $\phi^{\vdash}$. We have that
$\sigma^{\vdash} \times \tau^{\vdash} = \overline{\widehat{\ker \sigma}} \times \overline{\widehat{\ker \tau}} \subseteq \overline{\widehat{\ker \phi}} \subseteq \phi^{\vdash}$. The left and   right terms of these inclusions are vector spaces, and since $\sigma^{\vdash} \times \tau^{\vdash}$
is a closed codimension two subspace of  $\Aff S(K,(u,v))$ and $\phi^{\vdash}$ is codimension one, it follows that $\sigma^{\vdash} \times \tau^{\vdash}$ is a codimension one subspace of $\phi^{\vdash}$. (The proof does not stop herewe do not know that $\overline{\widehat{\ker \phi}}$ is a real vector space.)

The map $\Sigma$ induces $\overline {\widehat{\ker \phi}}/(\overline{\widehat{\ker \sigma}} \oplus\overline{\widehat{\ker \tau}})$ to a subgroup of the reals. However, this subgroup of the reals includes the dense subgroup $\sigma(G) \cap \tau(H)$, and as $\overline{\widehat{\ker \phi}}$ is a norm-complete abelian group, the image must be complete, and thus must be onto. In addition, since $\ker \Sigma = \sigma^{\vdash} \times \tau^{\vdash} = \overline{\widehat{\ker \sigma}} \times \overline{\widehat{\ker \tau}}$, it follows that $\ker \Sigma \cap \overline{\widehat{\ker \phi}} = \overline{\widehat{\ker \sigma}} \times \overline{\widehat{\ker \tau}}$. We thus have $\ker \Sigma \subset \overline{\widehat{\ker \phi}} \subseteq \phi^{\vdash}$, but $\Sigma$ induces equality
$\overline{\widehat{\ker \phi}}/(\ker \Sigma \cap \overline{\widehat{\ker \phi}}) = \phi^{\vdash}/\ker \Sigma$. It follows immediately that $\overline{\widehat{\ker \phi}} = \phi^{\vdash}$. \qed

Now we can show that if the closure of the images of $\ker \sigma$ and $\ker \tau$ are real vector spaces, and if $\ker \phi$ is order unit good, then $\sigma $ and $\tau$ are order unit good.

We wish to show $\overline{\widehat{\ker \sigma}} \times \overline{\widehat{\ker \tau}} = \sigma^{\vdash}\times \tau^{\vdash}$ (from this it follows trivially that $\overline{\widehat{\ker \sigma}} = \sigma^{\vdash}$ and $ \overline{\widehat{\ker \tau}} = \tau^{\vdash}$). Since the left thing is a vector space, and a complete normed abelian group (hence a closed vector subspace of $\Aff S(K, (u,v))$, if equality does not hold, there exists a bounded linear functional $\alpha$ on $\Aff S(K,(u,v))$ that kills $\overline{\widehat{\ker \sigma}} \times \overline{\widehat{\ker \tau}}$ but not $\sigma^{\vdash}\times \tau^{\vdash}$; in particular, $\alpha$ does not kill $\phi^{\vdash}$.

By composition with the affine representation, we \quotes{restrict} $\alpha$ to a real-valued bounded group homomorphism  $\Arrow \beta; G \oplus H.\R$ (for a treatment of bounded group homomorphisms on dimension groups, see [G]; their behaviour is just what you'd expect). Since $\alpha$ kills $\overline{\widehat{\ker \sigma}} \times \overline{\widehat{\ker \tau}}$, it follows that $\beta$ kills $\ker \sigma \oplus \ker \tau$. We form the normed abelian group $\ker \phi/(\ker \sigma \oplus \ker \tau)$, which via $\sigma$, we know to be  $\sigma(G) \cap \tau (H) \subset \R$. Thus $\beta$ induces a  bounded real-valued group homomorphism on $\ker \phi/(\ker \sigma \oplus \ker \tau)$, call it $\overline \beta$. We thus have two bounded group homomorphisms on the quotient, $\overline \beta$ and $\overline \sigma$, but as the quotient is isomorphic (as a normed abelian group) to a subgroup of the reals, there must be a positive real number $\lambda$ \st $\overline \beta = \lambda \overline \sigma$.

This forces $\beta = \lambda \cdot \sigma \circ \pi_G$ (as bounded group homomorphisms  on $\ker \phi$). Since $\ker \phi$ has dense image in its completion (!) which happens to be $\phi^{\vdash}$, we have that $\alpha|\phi^{\vdash} = \lambda \Sigma$. Thus $\alpha$ kills $\sigma^{\vdash} \times \tau^{\vdash}$, a contradiction. \qed

To summarize, we have the following results.

\Lem Proposition \fivnin. Suppose $(G,u,\sigma)$ and $(H,v,\tau) $ are approximately divisible dimension groups with order unit and distinguished trace, and form $K = G \oplus H$, and the trace $\Arrow \phi = \sigma \oplus \tau;K.\R$.
\item{(a)} If $\phi$ is order unit good (\wrt either the usual or the strict ordering on $K$), then $\sigma(G) \cap \tau(H)$ is a dense subgroup of the reals, and $\sigma \otimes 1_{\Q}$ and $\tau \otimes 1_{\Q}$ are order unit good as traces on $G \otimes \Q$ and $H\otimes \Q$ respectively.
\item{(b)} If the closure of the images of $\ker \sigma$ and $\ker \tau$ in $\sigma^{\vdash}$ and $\tau^{\vdash}$ respectively are real vector spaces, and if $\phi$ is order unit good, then both $\sigma$ and $\tau$ are order unit good.
\item{(c)} If $\sigma$ and $\tau$ are order unit good and $\sigma(G) \cap \tau(H)$ is dense in $\R$, then $\phi$ is order unit good.

Examples exist (given below) with $G$ and $H$ simple dimension groups to show that if $\phi$ is order unit good, then neither $\sigma$ nor $\tau$ (or exactly one of them) need be order unit good.


This method suggest what we should do with multiple traces. Let $(G_i,u_i,
\sigma_i)$, $i=1, 2, \dots, n$, be approximately divisible dimension
groups, each with order unit and (unnormalized) trace. Form $K = \oplus
G_i$ with the strict ordering, and $\Arrow \phi = \sigma_1 \oplus \sigma_2
\oplus \dots \oplus \sigma_n; K.\R $, and the map $\Arrow T; K.\R^n$
defined by $\phi((g_i)) = \sum \sigma_i(g_i)$ and $T((g_i)) =
(\sigma_1(g_1), \sigma_2(g_2), \dots, \sigma_n (g_n))$. Identify $\Aff
S(K,((u_i)))$ with the cartesian product $\Aff S(G_1,u_1) \times \dots
\times \Aff S(G_n,u_n)$.

If $(g_i) \in \ker \phi$, then $\sigma_n (g_n) = -\sum_{i=1}^{n-1}
\sigma_i (g_i)$, and we  can interchange $n$ with any other integer less
than $n$. In particular,
$V:= \sigma_n^{-1}\(\sigma_n(G_n) \cap \(\sum_{i=1}^n \sigma_i(G_i)\)\)$
is independent of permutations  and the range of $T$ on $\ker \phi$ is $T(V)$.

Extend $T$ to $\Arrow \Cal T; \Aff S(K, (u_i)).\R^n$ (sending $(j_i)$ to
$(j_i(\sigma_i))$. Restricted to $\phi^{\vdash}$, the range of $\Cal T$ is
exactly $(1,1,1,\dots,1)^{\perp}$, i.e., the entries add to zero.

Now we can form the  diagram analogous to the previous one.

$$
\diagram
0 &\rTo^{} & \ker \sigma_1 \oplus \dots \oplus \ker \sigma_n & \rTo^{} &
\ker \phi & \rTo^{T} &
T(V) &\rTo &0\\
&&\dTo^{}& & \dTo^{} &&\dTo^{} &&& \\
& & \overline{  \widehat{ \ker \sigma_1}} \times \dots \times \overline{
\widehat{ \ker \sigma_n}}  & \rTo^{} &
\overline{\widehat{\ker \phi}} & \rTo^{\Cal T} & \overline{T(V)} &&\\
&&\dTo^{}& & \dTo^{} &&\dTo^{} &&& \\
0&\rTo^{} & \sigma_1^{\vdash} \times \dots \times \sigma_n^{\vdash}& \rTo^{} &
\phi^{\vdash} & \rTo^{\Cal T} &  (1,1,\dots, 1)^{\perp} &\rTo &0\\
\enddiagram
$$

We quickly see that density of $T(V)$ (a subgroup of $\R^n$ contained in
$(1,\dots,1)^{\perp}$) in $(1,\dots,1)^{\perp}$ is necessary that $\phi$
be order unit good, that is, it is necessary in order that $\ker \phi$ have norm
dense image in $\phi^{\vdash}$.

Suppose all the $\sigma_i$ are order unit good and $T(V)$ is dense in
$(1,\dots,1)^{\perp}$. Then $\overline{ \widehat{\ker \sigma_1}} \times \dots \times \overline{
\widehat{ \ker \sigma_n}} = \sigma_1^{\vdash} \times \dots \times \sigma_n^{\vdash}$ is a closed subspace of $\phi^{\vdash}$, and  the middle line  yields that the codimension
of $  \overline{  \widehat{ \ker \phi}}  $ in $\Aff S(K)$ is $n -(n-1) =
1$, so being a closed subspace of the codimension one space
$\phi^{\vdash}$,   $\overline{  \widehat{ \ker \phi}} $ must equal it, and
therefore $\phi$ is order unit good.

Suppose $\phi$ is order unit good (hence we have density of $T(V)$ in
$(1,\dots,1)^{\perp}$) and each of $\overline{  \widehat{ \ker \sigma_1}}
$ is a vector space. To show $\sigma_i$ are all order unit good,
sufficient is that $\ker \sigma_i$ have dense image in $\sigma^{\vdash}$,
and it is easy to show sufficient for this is that $\overline{  \widehat{
\ker \sigma_1}} \times \dots \times \overline{  \widehat{ \ker
\sigma_n}}$ equals $\sigma_1^{\vdash} \times \dots \times
\sigma_n^{\vdash}$.

We note that the bounded real-valued group homomorphisms on $T(V)$ and of
course on its closure are linear combinations of the coordinate functions,
which correspond to the $\sigma_i$, with lack of uniqueness arising from
the relation that the sum of the coefficients is zero.

By assumption, each $\overline{  \widehat{ \ker \sigma_i}} $ is a vector
space (and of course closed in $\Aff S(G_i,u_i)$, whence the whole batch
$L:= \overline{  \widehat{ \ker \sigma_1}} \times \dots \times \overline{
\widehat{ \ker \sigma_n}}$
is a closed subspace of $M:= \sigma_1^{\vdash} \times \dots \times
\sigma_n^{\vdash}$ (which is itself a closed codimension $n$ subspace of
$\Aff S(K)$). If they are not equal, there exists a bounded linear
functional $\alpha$ on $\Aff S(K,(u_i))$ that kills $L$ but not $M$. This
induces a bounded real-valued group homomorphism $\beta$ on $\ker \phi$,
which kills $ \ker \sigma_1 \oplus \dots \oplus \ker \sigma_n$. Hence it
induces a bounded real-valued group homomorphism on the quotient, $T(V)$,
$\Arrow B; T(V).\R$.

Each $\sigma_i$ induces $\Sigma_i$ on $T(V)$, and these are the coordinate
functions. Hence there exist real $\lambda_i$ \st $B = \sum \lambda_i
\Sigma_i$.
Thus $\beta - \sum \lambda_i \sigma_i$ vanishes identically on $\ker
\phi$, and by density, $\alpha = \sum \lambda_i \sigma_i$ (where
$\sigma_i$ is now interpreted as the map $(j_i) \mapsto j_i(\sigma_i)$ on
$\Aff S(K)$). But this obviously kills $\sigma_1^{\vdash} \times \dots
\times \sigma_n^{\vdash}$, a contradiction. Hence each $\sigma_i$ is order
unit good.

To summarize what happens with multiple traces:

\Lem Theorem \fivten. Let $(G_i, u_i, \sigma_i)$ be approximately divisible dimension groups with order unit ($u_i$) and (unnormalized) trace $(\sigma_i)$. Form $K = \oplus G_i$ (with the strict ordering), and the trace $\phi = \oplus \sigma_i$ on $K$. Set $J = \sigma_n(G_n) \cap \(\sum_{i \leq n-1}\sigma_i(G_i)\)$, a subgroup of $\R$.
\item{(a)} If $\phi$ is order unit good, then $T(\sigma_n^{-1}(J))$ is dense in $(1,1,\dots,1)^{\perp}$.
\item{(b)} If the closure of the image of $\ker \sigma_i$ in $\sigma_i^{\vdash}$ is a real vector space for all $i$, and if $ \phi$ is order unit good, then all $\sigma_i$ are order unit good.
\item{(c)} If the  image of $\ker \sigma_i$ is dense in $\sigma_i^{\vdash}$ for all $i$ (that is, each $\sigma_i$ is order unit good), and if $T(\sigma_n^{-1}(J))$ is dense in $(1,1,\dots,1)^{\perp}$, then $\phi$ is order unit good.

The conditions for order unit goodness with $n$ direct summands are
slightly different, in that they involve the density of a subgroup of
$\R^{n-1}$ (identified with $(1,\dots,1)^{\perp}$), or simply that the closure of $T(V)$
is a vector space of dimension $n-1$ (in general, the closure need not be
a vector space). To some extent, this explains some of the phenomena
illustrated in the examples below, with direct sums of two not yielding an order
unit good trace, but direct sums of three doing so. In fact, the argument
in the example, $G_n = \Z + (\sqrt 3 + n\sqrt 2)\Z$, essentially boils
down to showing the closure of $T(V)$ is a two-dimensional vector space.
But actually calculating with $T(V)$ seems awkward.

However, in some cases computation is moderately feasible. Suppose $G_1 =
\Z + \sqrt 6 \Z$, $G_2 = \Z + \sqrt {15}\Z$, and $G_3 = \Z + \sqrt
{10}\Z$. Then $T(V)$ is discrete, so $\sigma_1 \oplus \sigma_2 \oplus
\sigma_3$ is not order unit good. However, if we add a fourth term, $G_4 =
\Z + (\sqrt 6 + \sqrt {15} + \sqrt{10})\Z$, then $\ker \phi \cap
\sigma_4^{-1}(G_4 \cap (\sum G_i)) = \Set{(a +b\sqrt 6,c+b\sqrt{15} , d+
b\sqrt{10}, -(a+c + d) - b((\sqrt 6 + \sqrt {15} + \sqrt{10})))}{a,b,c,d
\in \Z}$. Let $v_1 = (1,0,0,-1)$, $v_2 = (0,1,0,-1)$,  and $v_3 =
(0,0,1,-1)$; then $\ker \phi$ is the $\Z$-span of $\brcs{v_1, v_2, v_3,
\sqrt 6 v_1 + \sqrt {15} v_2 + \sqrt{10}v_3}$. The map from $\ker \phi$ to
$\R^3$ given by $v_i \mapsto e_i$ (standard basis elements) has range the
free abelian group on $\brcs{e_1, e_2, e_3,  \sqrt 6 e_1 + \sqrt {15} e_2
+ \sqrt{10}e_3}$. Since $\brcs{1, \sqrt 6, \sqrt {10}, \sqrt {15}}$ is
linearly independent over the rationals, this group is dense. It is
trivial that $\brcs{v_i}$ is a real basis for $\phi^{\vdash}$, so $\phi$
is good. In this example, all the $\ker \sigma_i$ are trivial, so $T(V) $
is all of $\ker \phi$.

On the other hand, if we omit any one or two  of the $G_i$, the resulting
trace is not order unit good, since  the resulting $T(V) $ will be
discrete.

We can similarly construct $(G_i, \sigma_i)$ (the $G_i$ subgroups of the
reals), $i = 1,\dots, n$ \st $\oplus_{i=1}^n \sigma_i$ is order unit good,
but for no proper subset $J$ of $\brcs{1,2,\dots, n}$ with $|J|> 1$ is
$\oplus_{i\in J} \sigma_i$ order unit good: Let $\brcs{p_i}_{i=1}^n$ be
distinct primes; set $G_i = \Z + \sqrt {p_i} \Z$ for $1 \leq i \leq n-1$,
and $G_n = \Z +\( \sum_{i=1}^{n-1}\sqrt{p_i}\)\Z$. The resulting $T(V) $
will be a critical group of rank $n$, so all  subgroups of lesser rank are
discrete.

\Lem Example \fivsev. Simple dimension groups $(G, \sigma)$ and $(H, \tau)$
with traces \st $\phi = \sigma \oplus \tau$ is (order unit) good on
the strict direct sum $K = G \oplus H$, but $\sigma$ is not good as a
trace on $G$ (and in one case, $\tau$ is good, in another case, it is
not).

\Pf For simple dimension groups (as $G$, $H$, and $K$ are), order unit
goodness is equivalent to goodness. Begin with three subgroups of the
reals,
$$\eqalign{
G_1 & = \Z + \sqrt 3 \Z + \sqrt 5 \Z\cr
G_2 & = \Z + \sqrt 2\Z + \sqrt 5 \Z\cr
G_3 & = \Z + (\sqrt 3 + \sqrt 2)\Z.\cr
}$$
Let $\tau_i$ denote the respective identifications of $G_i$ with subgroups
of the reals; these are traces on each of these three totally ordered
dimension groups. Each $\tau_i$ is the unique (up to scalar multiple)
trace, so is good. Now form $L = G_1 \oplus G_2$ with the strict order;
since both subgroups contain $\Z + \sqrt 5\Z$, which is dense, it follows
from the criterion above that $\tau_1 \oplus \tau_2$ is a good trace
thereon. Next, form $K = L \oplus G_3$, with the strict order; since the
value group of $\tau_1 \oplus \tau_2 (L)$ includes $\Z + (\sqrt  3 + \sqrt
2)\Z$ and the latter is dense, we can apply the criterion again, and so
deduce that $\tau_1 \oplus \tau_2 \oplus \tau_3$ is good, as a trace on
$K$.

However, we can obtain $K$ by proceeding in a different order. Set $G =
G_1 \oplus G_3$ with the strict order. Either by direct examination or by
the necessity of the density condition, $\tau_1
\oplus \tau_3$ is {\it not\/} good (note in particular, that the
intersection of the value groups is just $\Z$). Let $H = G_2$; then the
obvious permutation order isomorphism which takes $K$ to $G \oplus H$
takes $\tau_1 \oplus \tau_2 \oplus \tau_3$ to $\tau_1 \oplus \tau_3 \oplus
\tau_2$, hence the latter is good. But with $\sigma = \tau_1 \oplus
\tau_3$ and $\tau = \tau_2$, we have that $\sigma $ is not good whereas
$\sigma \oplus \tau$ (and $\tau$) is good.

To obtain an example wherein neither $\sigma $ nor $\tau$ is good, let
$G_4$ be another copy of $G_3$, set $G = G_1 \oplus G_3$ and $H = G_2
\oplus G_4$ (with the strict orderings of course); $\tau = \tau_2 \oplus
\tau_4$ is for the same reason as $\sigma= \tau_1 \oplus \tau_3$, not
good, but their direct sum is good.
\qed

\comment
A similar family of examples, perhaps simpler, is given below.

\Lem Example \fiveig. A family of dense subgroups of $\R$, $\brcs{G_n}_{n\in \Z}$, each
containing $1$, and normalized traces $\brcs{\tau_n}$ (which are the
inclusion maps $G_n \to \R$), \st for all $m\neq n$, $\phi_{m,n} = \tau_m
\oplus \tau_n$ is {\it not\/} an order unit good trace on the strict direct sum
$G_m \oplus G_n$, but for all $n(1) < n(2) < n(3) < \dots < n(k)$ with $k
\geq 3$, the trace $\phi_{n(1), n(2), \dots, n(k)} = \tau_{n(1)} \oplus
\dots \oplus \tau_{n(k)}$ is good  as a trace on $G_{n(1)} \oplus \dots
\oplus G_{n(k)}$.

\Pf Set $G_n = \Z + (\sqrt 3 + n\sqrt 2) \Z$, viewed as a subgroup of
$\R$; $\tau_n$ will simply be the identification of an element of $G_n$
with the real number it represents. We note that the intersection of any
two of these subgroups of $\R$ is just $\Z$; hence $\tau_m \oplus \tau_n$
is not an order unit good trace on the strict direct sum (and therefore
simple dimension group) $G_m \oplus G_n$.

Now we prove the second statement for $k = 3$, and then the rest follows
by induction, using Lemma~ \fivfou\ above (since if $k \geq 3$,  $\phi_{(n(i))}(\oplus
G_{(n(i))}) = \Z + \sqrt 3 \Z + \sqrt 2\Z$). Pick integers $m < n < p$,
set $\phi = \tau_m \oplus \tau_n \oplus \tau_p$ on $G  =G_m \oplus G_n
\oplus G_p$. Then $\tau_m$, etc, are the coordinate functions on $G$
arising from the obvious embedding in $\R^3$. With this identification,
$\phi$ is matrix multiplication by the row $(1,1,1)$  and $\phi^{\vdash}$
is thus the real linear  span of $v_1 = (1, 0, -1)$ and $v_2 = (0, 1,
-1)$. It is trivial to check that $\ker \phi$ contains the $\Z$-span of
$\brcs{v_1, v_2, v_3 = ((n-p)(\sqrt 3 + m \sqrt 2), (p-m)(\sqrt 3 + n
\sqrt 2), (m-n)(\sqrt 3 + p \sqrt 2))}$.
We note that $v_3 = \lambda v_1 + \mu v_2$ where $\lambda = (n-p)(\sqrt 3
+ m \sqrt 2)$ and $\mu = (p-m)(\sqrt 3 + n \sqrt 2)$.

The vector space isomorphism $\phi^{\vdash} \to \R^2$ obtained by  sending
$v_i \mapsto e_i$ (where $e_1 = (1,0)$ and $e_2 = (0,1)$ are the standard
basis elements), sends $\brcs{v_1, v_2, v_3}$ to $\brcs{e_1, e_2, \lambda
e_1 + \mu e_2}$. Since $\brcs{1,\sqrt 3, \sqrt 2}$ is linearly independent
over the rationals, so is
$\brcs{1, \lambda , \mu}$, and thus the $\Z$-span of the image is dense in
$\R^2$. Hence the $\Z$-span of $\brcs{v_i}$ is dense in $\phi^{\vdash}$,
and since $\ker \phi$ contains $\brcs{v_i}$, $\ker \phi$ has dense image
in $\phi^{\vdash}$, and thus $\phi$ is order unit good. Since $G$ is
simple, $\phi$ is good.
\qed

These examples are getting tiresome, but here is another one, this time
with the $\sigma$ and $\tau$ {\it pure\/} and not order unit good, but
with $\sigma \oplus \tau$ good.

Let $G \subset \R^2$ be the $\Z$-span of $(1,0), (0,1), (\sqrt 2, \sqrt
3)$, equipped with the strict ordering from $\R^2$; similarly, define $H$
by replacing the $3$ by $5$. Both $G$ and $H$ are dense in $\R^2$ (since
$\brcs{1, \sqrt 2, \sqrt 3}$ and $\brcs{1, \sqrt 2, \sqrt 5}$ are linearly
independent over the rationals), and thus both are simple dimension
groups. Let $\sigma$ and $\tau$ be the first coordinate maps on $G$ and
$H$ respectively; then $\sigma(G) = \tau (H) = \Z + \sqrt 2\Z $, and $\ker
\sigma$ and $\ker \tau$ are the cyclic group spanned by $\brcs{(0,1)}$.
Hence neither $\sigma$ nor $\tau$  is order unit good (they are in fact,
ugly). They are pure.

Let $\phi = \sigma \oplus \tau$ on $K = G \oplus H$ with the strict
ordering. Viewing $K \subset \R^4$ in the obvious way, we see that  $\phi$ is
given by taking the inner product with $(1,0,1,0)$. Thus
$\phi^{\vdash}$ is the real span of $\brcs{v_1 = (0,1,0,0), v_2 =
(0,0,0,1), v_3 = (1, 0,-1,0)}$ and  $\ker \phi$
contains the $\Z$-span of
$$
\brcs{v_1, v_2, v_3, v_4 = (\sqrt 2, \sqrt 3,
-\sqrt 2, -\sqrt 5)}.
$$
We can write $v_4 = \sqrt 3 v_1 + \sqrt 5 v_2 +
\sqrt 2 v_3 $. Hence under the isomorphism $\phi^{\vdash} \to \R^3$ given
by $v_i \mapsto e_i$ (the standard basis elements), the image of $\ker
\phi$ is the $\Z$-span of $\brcs{e_1, e_2, e_3, (\sqrt 3, \sqrt 5, \sqrt
2)}$. Since $\brcs{1, \sqrt 2, \sqrt 3, \sqrt 5}$ is linearly independent
over the rationals, the resulting group is dense (a critical group, that
is, a free abelian group dense in a real vector space and of minimal
possible rank). Hence $\ker \phi$ has dense range in $\phi^{\vdash}$, and
thus $\phi$ is order unit good, and since $K$ is simple, $\phi$ is good.
\qed

What makes the examples work is that the trace $\sigma$ has $\ker \sigma$
nonzero, but not dense in $\R$ (in particular, $\sigma \otimes 1$ on $G
\otimes \Q$ is good).

If $G$ and $H$ are simple (noncyclic) dimension groups with finite
dimensional trace spaces, and the obstruction of the example is avoided,
the $\sigma \oplus \tau$ being good does imply $\sigma$ and $\tau$ are
good. Specifically, we assume $\ker \sigma$ and $\ker \tau$ have closure
in their respective affine spaces, which are real vector spaces (in
particular, no discrete direct summands).
\endcomment

\comment
\Lem Lemma. Suppose $(G,u)$ and $(H,v)$ are simple dimension groups with
finitely many pure traces. Let $\sigma$ and $\tau$ be (unnormalized)
traces on $G$ and $H$ respectively, and suppose that $\phi = \sigma \oplus
\tau$ is an order unit good trace on the strict direct sum $ G\oplus H$.
If the closures of the images of $\ker \sigma$ and $\ker \tau$ in $\Aff
S(G,u)$ and $\Aff S(H,v)$ respectively are both vector spaces, then
$\sigma$ and $\tau$ are good. In particular, this applies when $G$ and $H$
are rationally divisible.

\Lem Corollary. Suppose $(G,u)$ and $(H,v)$ are simple dimension groups
with finitely many pure traces. Let $\sigma$ and $\tau$ be (unnormalized)
traces on $G$ and $H$ respectively, and suppose that $\phi = \sigma \oplus
\tau$ is an order unit good trace on the strict direct sum $K=  G\oplus
H$. Then $\ker \phi \otimes \Q$ is a good trace of $G \otimes \Q$.

\Pf (Lemma). Let $J = \sigma^{-1}(\sigma(G) \cap \tau(H))$; this is a
subgroup of $G$ (which contains $\ker \sigma$). The projection map $\Arrow
\pi_G; G \oplus H.G$ restricted to  $\ker \phi$ maps onto $J$, its kernel
is $0 \oplus \ker \tau$, and we have a short exact sequence,
$$
0 \to 0 \oplus \ker \tau \to \ker \phi \to J \to 0.
$$
Now take the closures of the images, identifying $\Aff S(G \oplus H,
(u,v))$ with $\Aff S(G,u) \times \Aff S(H,v)$; we obtain another sequence
$$
0 \oplus \overline{\widehat{\ker \tau}} \to  \overline{\widehat{\ker
\phi}} \to \overline {J}
$$
Now the middle term, by assumption is $\phi^{\vdash}$, since $\phi$ is
order unit good. The term at left is (by assumption) a real subspace of
$\tau^{\vdash}$. The kernel of the projection map is still the thing at
the left. Assume $G$ has $m$ pure traces, and $H$ has $n$ pure traces, so
that $\Aff K$ is $m+n$-dimensional, and thus $\phi^{\vdash}$ is $m+n-1$
dimensional. The thing at the left has dimension at most $n-1$ (since it
is contained in the $n-1$-dimensional subspace, $\tau^{\vdash}$, with
equality if and only if it equals $\sigma^{\vdash}$.

We thus have the inequality, $m+n -1 \leq \dim \overline{\widehat{\ker
\sigma}} + \dim \overline {J}\R$. Since $J$ is a subgroup of $G$, its
dimension cannot exceed $m-1$ (the completion is \wrt the pseudonorm on
$G$, which is the right norm on the quotient group).
Hence $\dim \overline{\widehat{\ker \sigma}} \geq m-1$; since
$\overline{\widehat{\ker \tau}}$ is a subspace of the $n-1$-dimensional
vector space $\tau^{\vdash}$, and this is equal to it. Hence $\tau$ is
order unit good, and the same argument applies to $\sigma$.
\qed

\Pf (Corollary) Since $\phi$ is order unit good, so is $\phi \otimes 1$ on
$K \otimes \Q = (G \otimes \Q) \oplus (H \otimes \Q)$ (strict order), and
we obtain $\sigma \otimes 1$ is order unit good, etc.\qed

Here is an alternative approach to most of these results, including an improvement.
\endcomment

\SecT 6 Good sets of traces

As in [BeH], a compact convex subset $Y$ of $S(G,u)$ is
order unit good (\wrt $(G,u)$) if given $b \in G^+\setminus \brcs{0} $
($b$ is an order unit of $G$) and $a \in G$ \st $0 \ll \hat a|K \ll \hat
b|K$, there exists $a' \in G$ \st $\hat a |K = \hat a'|K$ and $0 \leq a'
\leq b$. When $Y$ is a face (it need not be; e.g., for any singleton
subset of $S(G,u)$, $\brcs{\tau}$ is good iff the trace $\tau$ is good as
defined for individual traces), $Y$ is order unit good iff $\ker Y:=
\cap_{\tau \in K} \ker \tau$ has dense range in $Y^{\perp} = \Set{h \in
\Aff S(G,u)}{h|K \equiv 0}$. When $G$ is simple, this was defined as good in [BeH].
When $G = \Aff K$ (where $K$ is a Choquet simplex), equipped with the strict ordering, goodness of subsets of $K$ is an interesting geometric property. In Appendix B, we
show that at least when $K$ is finite-dimensional, the good subsets of $K$ are as conjectured in [BeH, Conjecture, section 7].

There is a problem in defining what a good subset $Y$ should be in the non-simple
case. It should be consistent with what has been defined in the simple case (where
good  = order unit good), and the singleton case (whence came the original definition of good); additionally, it would be desirable that  if $Y = S(G,u)$, then $Y$ should be good whenever $G$  is a dimension group \st $\Inf G = \brcs{0}$.

We give a definition of good in more complicated situations, including for
a set of traces; this extends some of the definitions in [BeH]. For any
partially ordered abelian group $H$ and $h \in H^+$, recall the definition of the {\it
interval generated by $h$}, denoted $[0,h]$ (possibly with a subscript $H$
if necessary to avoid ambiguity about the choice of group), to be $\Set{j
\in H}{0 \leq j \leq h}$. Let $(G,u)$ be a dimension group (at this stage,
we really only require that it be a partially ordered unperforated group)
with order unit. Let $L$ be a subgroup of $G$; we say {\it $L$ is a good
subgroup of $G$} if the following hold:
\item{(i)} $L$ is convex (that is, if $a \leq c \leq b$  with $a,b \in L$
and $c \in G$, then $c \in L$), and $G/L$ is unperforated
\item{(ii)} using the quotient map $\Arrow \pi; G. G/L$, the latter
equipped with the quotient ordering, for every $b \in G^+$, $\pi([0,b]) =
[0,\pi(b)]$.

Convexity is required in order that the quotient positive cone be proper,
that is, the only positive and negative elements are zero. Unperforation
is often redundant; it may always be (there are no counter-examples; see
the discussion concerning refinability in [BeH]). The second property says
that for all $b \in G^+$, and for all $a \in G$ \st $0 \leq a + L \leq b +
L$ (or equivalently, $(a + L) \cap G^+$ and $(b-a +L) \cap G^+$ are both
nonempty), there exists $a' \in G$ \st $a - a' \in L$ and  $0 \leq a' \leq
b$. This is  a strong lifting property.

For example, if $\tau$ is a trace, set $L = \ker \tau$; this is convex,
and is a good subgroup of $G$ iff $\tau$ is good (as a trace); in this
case, $G/L$ is naturally isomorphic to a subgroup of the reals, so
unperforation is automatic.

For a subset of $S(G,u)$, $U$, define $\ker U = \cap_{\sigma \in U} \ker \sigma$; for
a subset of $G$, $W$, define $Z(W) = \Set{\sigma \in S(G,u)}{\sigma (w) = 0}$. For good sets,
we may as well assume that $Y = Z(\ker Y)$ at the outset, in other words, $\sigma \in Y$ iff
$\sigma(\ker K) = 0$, since in any reasonable definition for good or order unit good, the candidate
set will satisfy $Y = Z(\ker Y)$. As explained in [BeH], these form the collection of closed sets
\wrt a Zariski-like topology, and also extend the definition of facial topology (relative to $G$) on $\partial_e S(G,u)$
to $S(G,u)$.
If $Y \subset S(G,u)$, set $\tilde{Y} =Z(\ker
Y) = \Set{\sigma \in S(G,u)}{\sigma(\ker Y) = \brcs{0}}$; this is a
closure operation, corresponding to the facial topology and analogous to
the Zariski topology from algebraic geometry. In many cases, we just
assume $Y = \tilde{Y}$ already, since $\ker Y = \ker \tilde{Y}$.

We say
{\it $Y$ is a good subset of $S(G,u)$ \wrt $(G,u)$} if $Y = \tilde{Y}$ and
$\ker Y$ is a good subgroup of $G$.
If $Y = \brcs{\tau}$, and $\tau$ is merely an order unit good trace, then
$\ker \tau$ has dense image in $\tau^{\vdash}$, and this implies $Y =
\tilde{Y}$.

If $L$ is a subgroup of $G$, then we may form $Y \equiv Z(L) =
\Set{\sigma \in S(G,u)}{\sigma(L) = \brcs{0}}$. Then $Y$
satisfies $\tilde {Y} = Y$. However, it can happen that $L$ is a good
subgroup of $G$, but $Z(L)$ is not a good subset of $S(G,u)$ \wrt $G$.

For example, let $(H,[\chi_X])$ be the ordered \v Cech cohomology group of any
noncyclic primitive subshift of finite type. It is known not to be a
dimension group, but is unperforated and has numerous other properties
[BoH1, BoH2]. There exists a dimension group $(G,u)$ \st $H \iso G/\Inf G$ with
the quotient ordering. Set $Y=
S(G,u)$, so that $\ker Y = \Inf G$. Since the quotient $H = G/\Inf G$ is
not a dimension group, it follows from results below that property (ii)
fails. However, $L = \brcs{0}$ is clearly a good subset of $G$, and $Z(L)
=Y$, but $\ker Y = \Inf G$. So $Y$ is not a good subset of $S(G,u)$.

In the definition of a good subgroup, it may be that the relatively strong condition that
$G/L$ is unperforated can be replaced by the much weaker $G/L$ is torsion-free, in the presence of (ii), the lifting property. This is the case in the situation described in [BeH, Proposition 7.6], dealing with simple dimension groups and $L = \ker Y$. There are criteria for the quotient $G/L$ to be unperforated [BeH, Lemma B1], but these are not always easy to verify.

The following is implicit
in [BeH, Proposition 7.6].

\Lem Lemma \sixone. Suppose $(G,u)$ is a dimension group and $L$ is a convex
subgroup of $G$ satisfying (ii). Then $G/L$ with the quotient ordering is
an interpolation  group, and its trace space is canonically affinely
homeomorphic to $L^{\vdash}$. The latter is a Choquet simplex.

\Pf If $0 \leq a + L  \leq (b+ L) + (c+L) $ in $G/L$, first lift $b$ and
$c$ separately to positive elements of $G$; it doesn't hurt to relabel
them $b$ and $c$. Applying (ii) to $0 \leq a+L \leq (b+c ) +L$, we can
find $a' \in [0,b+c]$ \st $a-a' \in L$. Hence $0 \leq a' \leq b+c$; by
interpolation in $G$, we may find $a_1 \in [0,b]$ and $a_2 \in [0,c]$ \st
$a' = a_1 + a_2$. Then $a+L = a' + L = (a_1 + L) + (a_2 + L)$ and $a_i +L$
are positive elements of $G/L$, and each is contained less $b+L$, $c+L$
respectively. Thus $G/L$ satisfies Riesz decomposition. The rest is standard.
\qed

\comment
Functions ...

\Lem Lemma. Suppose $(G,u)$ is a dimension group, and $L$ is a convex
subgroup satisfying (ii). If $G/L$ is simple, then $G/L$ is unperforated.

ARgument in [BeH] ...


The following is an improvement on [BeH, proposition B1], where the
condition (*) was introduced. For a subgroup $L$ of the dimension group
$(G,u)$, define the following property, (*):
\item{(*)} for all $k \in L$, for all $\epsilon > 0$, there exist $k_1,
k_2, k_3 \in L$ \st
$$
\hat k = 2k_1 + k_2 + k_3  \qquad \text{and} \qquad \hat k_2, \hat k_3
\geq -\epsilon \pmb 1.
$$

In [BeH, Prop B1], it was shown that if $L$ is a convex subgroup
satisfying (*) \st $G/L$ is torsion-free and $G$ was a simple dimension
group, then $G/L$ is unperforated.
Here we weaken the hypotheses, by assuming merely that $G/L$ is simple
(rather than $G$). This seems like a minor point, but if $Y \subset
S(G,u)$ consists of faithful traces, then $G/\ker Y$ is simple. The
condition (*) applies if the image of $L$ in $\Aff S(G,u)$ is dense in
$Z(L)^{\vdash}= \Set{h \in \Aff S(G,u)}{h(\sigma) = 0 \text{ for all
$\sigma \in Z(L)$}}$. However, $(*)$ does not behave particularly well on
going down to order ideals, so we extract stronger property, density in a
real vector space in the affine representaiton. This is a technical
modification of the argument in [opcit].

\Lem Lemma. Suppose that $(G,u)$ is an approximately divisible dimension
group, and $L$ is a convex subgroup \st $G/L$ is torsion-free and simple
(as a partially ordered abelian group). Suppose that for every order ideal
with order unit $(I,w)$ of $G$, the closure of the image of $I \cap L$ in
$\Aff S(I,w)$ is a real vector space. Then $G/L$ is unperforated.

\Pf Let $a \in G/L$ and suppose $P(a):= \Set{n \in \N}{na \geq 0}$ is
nonempty; as in
[BeH], since $G/L$ is simple, $P(a)$ must be cofinite in $\N$, and as
there, it suffices to show $2a \geq 0$ implies $a\geq 0$. Suppose $2a
\geq0$; there thus exists $p \in G^+$ \st $2a =  p +L$. Form the order
ideal $I$ of $G$ generated by $p$. As $I/(I \cap   \ker L)$ is thus an
order ideal in $G/\ker L$ (done earlier, or it should have been), and as
the latter is simple, we have that the natural inclusion $I/(I \cap   \ker
L) \subseteq G/\ker L$ is onto. Hence we may write $a = g + (I \cap L)$
where $g \in I$. Also, $h = w$ is an order unit for $I$  and $2a =p + (I
\cap L)$. Thus $2g + k = p$.

Now $p$ is an order unit of $I$, so \wrt to the affine representation
$(I,p) \to \Aff S(I,p)$, $\hat p = 1$; choose $\epsilon < 1$. Since $k \in
L \cap I$ and the latter has dense range in a vector space, there exists
$k_1 \in L \cap I$ \st $\|\hat k - 2 \hat k_1\| < \epsilon$, say $k = 2k_1
+ k_2$ where $k_1, k_2 \in L \cap I$, and $\| k_2\| < \epsilon$. Set $h =
g  + k_1 + k_2$; then
$$\eqalign{
2\hat h &= 2\hat g + 2\hat k_1 + 2\hat k_2 \cr
& = 2\hat g + \hat k + \hat k_2 \cr
&= \hat p + \hat k_2 \cr
& = \pmb 1 + \hat k_2 \gg 0.\cr
}$$
Hence $2h$ is an order unit of $I$, and thus so is $h$. Moreover, $h + (L
\cap I) = g + L = a$, so $a \in (I/(I \cap L))^+$.\qed

Better: just assume that every order ideal contains one of these nicer
order ideals (\wrt which the image is dense in a vector space; begin with
an order ideal drop down no doesn't quite work, since must know at the
outset that 2a lifts to an order unit.

all depends on lifting order units.

unperforation occurs if $L = \ker Y$ for reasonable family Y, e.g., good
in lifting sense.

Next, lifting should imply density for all order ideals with order unit.

general criterion at beginning, [BeH,4.3] w/o purity. There argument
outlined with purity, which doesn't seem to play a role [check most recent
version].

probably have this already: if I is order ideal and for all traces t in
Z(L), t(I) is not zero, then I +  L = G +  L:  first, I +  L is order
ideal in G/L; next, Z is the trace space of  G/\ker L, so if I +\ker L is
a proper order ideal, there exists a trace on G/L that kills it, but this
trace must belong to Z, contradiction. what about order iso property;
given g \in G^+, there exists i \in I \st  i + L = g + L; can we choose i
\in I^+. ? OK if unperforated??

observe that if G/L is simple, to show order units lift [to order units],
it is sufficient to show that for a+L an order unit, there exists v \in
G^{++} \st $v \leq a$ (lift a-v to A \in G^+; then A+v \mapsto  a+L, and
is an order unit. When G/L is simple, this is enough to show the latter is
unperforated (proof: suppose $n(a+L) = na + L \in (G/L)^+$; as $G/L$ is
simple, there exists $N \in \N$ \st $u + L \leq Nna + L$. Hence each of
$(Nn+j)a + L$ lifts to an order unit (lift $(Nn + j)a - u + L$ to a
positive element and add u). rest??
\endcomment

If we attempt the simplest definition possible for goodness of a set, that is, $Y$ is
{\it better\/} (a facetious, but not inapt name) for $(G,u)$ if (i) $Y= Z(\ker Y)$ and (ii) whenever $a \in G$, $b \in G^+$
and $0 \leq \hat a|Y \leq \hat b|Y$, there exists $a' \in G^+$ \st $\hat a'|Y = \hat a|Y$ and $a' \leq b$.
This turns out to be much too restrictive (although it is an interesting property); for example, if $Y = S(G,u)$,
then $Y$ is better implies $G/\Inf G$ (with the quotient ordering; this need not be a dimension group) is
archimedean, which hardly ever occurs; and if $G$ is simple, this is generally stronger than order unit good. If
$Y$ is a singleton, then strong goodness agrees with the original definition of good.

\comment
For a convex (but not necessarily directed) subgroup $L$ of a partially ordered abelian group $G$, we can define
an ordering on $G/L$ via $a + L \geq 0$ iff $a + L \cap G^+ \neq \emptyset$, that is, there exists $a' \in G^+$
\st $a - a' \in L$. Convexity guarantees this is a proper ordering (that is, any element of the quotient that is
both positive and negative is zero), but in general, there is no guarantee of other properties, e.g., $G/L$ need not
be torsion-free, and if it is torsion-free, it need not be unperforated. When $L = \ker Y$ for some subset $Y$ of $S(G,u)$,
then $G/L$ is torsion-free (but need not be unperforated; [BeH, Appendix 2]); moreover, we may replace $Y$ by $Z(L)$, so we
may have assumed at the outset that $Y = Z(\ker Y)$.

Instead, we define a subset $Y$ of $S(G,u)$ to be {\it good\/} if it satisfies
\item{(i)} $Y = Z(\ker Y)$
\item{(ii)} if $a \in G$, $b \in G^+$, and $0 + \ker Y \leq a + \ker Y \leq b + \ker Y$, then there exists $a' \in G^+$ \st
$a' - a \in \ker Y$ and $a' \leq b$.

In the definition, (i) is really redundant, since (ii) depends only on $\ker Y$; it is however, convenient to restrict our attention to
sets with $Y = Z(\ker Y)$. In general, such sets are compact and convex, and can be written as $S(G,u) \cap \Cal L$, where $\Cap L$ is a closed flat; when $Y$ is good (for $G$), it is also a Choquet simplex. Part (ii) of the definition is a strong lifting property. When $Y = S(G,u)$, $Y$ is good if $\Inf G = 0$ (interestingly, there is a dimension group $G$, together with an order isomorphism $G/\ker S(G,u) = G/\Inf G \iso G_0$ where $G_0$ is not a dimension group, and property (ii) fails; we can construct $G_0$ as the ordered Grothendieck group of any noncyclic primitive subshift of finite type [BoH]). If $Y$ is singleton, this notion of good agrees with the original. If $G$ is simple, the definition also agrees with order unit good.

In the simple case, there was also a notion of good \wrt $\Aff S(G,u)$ (in fact, order unit good, since the corresponding dimension group was simple), by imposing the strict ordering on the latter. Then a basic result from [BeH] is that $Y$ is order unit good \wrt $G$ iff $Y$ is order unit good \wrt $\Aff S(G,u)$ (this is a strictly geometric property) and the image of $\ker Y$ in $Y^{\vdash}:= \Set{h \in \Aff S(G,u)}{h|Y \equiv 0}$ under the affine representation
$G \to \Aff S(G,u)$ ($g\mapsto \hat g$, where $\hat g(\tau) = \tau (g)$ for normalized traces $\tau$) is dense.
With respect to $\Aff S(G,u)$, singleton subsets and closed faces are automatically (order unit) good, so if $Y = F$ a closed face of $S(G,u)$, $Y$ is order unit good (\wrt $G$) iff $\ker Y$ has dense image in $F^{\perp}$ (for subsets of $S(G,u)$ that are faces, we use the notation $F^{\perp}$ rather than $F^{\vdash}$, to signal that $F^{\perp}$ is a closed order ideal in $\Aff S(G,u)$ \wrt the usual ordering).


We give a definition of good in more complicated situations, including for
sets of traces; this extends some of the definitions in [BeH]. For any
partially ordered abelian group $H$ and $h \in H^+$, define the {\it
interval generated by $h$}, denoted $[0,h]$ (possibly with a subscript $H$
if necessary to avoid ambiguity about the choice of group), to be $\Set{j
\in H}{0 \leq j \leq h}$. Let $(G,u)$ be a dimension group (at this stage,
we really only require that it be a partially ordered unperforated group)
with order unit. Let $L$ be a subgroup of $G$; we say {\it $L$ is a good
subgroup of $G$} if the following hold:
\item{(i)} $L$ is convex (that is, if $a \leq c \leq b$ with $a,b \in L$
and $c \in G$, then $c \in L$), and $G/L$ is unperforated
\item{(ii)} using the quotient map $\Arrow \pi; G. G/L$, the latter
equipped with the quotient ordering, for every $b \in G^+$, $\pi([0,b]) =
[0,\pi(b)]$.

Convexity is required in order that the quotient positive cone be proper,
that is, the only positive and negative elements are zero. Unperforation
is often redundant; it may always be (there are no counter-examples; see
the discussion concerning refinability in [BeH]). The second property says
that for all $b \in G^+$, and for all $a \in G$ \st $0 \leq a + L \leq b +
L$ (or equivalently, $(a + L) \cap G^+$ and $(b-a +L) \cap G^+$ are both
nonempty), there exists $a' \in G$ \st $a - a' \in L$ and $0 \leq a' \leq
b$. This is a strong lifting property.

For example, if $\tau$ is a trace, set $L = \ker \tau$; this is convex,
and is a good subgroup of $G$ iff $\tau$ is good (as a trace); in this
case, $G/L$ is naturally isomorphic to a subgroup of the reals, so
unperforation is automatic. If $Y \subset S(G,u)$, set $\tilde{Y} =Z(\ker
Y) = \Set{\sigma \in S(G,u)}{\sigma(\ker Y) = \brcs{0}}$; $Y \mapsto \tilde{Y}$ is a
closure operation, corresponding to the facial topology and analogous to
the Zariski topology from algebraic geometry. In many cases, we just
assume $Y = \tilde{Y}$ already, since $\ker Y = \ker \tilde{Y}$. We say
{\it $Y$ is a good subset of $S(G,u)$ \wrt $(G,u)$} if $Y = \tilde{Y}$ and
$\ker Y$ is a good subgroup of $G$.
If $Y = \brcs{\tau}$, and $\tau$ is merely an order unit good trace, then
$\ker \tau$ has dense image in $\tau^{\vdash}$, and this implies $Y =
\tilde{Y}$.

If $L$ is a subgroup of $G$, then we may form $Y \equiv Y(L) =
Z(L):= \Set{\sigma \in S(G,u)}{\sigma(L) = \brcs{0}}$. Then $Y$
satisfies $\tilde {Y} = Y$. However, it can happen that $L$ is a good
subgroup of $G$, but $Y(L)$ is not a good subset of $S(G,u)$ \wrt $G$.

For example, let $(H,[\chi_X])$ be the \v Cech cohomology group of any
noncyclic primitive subshift of finite type, denoted $X$. It is known not to be a
dimension group, but is unperforated and has numerous other properties
[BoH]. There exists a dimension group $(G,u)$ \st $H \iso G/\Inf G$ with
the quotient ordering (unpublished results of Mike Boyle and me). Set $Y=
S(G,u)$, so that $\ker Y = \Inf G$. Since the quotient $H = G/\Inf G$ is
not a dimension group, it follows from results below that property (ii)
fails. However, $L = \brcs{0}$ is clearly a good subset of $G$, and $Y(L)
=Y$, but $\ker Y = \Inf G$. So $Y(L)$ is not a good subset of $S(G,u)$.

\Lem Lemma \sixone. Suppose $(G,u)$ is an approximately divisible dimension group, and $L$ is a convex
subgroup satisfying (ii). If $G/L$ is simple, then $G/L$ is unperforated.

\Pf This is a minor modification of the argument in [BeH, Proposition B1] (which requires simplicity of $G$).

The lifting property (ii) is extremely strong. The following is implicit
in [BeH, xxx].

\Lem Lemma \sixtwo. Suppose $(G,u)$ is a dimension group and $L$ is a convex
subgroup of $G$ satisfying (ii). Then $G/L$ with the quotient ordering is
an interpolation group, and its trace space is canonically affinely
homeomorphic to $L^{\vdash}$. The latter is a Choquet simplex.

\Pf If $0 \leq a + L \leq (b+ L) + (c+L) $ in $G/L$, first lift $b$ and
$c$ separately to positive elements of $G$; it doesn't hurt to relabel
them $b$ and $c$. Applying (ii) to $0 \leq a+L \leq (b+c ) +L$, we can
find $a' \in [0,b+c]$ \st $a-a' \in L$. Hence $0 \leq a' \leq b+c$; by
interpolation in $G$, we may find $a_1 \in [0,b]$ and $a_2 \in [0,c]$ \st
$a' = a_1 + a_2$. Then $a+L = a' + L = (a_1 + L) + (a_2 + L)$ and $a_i +L$
are positive elements of $G/L$, and each is contained less $b+L$, $c+L$
respectively. Thus $G/L$ satisfies Riesz decomposition.

\endcomment
\Lem Lemma \sixthr. Let $(G,u)$ be a dimension group with order unit $u$. If $Y \subseteq S(G,u)$ is good, then $G/\ker Y$ is a dimension group, with
trace space canonically affinely homeomorphic to $Y$.

\Pf As good implies order unit good, $\ker Y$ has dense image in
$Y^{\vdash}$, and thus its closure is a vector space, so that by [BeH;
Corollary B2], $G/\ker Y$ is unperforated. Now suppose $0 \leq a + \ker Y \leq (b
+ \ker Y) + ( b' + \ker Y)$, where the latter two terms are nonnegative.
Hence we may assume $b,b' \geq 0$, and thus $0 \leq a + \ker Y \leq (b+b')
+ \ker Y$ implies there exists $a' \in G^+$ \st $a' + \ker Y = a + \ker Y$
and $a' \leq b+b'$. Riesz interpolation in $G$ yields a decomposition $a'
= a_1 + a_2$ where $0 \leq a_1, a_2$ and $a_1 \leq b$ and $a_2 \leq b'$.
Hence $a+\ker Y = a'+ \ker Y = (a_1 + \ker Y) + (a_2 + \ker Y)$, and $a_1
+ \ker Y \leq b  + \ker Y$, and $a_2  + \ker Y \leq b' + \ker Y$. Thus
$G/\ker Y$ satisfies interpolation.

Any trace $\tau$ of $G/\ker Y$,   normalized at $u + \ker Y$, induces a trace $\tilde \tau$
of $(G,u)$ by composing with the quotient mapping. Conversely, if $\sigma$ is a trace
that kills $\ker Y$, then from the definition, $\sigma \in Y$. Hence the map $S(G/\ker Y, u + \ker Y) \to
S(G,u)$ is one to one and onto, and it is easy to see that it is an affine homeomorphism.
\qed

\Lem Lemma \sixfou. If $Y$ is a  good subset of $S(G,u)$,  $(I,w)$ is an order
ideal of $G$ with its own order unit, and for all $\sigma \in Y$,
$\sigma|I \not \equiv 0 $,  then the map $I/(I \cap \ker Y) \to G/\ker Y$
is an order isomorphism.

\Pf First  we show $I/(\ker Y \cap I)$ is unperforated, by showing the  image of $I$ is  an order ideal in $G/\ker Y $ (which is unperforated, by the preceding). Select $0 \leq a + \ker Y \leq b + \ker Y$, where $b \in I$; we can write $b = b_1 - b_2$ where $b_i \in I^+$, and thus $0 \leq a + \ker Y \leq b_1 + \ker Y$, and now $b_1 \in I^+$. There thus exists $a' \in [0,b_1]$ \st $a-a_1 \in \ker Y$. As $0 \leq a' \leq b_1$ and $b_1 \in I$, it follows that $a_1 \in I^+$, so that $a_1 + \ker Y$ is in the image of $I$; the latter is thus a convex subgroup of $G/\ker Y$. Directedness of the image is trivial, so $I/(I\cap  \ker Y)$ is an order ideal in $G/\ker Y$.

Any order ideal in an unperforated partially ordered group is itself unperforated, so $I/(\ker Y \cap I)$ is unperforated.

If $\sigma \in Y$ and $\sigma(w) = 0$, then $\sigma (I) = 0$, contradicting the property of $Y$; hence $\hat w|Y \gg \delta$ for some $\delta > 0$. Since $G/\ker Y$ is unperforated and its trace space is canonically identified with $Y$, it follows that $w + \ker Y$ is an order unit for $G/\ker Y$. Hence the order ideal generated by $w + \ker Y$ is all of $G/\ker Y$. Hence the image of $I$ in $G/\ker Y$ is onto.

So far, the map $I/(I \cap \ker Y) \to G/\ker Y$ is one to one (by construction), order-preserving (by definition), and now we know that it is onto. To show it is an order-isomorphism, it suffices to show that the image of $I^+$ is all of the positive cone.

Select $b \in G^+$. Then $\hat b|Y \ll m$ for some integer $m$, so there exists an integer $N$ \st $\hat b \ll N\hat w$, and thus $0 \leq b + \ker Y \leq Nw + \ker Y$ (the latter by unperforation, again). By goodness, there exists $a \in [0,Nw] $ \st $a - b \in \ker Y$; $0 \leq a \leq Nw$ implies $a \in I^+$, and it  maps to $b + \ker Y$.
\qed

\comment
The map is clearly a one to one and order-preserving map. It suffices
to show that every element of the positive cone of the right side is in
the image of the positive cone of the left. Select $g \in G^+$. Consider
the restriction to $Y$ of $w$, $\hat w|Y$. It is nonnegative, and if $\hat
w|Y$ had a zero, say $\hat w(\sigma) = 0$, then $\sigma(w) =0$; this would
force $\sigma(I) = \brcs{0}$, a contradiction. Hence $\hat w|Y$ is bounded
below (away from zero). There thus exists a positive  integer $n$ \st
$\hat g|Y \ll n \hat w|Y$. Because $G/\ker Y$ is unperforated, and its trace space is canonically $Y$, this forces $g + \ker Y
\leq n(w + \ker Y)$. By goodness, there exists $g' \in G^+$ \st $g' + \ker
Y = g + \ker Y$ and $g' \leq nw$. The latter forces $g' \in I$, so that
the map $I^+ \to (G/\ker Y)^+$ is onto, and we are done.
\qed
\endcomment

The latter property is the analogue of $\tau(I) = \tau(G)$ for a single
good trace $\tau$ of $G$. If we weakened the hypotheses, say to simply
$\ker Y $ does not contain $I$, then the result is unclear. We have
similar problems with the following characterization when some points of
$Y$ are not faithful.

\Lem Lemma \sixfiv. Let $(I,w)$ be an order ideal of $G$ with its own order unit,
and suppose that every point of $Y$ does not kill $I$. Then the map $\Arrow
\phi_I; Y. S(I,w) $ given by $\sigma \mapsto \sigma/\sigma(w)|I$ is
continuous. If $Y$ is good \wrt $(G,u)$, then $\phi_I (Y) $ is good \wrt
$(I,w)$.

\Pf The restriction map on traces sends every point to a nonzero trace of $I$,
and thus  the map is continuous (as $Y$ is compact, $\inf_{\sigma\in Y}
\sigma(w) > 0$). Suppose $\rho$ is  a normalized trace on $(I,w)$ \st
$\rho|(I \cap \ker \rho) $ is identically zero. Then $\rho$ induces a
trace on $I/(I \cap \ker Y)$, hence is a trace on $G/\ker Y$, and
therefore $\rho$ is the restriction of a trace from $G$, necessarily
killing $Y$. If $r$ is the lifted trace, we must have $r \in Y$, and thus
$\rho \in \phi_I(Y)$. In particular, relative to $(I,w)$, $\phi_I(Y) =
Z(\ker \phi_I (Y))$, and it follows immediately that  $\phi_I (Y)$ is good
\wrt $(I,w)$. \qed

The condition on $Y$ in the next result, that every point be faithful, is rather strong, but makes things easier to deal with. The much weaker
faithfulness condition ($\ker Y \cap G^+ = \brcs{0}$) is innocuous, as we can factor
out the maximal order ideal contained in $\ker Y$.

\Lem Lemma \sixsix. Let $(G,u)$ a dimension group, and $Y$ a subset of $S(G,u)$
for normalized traces $\sigma$, $\sigma |\ker Y \equiv 0$  iff $\sigma \in
Y$, and $\ker Y \cap G^+ = \brcs{0}$.
\item{(a)} The trace space of the quotient abelian group $G/\ker Y$ is
canonically affinely homeomorphic to $Y$.
\item{(b)} If $G/\ker Y$ is unperforated and $Y$ satisfies the additional
condition that every element of $Y$ is faithful, then $G/\ker Y$ is
simple.

\comment
\item{(c)} If the hypotheses of (b) apply, and in addition, the map $I \to G/\ker Y$
is onto, then
\endcomment

\Pf  Let $\phi$ be a normalized trace of $(G/\ker Y,u+\ker Y)$, and let
$\Arrow \pi; G. G/\ker Y$ be the quotient map. Then $\sigma':= \sigma
\circ \pi$ is a normalized trace of $(G,u)$ satisfying $\sigma(\ker Y) =
0$, so $\sigma \in Y$. Thus the map $S(G/\ker Y, u + \ker Y) \to S(G,u)$
given by $\sigma \mapsto \sigma \circ \pi$ has image in $Y$, and is
clearly onto $Y$.

\noindent (a) The map  is obviously one to one, affine, and continuous,
with continuous inverse $Y \to S(G/\ker Y,u+\ker Y)$, so is an affine
homeomorphism.

\noindent (b) Suppose nonzero $ a + \ker Y \geq 0 + \ker Y$; there thus
exists $a' \geq 0$ \st $a' - a \in \ker Y$ (from the definition of the
ordering on the quotient group). If $a + \ker Y$ is not an order unit,
then there exists a normalized trace $\sigma $ on $G/\ker Y$ \st $\sigma
(a+ \ker Y) = 0$ (otherwise, $\hat a|Y $ is strictly positive, and as
$G/\ker Y$
is unperforated, this would imply $a + \ker Y$ is an order unit in $G/\ker
Y$). Then $\sigma' = \sigma \circ \pi$ belongs to $Y$ and $\sigma'(a') =
0$, contradicting $\ker \sigma' \cap G^+ = \brcs{0}$.

Hence every nonzero  element of $G/\ker Y$ is an order unit.
\qed

If in part (b), we drop the unperforated hypothesis, then we can still say
something. From $ a + \ker Y \geq 0 + \ker Y$, we have $0 \leq \hat a|Y $;
if for all positive integers  $m$, $ma + \ker Y$ is not an order unit in
$G/\ker Y$, then there must exist a trace $\phi$ on $G/\ker Y$ \st
$\phi(a') = 0$. As in the argument above, this leads to a contradiction.
So in the perforated case, we obtain there exists $m > 0$ \st $m(a+ \ker
Y)$ is an order unit. If we define simple to mean no proper order ideals,
then the quotient group is simple. (We normally deal with unperforated
order groups, wherein the lack of order ideals is equivalent to every
nonzero nonnegative element being an order unit.)

The following contains a slightly different proof of a variant of [BeH; Lemma 7.1].

\Lem Lemma \sixsev. Let $(G,u)$ be an approximately divisible dimension group, and
let $L$ be a convex subgroup. If $G/L$ is unperforated, then order units
lift. [That is, given $a$ \st $a + L$ is an order unit of $G/L$, there
exists an order unit $v$ of $G$ \st $a - v \in L$.]

\Pf The traces of $G/L$ are the traces of $G$ that kill $L$,
$Z:=Z(L)\subset S(G,u)$. As $a + L$ is an order unit, $\hat a|L \gg
\delta$ for some $\delta > 0$. As $G$ is approximately divisible, there
exists $w \in G$ \st $\delta/3 <\hat w < \delta/2$. Then $(\hat a - \hat
w)|Z \gg \delta/2$; as $G/L$ is unperforated, $a - w + L$ is in $(G/L)^+$.
From the definition of quotient ordering, there exists $ c\in G^+$ \st $c
+ L = a-w +L$. Set $v = c + w$. Then $v+L = a-w+w +L = a+L$; since $v \geq
w$ and $w$ is an order unit, it follows that $v$ is an order unit.
\qed

If we drop   approximate divisibility, we obtain that for all order units $a+L$ of $G/L$, there exists an integer $N$ \st for all $n \geq N$, there exist order units $v_n$ of $G$ \st $v_n - na \in L$. (Instead of using a small order unit $w$, we take $u$ or any other order unit we can find.)

The following gives a general result (without assuming $G/\ker Y$ is
unperforated, but instead, that $Y$ is a face) about lifting order units.

\Lem  Lemma \sixeig. Suppose $Y = Z(\ker Y)$ is a face of $S(G,u)$ \st the image
of $\ker Y$ is dense in $Y^{\perp}$. Let $a \in G$ satisfy $a + \ker Y \geq
0$ and $\hat a|Y \gg \delta$ for some $\delta > 0$. Then there exists $a'
\in G^{++}$ \st $a' + \ker Y = a + \ker Y$.

\Pf From the quotient ordering, there exists $c \in G^+$ \st $c - a \in
\ker Y$. Let $F = \Set{\tau \in S(G,u)}{\tau(c) = 0}$; because $c \in
G^+$, $F$ is a face, and is obviously closed. Since $\hat c|Y = \hat a|Y$,
we must have $F \cap Y = \emptyset$. There thus exists $h \in \Aff
S(G,u)^+$ \st $h|Y \equiv 0$ and $h|F\equiv 1$.

As $h \in Y^{\perp}$, there exist $g_n \in \ker Y$ \st $\hat g_n \to h$
uniformly. Hence $\widehat{g_n + c} \to h + \hat c$ uniformly. The latter
however is strictly positive (since $\hat c \geq 0$ and $\hat c^{-1}(0) =
F$). Hence there exists $n$ \st $\widehat {g_n + c}$ is strictly positive;
as $G$ is unperforated, $a':= g_n + c$ is an order unit of $G$. Its image
modulo $\ker Y$ is $c + \ker Y = a+\ker Y$.
\qed

\Lem Proposition \sixnin. Suppose that $(G,u)$ is a nearly divisible dimension group,
and $Y = Z(\ker Y)$ is a subset of $S(G,u) $ \st for all $\sigma \in Y$,
$\ker \sigma \cap G^+ = \brcs{0}$. Suppose that either $Y$ is a face or
$G/\ker Y$ is unperforated. Then $Y$ is  good (\wrt $(G,u)$) iff
\item{(a)} $\phi_I (Y)$ is order unit good for all order ideals $I$ having
their own order unit, and
\item{(b)}  for every nonzero  order ideal $I$, $I + \ker Y = G + \ker Y$.

\Rmk Condition (b) is just a restatement of the map $I/(I\cap \ker Y) \to
G/\ker Y$ being onto. It does not require the stronger property, that it is an order
isomorphism.

\Pf  {\it Sufficiency of the conditions.} Suppose $b \in G^+$ and $a \in
G$ and in addition, $0 \leq a + \ker Y \leq b + \ker Y$. Let $I \equiv
I(b)$ be the order ideal generated by $b$ ($I(b)= \Set{ g\in G}{\exists N
\in \N \text{ \st } -N b \leq g \leq Nb}$). By (b), there exists $a_1 \in
I(b)$ \st $a_1+ \ker Y =  a + \ker Y$. Since $I/(I \cap \ker Y)$ is simple,
$ 0 \leq a_1 + \ker Y \leq b + \ker Y$ entails either $a_1 + \ker Y = 0 +
\ker Y$ or $a_1 + \ker Y$ is an order unit. In the former case, set  $a' =
0$.

Otherwise, if  $Y$ is a face, there exists $a_2 \in I^{++} $
\st $a_2 + \ker Y = a_1 + \ker Y$. Similarly, either $b  + \ker Y = a_1 +
\ker Y$ (in which case, we take $a' = b$) or the difference $b + \ker Y -
(a_2 + \ker Y)$ is an order unit in $I/(\ker Y \cap I)$.

If $G/\ker Y$ is unperforated, then $I/(I \cap \ker Y)$ is unperforated (follows from
$I$ being an order ideal in $G$), and applying Lemma \sixsix(b) to $\phi_I (Y)$, is simple with
trace space canonically $\phi_I (Y) $. This means that the order-preserving one to one
and onto map $I/(I \cap \ker Y)\to G/\ker Y$ induces an affine homeomorphism on their respective
trace spaces; since the image in their affine function representations are the same, that of $I/(I \cap \ker Y)$
has dense range, and being simple (and $\phi_I(Y)$ being a simplex), the latter is a simple dimension group.
A one to one order-preserving group homomorphism between simple dimension groups which induces an affine homeomorphism
on the trace spaces is necessarily an order isomorphism.

Thus in either case, we have $0 \ll \hat a |Y \ll \hat b|Y$; now order unit
goodness of $(I(b), b)$ yields $a' \in I^+$ \st $a' \leq b$.

Necessity of the conditions follows from the preceding results.
\qed

Now we briefly examine examples in $R_P$. When $R$ is a partially ordered commutative unperforated ring with $1$ as
an order unit, every closed face of $S(R,1)$ is uniquely determined by its
extreme points and these form a compact subset of $X = \partial_e S(R,1)$
(and conversely, every closed subset of $X$ yields a closed face in this
way). Thus, as a preliminary question, we can ask when the closed face
obtained from the closed subset $Y$ of $X$ is good (for $R$) or order
unit good. We say $Y$ {\it generates an\ \ \  \/\paren{order unit\/} good face\/}
when this occurs.

It is easy to see that $Y$ generates an order unit good face for $R$ iff
for all pure traces $\sigma \not\in Y$, $\sigma|\ker Y $ is not
identically zero (we define $\ker Y = \cap_{\tau\in Y}\ker \tau$, as
usual).

To see this, if $Y$ generates an order unit good face for $R$, then $\ker
Y$ has dense image in $\Ann Y:= \Set{f \in C(X,\R)}{f|Y \equiv 0}$. There
exists $f \in \Ann Y$ \st $f(\sigma) \neq 0$, and there exist $a_n \in
\ker Y$ \st $\hat a_n \to f$ uniformly, so there exists $a \in
\brcs{a_n}$ \st $ 0 \neq \hat a(\sigma) = \sigma(a)$, hence $\sigma|\ker
Y $ is not identically zero.

Conversely, suppose $\sigma(\ker Y) \neq \brcs{0}$ for every $\sigma \in X
\setminus Y$. It is trivial that $\ker Y$ is an ideal of $R$ (not
generally an order ideal), so its closure in $C(X,\R)$ is a closed ideal
thereof, hence of the form $\Ann Z$ for some closed $Z \subset X$.
Obviously $Y \subset X$, but if $\sigma \in Z \setminus Y$, there exists
$a \in \ker Y$ \st $\sigma(a) \neq 0$, so that $\hat a \not \in \Ann Z$, a
contradiction. Hence $Z = Y$, so $\ker Y$ has dense image in $\Ann Z$, and
thus $Y$ is order unit good for $R$.

\Lem Lemma \sixten. Let $R$ be a partially ordered unperforated approximately
divisible commutative ring, and let $Y$ be a compact subset of the set of
faithful pure traces. Let $(I,v)$ be a nonzero order ideal with its own
order unit.
\item{(a)} The set $Y$ maps by normalized restriction to a compact set of
pure faithful traces on $(I,v)$, $Y_I$.
\item{(b)} If the closed face generated by $Y$ is order unit good for $R$,
then the closed face of $S(I,w)$ generated by $Y_I$ is order unit good for
$(I,v)$.

Good sets for $R_P$ (several variables) corresponding to faces (that
is, closed subsets of  the pure trace space are highly dependent on the
choice of coefficients. For example, as we will see below, if $V$ is the variety given by $f =
(x-3)^2 + (y-3)^2 - 1$, the circle of radius one centred at $(3,3)$ and $P
= c_0 + c_1 x + c_2 y$, then $V$ (or its corresponding face in $S(R_P,1)$
is order unit good, but not good, no matter what the choice of (positive)
integers $c_0, c_1, c_2$. On the other hand, if $P_1 = P\cdot Q$ where $Q
= c + x f + y g + xy h$ where $f$ is a polynomial in $x$ with no negative
coefficients  \st $(x-3)^2 + 8$ divides some power of $c+ xf$ (such
exist!),  $g$ is a polynomial in $y$ \st $(y-3)^2 + 8$ divides some power
of $c + y g$, and $h$ is a polynomial in $xy^{-1}$ \st $ (1 + X^2)$
divides some power of $h(X)$, then $V$ is a good subset for $R_{P_1}$ (the
conditions on the coefficients of monomials appearing in the  faces of the
Newton polytope described by the divisibility condition are necessary and
sufficient for \sixnin(b) to apply; however they are also extremely
complicated).


Now we specialize to $R = R_P$ or $R_P \otimes \Q$, and to avoid severe
complications, also assume that the compact $Y$ consists of pure faithful
traces (that is, $Y $ is a compact subset of the positive orthant,
$(\R^d)^{++}$, after identifying the pure faithful traces with points of
the orthant). Then $\ker Y = \Set{f/P^k \in R_P}{f|Y \equiv 0}$. Recalling
that for $f \in \Z[x_1, \dots, x_d]$, $f/P^k \in R_P$ means there exists
$l$ \st $\Log fP^l \subseteq \Log P^{k+l}$, which means we may as well
assume $\Log f \subseteq \Log P^k$.

Hence $Y$ is order unit good for $R$ iff whenever $\sigma$ is a pure trace
not in $Y$, $\sigma|\ker Y \neq 0$. If we restrict $\sigma$ to the
faithful pure traces, then we deduce a necessary condition:
If $Y \subset (\R^d)^{++}$ is compact, then $Y$ is order unit good for
$R_P$ implies
$$
ZI(Y) \cap (\R^d)^{++}) = Y.
$$
That is, intersecting the Zariski closure of $Y$ with $(\R^d)^{++}$ gives
no new points. In the singleton case, we have see that this condition,
real isolation, is sufficient. However, for general compact $Y$, it is no
longer sufficient.

In fact, examples to illustrate this are ubiquitous (when $d > 1$). The
very simplest one I could think of is the following. Let $P = 1 + xy + x$
(the coefficients, here all ones, are not terribly important); then $\Log P$ is the
triangle with vertices $\brcs{(0,0), (1,1),(1,0)}$, and as rings $R_P \iso
\Z[X,W]$ (the pure polynomial ring in two variables) via the
transformation $X = x/P$ and $W = xy/P$. Let $f = (x-3)^2 + (y-3)^2 -1$,
so $Z(f)\cap \R^2$ is the circle of radius one centred at $(3,3)$, and we
set $Y$ to be this circle, sitting inside the positive quadrant of $\R^2$.
In particular, $\Log f= \brcs{(0,0), (1,0), (2,0), (0,1), (0,2)}$.
It is trivial that $ZI(Y) \cap (\R^2)^{++} = Y$. However, there exists
$\sigma \in \partial_e S(R_P,1)\setminus Y$ \st $\sigma | \ker Y =0$.

Explicitly, $\sigma$ is the pure trace corresponding to the extreme point
of $\cvx \Log P$ given by $(0,0)$; $\sigma (g/P^k) =
(g,x^{0,0})/(P,x^{0,0})^k$. Suppose $a = h/P^k \in R_P$; we may assume $\Log h
\subseteq \Log P^k$. If $r \in Y$ implies $h(r)/P^k (r) = 0$, that is,
$\tau(a) = 0$ for all $\tau \in Y$, then $h|Y \equiv 0$ (since $Y $ is in
the positive quadrant, $P|Y$ vanishes nowhere). Hence there exists $e \in
\Q[x,y]$ \st $h = e\cdot f$ (as $I_{\Q} (f) = f\Q[x,y]$); multiplying by a
positive integer $N$, we may assume $N h = e\cdot f$ where $e \in
\Z[x,y]$.

We claim that this forces $h(0,0) = 0$, that is, its constant term must be
zero, from which it would follow that $\sigma(a) =0$, showing that $\ker Y
\subset \ker \sigma$, as desired. If $h(0,0) \neq 0$, then as $\Log h
\subseteq \Log P^k$, we would have to have $(0,0) \in \Log h$, and in
particular, this point is an extreme point of $\cvx \Log h$. Since $(0,0)$
is also an extreme point of $\cvx \Log f$, it easily follows that $(0,0)$
is an extreme point of $\cvx\Log e$ (I am used to working with Laurent
polynomials, hence this complicated argument, rather than the simple
observation about evaluation). Now consider the coefficients of $e$ and of
$f$ restricted to the line $x = 0$ (that is, throw away all the monomials
with a power of $x$), $e_0$ and $f_0 = y^2 - 6y + 17$. The product is not
zero, and cannot be a single monomial (since $f_0$ is not), hence there
must be, in addition to the constant term, a term of the form $y^j$ in the
product. It is easy to check that this forces $(0,j) \in \Log e\cdot f =
\Log h$. However, $\Log P^k$ is contained in the lattice cone generated by
$\brcs{(0,0), (1,1), (1,0)}$, which does not contain $(j,0)$. This
contradicts $\Log h \subseteq \Log P^k$.

This example does not depend on the coefficients in $P$, that is, we could
just as well have taken $P = 2 + 3xy + 5x$ (which guarantees that $R_P$ is
approximately divisible), nor whether we take $R_P$ or $R_P \otimes \Q$.

In contrast, if we take the same $f$, but $P = 2 + 3 x + 5y$ (or with any
other positive coefficients), then $f/P^2 \in R_P$ and for all
non-faithful pure $\sigma$, $\sigma (f/P^2) > 0$, hence the same $Y$ is
now order unit good for $R_P$. This is part of a more general criterion.

Let $h$ be a polynomial in $d$ variables, and let $S$ be a finite set of
lattice points in $\Z^d$, and $K(S) = \cvx S$. Suppose $F$ is a proper
face of $K(S)$, and $\Log h \subseteq kS$ (the set of sums of $k$ elements
of $S$). We define the {\it facial polynomial of $h$ relative to $F$ and
$k$, $h_{F,k} $} by throwing away all the terms in $x^w$ of $h$ for which
$w \not \in kF$. In case $ S = \Log P$, we can form the element
$h_{F,k}/(P_F)^k \in R_{P_F}$ (in fewer variables, the number being the
dimension of $F$). This yields a positive homomorphism $R_P \to R_{P_F}$
as described in [H1].

Let $Y$ satisfy $ZI(Y) \cap (\R^d)^{++}$, and form the ideal $I(Y)$ of
$\Z[x_1, \dots, x_d]$. Let $P$ be a projectively faithful polynomial in
$\Z[x_1, \dots, x_d]$. We say that {\it $Y$ can be fitted \wrt $P$} if
there exists a polynomial $h \in I(Y)$ \st
\item{(a)} $\Log h \subseteq \Log P^k$ for some $k$
\item{(b)} for every proper face $F$ of $\cvx \Log P$, $h_{F,k}$ has no
negative coefficients.

This depends on $\Log P$, but not very much on the coefficients $P$, as follows
from [H2, Proposition II.5].

Condition (b) can be somewhat weakened, since we are permitted to multiply
numerator and denominator of $h/P^k$ by powers of $P$, and apply eventual
positivity criteria, e.g., [H1A]. The condition is equivalent to
for all pure $\sigma$ that is not faithful, there exists $h \in I(Y)$ \st
$\sigma(h/P_k) > 0$. For example, with $\Log P = \brcs{(0,0), (0,1),
(1,0)}$ and $Y$ the circle in $(\R^2)^{++}$ of radius $1$ centred at
$(3,3)$, $Y$ is fitted \wrt $P$. Just observe that $f$ has the three
facial polynomials (corresponding to the three edges of $\cvx \Log P$ (the
extreme points take care of themselves, so we need not worry about the
zero-dimensional faces), $(x-3)^2 + 17$, $(y-3)^2 + 17$, $x^2 + y^2$. If
we multiply the first two by a sufficiently high power, say $N$, of $1+x$
(respectively $(1+y)$), the outcome will have no negative coefficents. It
follows that if $h = P^N f$, then $h$ will be positively fitted \wrt $P$,
with $k = N+2$.

Now the following is practically tautological.

\Lem Proposition \sixele. Let $P$ be a faithfully projective element of $\Z[x_i]
$, and $Y$ a compact subset of $((\R)^d)^{++}$. Then $Y$ generates an
order unit good face for $R_P$ (and simultaneously for $R_P \otimes \Q$)
iff
\item{(i)} $ZI(Y) \cap (\R^d)^{++} = Y$ and
\item{(ii)} $Y$ can be fitted \wrt $P$.

Conditions on $Y$ to guarantee property (b) of Proposition \sixnin\ seem to be very difficult, involving divisibility of polynomials (and so are highly dependant on the actual coefficients). So goodness of subsets of $\partial_e S(R_P,1)$ is still problematic.

\comment
\Lem Lemma. Let $R$ be as above. Suppose that $K$ is a compact convex
subset of $S(G,u)$.
\item{(a)} If $K$ is order unit good for $G$, then for every order ideal
$(I,v)$ with its own order unit \st $K|Y$ is not zero, the normalized set
of nonzero traces appearing as the restrictions to $I$, $K_I$, is order
unit good.
\item{(b)} If
\endcomment

\SecT Appendix 1.
Order unit good traces on $\Z^k$

The criteria for goodness of traces on nearly divisible dimension groups
depend on  order unit goodness; and the usefulness of the former is a
consequence of the relatively simple characterization of order unit good
traces on approximately divisible dimension groups, namely density of the
image of $\ker \tau$ in $\tau^{\vdash}$ via the affine representation of
$(G,u)$.

To obtain useful criteria for goodness on a larger class of dimension
groups, it would be helpful to find an analogous characterization  of order unit
goodness in the presence of discrete traces. In this appendix, we consider
the  extreme dimension groups with discrete traces, namely
the simplicial ones, $\Z^k$, with the usual ordering. It is already known
that up to scalar multiple, the only good traces are given by left
multiplication by a $0-1$ vector (that is, the entries consist only of
zeros and ones) [H6, Lemma 6.2].

With the current definition of order unit good (really intended for
approximately divisible groups), the order unit good traces on $\Z^k$ can
be characterized, but the characterization makes it difficult to see how
to obtain goodness criteria, as we did in the nearly divisible chase.

Let $v \in (\R^{k \times 1})^+ \setminus \brcs{\pmb 0}$; then $v$ induces
a trace on $\Z^k$, via left multiplication, $\Arrow \phi_v; \Z^k . \R$
sending $w \mapsto vw$ (we think of $\Z^k$ as a set of columns, so matrix
multiplication makes sense). Obviously we can replace $v$ by any positive
real multiple of itself without changing properties such as goodness or
order unit goodness. In addition, we may apply any permutation to the
entries, with the same lack of bad consequences. We may also discard any
zeros (reducing the size of the vectors, that is, decreasing $k$)

Suppose $v $ has only integer entries; then we may order the nonzero
entries, so that
$$v = (n(1), n(2),\dots, n(r); 0, 0, \dots, 0) \quad\text{ where }
n(1) \leq n(2) \leq \dots.$$ We may also assume that $\gcd\brcs{n(i)} = 1$.

\Lem Lemma \Aone. With this choice of $v$, $\phi_v$ is order unit good iff $n(1)
= 1$ and for all $r \geq j > 1$, $n(j) \leq 1 + \sum_{i<j} n(i)$.

\Pf Assume $v$ is in the form indicated, and $\phi_v$ is order unit good. Since
$\gcd\brcs{n(i) }=1$, there exists a vector $w$ \st $vw =1$. Set $u = (1,1,1,\dots,1)$;
we have that $u$ is an order unit, hence it is $\phi_v$-order unit good. Since $vu > 1$ (unless
$v = (1,0,0,\dots,0)$ which is trivially good), there must exist $w_0 \in (\Z^{d})^+$
\st $v w_0 = vw = 1 < vu$. Since the nonzero entries of $v$ are increasing, this forces the smallest
one, $n(1)$, to be $1$. Hence $n(1) = 1$.

Since $vu = \sum n(i):= N$, and there exists $w \in \Z^k$ \st $vw = 1$, there exists for each $s$
with $1 < s < N$
$w_s \in \brcs{0,1}^k$ (as $0 \leq w_0 \leq u$) \st $vw_s = s$, by order unit goodness of $u$.
Now suppose that for some $j$, $n(j) > 1 + \sum_{i < j} n(i)$. Then $n(j) -1$ cannot be realized as a
sum of $n(i)$s (using at most one for each choice of $i$), since $n(j) -1 > \sum_{i < j} n(i)$, and $n(j) \leq n(j')$ for all $j' > j$ (if
there are any such $j'$). Hence no such $w_0$ can exist.

Thus, if $u$ is $\phi_v$-order unit good, then the constraint on growth must hold.

Conversely, suppose the inequalities hold. It is then an easy induction argument (on $r$, augmenting the vector by adjoining $n(k+1)$)
to show that $u$ is $\tau_v$-order unit good, by realizing every integer in the interval
$(0,N)$. Finally, to show that every order unit is $\phi_v$-order unit good ($u$ was the smallest choice), it suffices to show that if we add
a single one to a $\phi_v$-order unit good vector, the outcome is again $\phi_v$-order unit good.
\qed

In particular, the choices for $v$,  $(1,2, 4, 8, 16)$ and $(1,1,1,4)$
yield order unit good traces, but $(1,3) $ and $(1,1,1,5) $ do not. This
rather complicated set of conditions, when applied to order ideals in
dimension groups that have a simplicial quotient by an order ideal, makes
order unit goodness likely unusable for the purposes we had in mind.

\Lem Lemma \Atwo. If $\Arrow \phi_v; \Z^k. \R$
is an order unit good trace, then up to scalar multiple, $v \in (\Z^{k
\times 1})^+$.

\Pf In $\Z^k$, all intervals of the form $[0,u]$ (where $u$ is an order
unit) are finite sets. If there were an irrational ratio among the nonzero
entries of $v$, we would obtain $\phi_v(\Z^k) \cap [0,N]$ is infinite, for
any positive integer $N$. If order unit goodness held, this would be
impossible. Hence all the ratios are rational, and it easily follows that
after suitable scalar multiplication, we can convert $v$ to an integer
row.
\qed

\Lem Proposition \Athr. Let $v$ be an element of  $ \R^{k})^+\setminus
\brcs{\pmb 0}$. Then $\phi_v$ is an order unit good trace, iff up to
scalar multiple and after rearrangement so that $v = (n(1), \dots, n(r);
0, 0, \dots)$ with $n(i-1) \leq n(i)$, we have $n(i) \in \N$, $n(1) = 1$,
and for all $1 < j \leq r$,
$$
n(j) \leq 1 + \sum_{i<j} n(i).
$$

\SecT Appendix 2. Good simplices

In the finite-dimensional case, we verify a conjecture from [BeH, section 7] that good subsets of Choquet simplices are obtained as coproducts of faces with singleton subsets of disjoint faces.

Let $K$ be a Choquet simplex. A nonempty subset $J$ of $K$ is said to be
{\it good\/} (following [BeH]) if it satisfies the following (redundant set of) properties:
\item{(i)} $J$ is a (compact) Choquet simplex
\item{(ii)} there exists a closed flat $\Cal L$ \st $J = \Cal L \cap K$
\item{(iii)} if $a \in \Aff (J)^{++}$ and $b \in \Aff(K)^{++}$
are such that $a \ll b|J$, then there exists $a' \in \Aff(K)^{++}$ \st
$a'|J = a$ and $a' \ll b$.

\noindent We denote this relationship between $J$ and $K$, $J \Gd K$ (there is an uppercase {\it G\/} inside the inclusion sign). If $F$ is a
closed face of $K$, we denote it $F \ideal K$. A question arising out of
[BeH] is to characterize good subsets of Choquet simplices. For example,
closed faces are good, and singleton sets are also good, and coproducts
(within the category of simplices and good subsets) preserve these
properties. A conjecture was made concerning the structure of good
subsets; we verify this in the case that $K$ is finite-dimensional.

Now (ii) is redundant, and   only the compact convex part of (i) is necessary. This is based on the following simple construction.

If $X$ is a subset of a real vector space, define the {\it affine span of $X$},
denoted $\Asp X$, as the  set of finite sums $\Set{\sum r_i
x_i}{r_i \in \R,\ \sum r_i = 1,\ x_i \in X}$.

If $J$ is a singleton or a line segment, there is (almost) nothing to do. Define $\Cal L_0 = \Asp J$. If there exists $v \in (K \cap \Cal L_0)\setminus J$, we can write $v = \sum \alpha_i v_i - \beta_j w_j$ where $v_i, w_j \in J$, and $\alpha_i, \beta_j > 0$, and $\sum \alpha_i - \sum \beta_j = 1$. We can also arrange that $\cvx \brcs{v_i} \cap \cvx \brcs{w_j} = \emptyset$. Hence for any positive $\eta < 1$, there exists $a \in (\Aff J)^{++}$ \st $1-\eta < a |\cvx\brcs {w_j} < 1$ and $a|\cvx\brcs{a_i} < \eta$. Since $a$ is continuous, it is bounded above, so (iii) applies with some constant $b \in \Aff K$.

Hence there exists $a' \in (\Aff K)^{++}$ \st $a = a'|J$. Evaluating the equation at $a'$, we obtain $0 < a'(w) = \sum \alpha_i a(v_i) - \sum \beta_j a(w_j) < \eta \sum \alpha_i - (1-\eta) \sum \beta_j$. This entails
$\eta \(\sum \alpha_i + \sum \beta_j\) > \sum \beta_j$. Now $\sum \beta_j >0$, since otherwise $v \in J$. Hence we can choose at the outset positive $\eta < \sum \beta_j/\(\sum \alpha_i + \sum \beta_j\)$, which yields a contradiction.

Thus $\Cal L_0 \cap K = J$. If $x_n \in \Cal L_0$ and $x_n \to x \in K$, but $x \not\in J$, there exists a line segment joining $x$ to an element of the relative interior of $J$; it must pass through at least two points in $J$, hence $x \in \Cal L_0$. In other words, with $\Cal L$ equalling the closure of $\Cal L_0$, we have $J = \Cal L_0 \cap K = \Cal L \cap K$.

To check that the compact convex set $J$ must be a simplex if (iii) is satisfied, note that the quotient $\Aff K/J^{\vdash}$ (with the strict ordering on $\Aff K$, $J^{\vdash} = \Set{a \in \Aff K}{a|J \equiv0}$, and the quotient ordering) is order isomorphic to $\Aff J$ (with the strict ordering). But goodness imples ([BeH]) that it satisfies Riesz interpolation, which of course forces $J$ to be a Choquet simplex.

\comment
If there is a point $v \in \partial_e J \setminus \partial K$ and $J \neq \brcs{v}$, then there is an edge from a vertex $v'$ of $J$ to $v$; the line segment joining $v'$ to $v$ can be extended a bit beyond $v$ and still remain in $K$. Hence we obtain $w \in K$ together with positive real $0 < \lambda < 1$ \st $v = \lambda w + (1-\lambda)v'$.

We may find $a \in \Aff (J)^{++}$ \st $a(v') = 1$ and $a(v) = \eta$ (where $\eta $ is to be chosen later). There obviously exists a constant \st $a$ is everywhere less than some constant $M$, so set $b $ to be the constant function $M$, defined on all of $K$. By (iii), there exists $a' \in (\Aff K)^{++}$ \st $a'|J = a$ (we really only require that positive elements on $J$ extend to positive elements on $K$).
Applying $a'$ to the convex linear equation, we deduce $\lambda a'(w) = \eta - (1-\lambda)$. If we choose $\eta < 1 - \lambda$, we  obtain a contradiction.
\endcomment

Let $K'$ and $K''$ be simplices (simplices mean Choquet simplices; but
most of the time we will working in finite dimensions, so simplex means
the usual simplex) sitting inside some common simplex $K$ which in turn is
contained in some topological vector space. Suppose that $\Asp K' \cap
\Asp K'' = \emptyset$; we write this as $K' \wedge K'' = \emptyset$. Then
the closure of $\cvx (K', K'')$ is itself a simplex, and we refer to this
as the coproduct, written $K' \dot \vee K''$ (this is more an internal
coproduct, but we shall not distinguish internal from external). If $K'$
and $K''$ are faces of $K$, sufficient for $K' \wedge K'' = \emptyset$ is
that their intersection be empty (since $K$ is a simplex); in this case,
we say that $K'$
and $K''$ are {\it disjoint.} If $\brcs{K^i}$  is a finite family of subsimplices, then disjointness of the set is defined inductively in the obvious way, so that $\dot\vee_i K^i$ makes sense and is a simplex.

We record elementary properties related to goodness.

\Lem Lemma \Bone. (a) Suppose $J \Gd K$ and $K \Gd L$; then $J \Gd L$.
\item{(b)} If $F \Ideal K$, then $ F \Gd K$
\item{(c)} If $J \Gd K$ and $F \Gd K$, then $J \cap F \Ideal J$ and $J
\cap F \Gd K$ whenever $J \cap F \neq \emptyset$.
\item{(d)} If $J_i \Gd K_i$ for $i=1,2$ and $K_1 \wedge K_2 =\emptyset$,
then $J_1 \dot\vee J_2 \Gd K_1 \dot \vee K_2$.

The crucial result is the following. Its proof rests heavily on
finite-dimensionality, but is  a minor modification of the previous argument.

\Lem Lemma \Btwo. Let $K$ be a finite dimensional simplex, and suppose $J \Gd
K$. Let $J_1$ and $J_2$ be disjoint faces of $J$. Set $F_i$ ($i=1,2$) to
be the smallest face of $K$ that contains $J_i$. Then $F_1$ and $F_2$ are
disjoint.

\Pf It suffices to show that $F_1 \cap F_2 = \emptyset$. If not, the
intersection is a face, hence contains a vertex (that is, extreme point)
of $K$, call it $v$. We may suppose that $ v \not\in J_2$ (since $J_1 \cap
J_2 = \emptyset$). Since $J$ is itself a finite dimensional simplex and
$J_i$ are disjoint faces, for any $\eta > 0$ (which we will specify
later), we may find $a \in \Aff(J)^{++} $ \st $a|J_2 \ll 1-\epsilon$,
$a|J_1 \ll \eta$, and $a \ll 1$ (on all of $J$). Set $b$ to be the
constant function $\pmb 1$ on all of $K$, so that $0 \ll a \ll b|J$.

By goodness, there exists $a' \in \Aff(K)^{++}$ \st $a' \ll b$ and $a'|J =
a$. It is now easy to show that for suitably small $\eta$ (depending on
the boundary measure of elements of $J_i \subset F_i$), this leads to a
contradiction.

Since $v \not\in J$ and $F_2$ is the smallest face containing $J_2$, there
must exist $w \in J_2$ \st $w = \lambda v + \sum_s \lambda_s v_s$ where
$v_s \in \partial_e F_2$, $\lambda > 0$, $\lambda_s \geq 0$ and $\lambda =
1 - \sum \lambda_s$. Evaluating at $a'$, we obtain $\lambda a'(v) = a(w) -
\sum \lambda_s a'(v_s) \geq 1 - \eta - (1 - \lambda)$ (since $a'(v_s) \leq
b(v_s) = 1$).
Thus $a'(v) \geq 1 - \eta/\lambda$.

Now working within $F_1$, again since $F_1$ is the smallest face
containing $J_1$, there must exist $y \in J_1$ \st $y = \mu v + \sum_t
\mu_t y_t$ where $\brcs{v, y_t} \subseteq \partial_e F_1$, $\mu > 0$,
$\mu_s \geq 0$, and $\mu = 1 - \sum \mu_s$. Applying $a'$, we obtain $\mu
a'(v) = a(y) - \sum \mu_t a'(y_t) < \eta$. Hence $a'(v) < \eta/\mu$.

Thus the two inequalities force $\eta/\mu + \eta/\lambda > 1$. We reach a
contradiction if we choose $\eta < 1/(1/\mu + 1/\lambda)$.
\qed

One  obstruction (among several) to extending this to infinite-dimensional simplices is the
fact that the representing measures of relative interior points might
vanish on the intersection of the faces. We would also have to restrict to
closed faces in this case (since otherwise it is not clear that the
smallest face exists), and this will present problems when we want to use
it.

Let $\brcs{F_i}$ be a disjoint collection of  faces  (that is, for all $i$,  $F_i \wedge (\dot\vee_{j\neq i} F_j) = \emptyset$) of the simplex $K$,
and for each $i$, let $v_i$ be a point in the relative interior of $F_i$;
we also assume that the $F_i$ are not themselves singletons. We may form
$J_0:= \cvx \brcs{v_i}$ and $F_0 := \cvx \brcs{F_i}$; of course, this is
the coproduct of $(\brcs{v_i}, F_i)$, and $J_0 $ is thus a good subset of
$F_0$ (since each $\brcs{v_i} \Gd F_i$). As in [BeH], we call the
$(v_i,F_i)$, together with $(F,F)$ (that is, the face $F \Gd F$) {\it building
blocks.} It was conjectured (in the finite-dimensional case) that if $J \Gd
K$, then there exists a face $F$ of $K$, together with a disjoint face
$F_0$ obtained as the coproduct, \st $J = F \dot \vee J_0$; in other
words, that coproducts of the building blocks yield all good subsets;
alternatively, that there is a face maximal $F $ of $K$ sitting inside
$J$, and $J$ is obtained by taking coproducts with respective singleton
sets sitting inside pairwise disjoint faces. This now follows easily.

\Lem Corollary \Bthr. Suppose $K$ is a finite-dimensional simplex and $J \Gd K$.
Then there exist a (possibly empty) face $F$ of $K$ together with a
finite set of faces $F_i$ of dimension at least one \st $\brcs{F,F_1,
\dots}$ is  disjoint, together with $v_i $ in the relative
interior of $F_i$ \st $J = \cvx\brcs{F,v_i}$.

\Pf We proceed by induction on the dimension of $J$. Let $F$ be the convex
hull of all the vertices of $K$ that lie in $J$; these are automatically
vertices of $J$. If this exhausts the vertices of $J$, then $F = J$ and
$F$ is a face (since $K$ is a finite-dimensional simplex), and there is
nothing to do. Of course, $F$ can be empty.

Otherwise, there exists a vertex $v_1$ of $J$ that is not in $\partial_e
K$; necessarily this belongs to a proper face (it cannot be in the
interior, in fact by property (ii), but this can also be proved using only
(i) and (iii)) of $K$, and let $F_1$ be the smallest face of $K$
containing $v_1$. Then $v_1$ is in the relative interior of $F_1$. Let
$J^1$ be the complementary face to $\brcs{v_1}$ in $J$ (that is, the
convex hull of all the other vertices of $J$).

If $J^1$ is empty, then $J = J^1$ is already a singleton, and we are done.

If $J^1$ is not empty, then $J^1 \ideal J $, so $J^1 \Gd J$, and thus by
transitivity, $J_1 \Gd K$. We can apply the previous lemma. Let $F^1$ be
the smallest face of $K$ containing $J^1$; then $F^1 \cap F_1 =
\emptyset$, and thus $J$ decomposes as the coproduct of $J^1$ and
$\brcs{v_1}$ (using faces $F^1$ and $F_1$), so by induction on the
dimension of $J$, we are done.
\qed

The conjecture in the case that $K$ be infinite-dimensional is more complicated, and I have no idea how to proceed.

\long\def\Rf[#1] #2, #3. #4\par%
{\vskip 2pt \itemitem{[#1]} #2, {\it #3,} #4\par\vskip2pt}

\SecT Acknowledgment

Discussions with my colleague Damien Roy concerning the material in section 3 were very helpful.

\SecT References

\Rf [Ak1] E Akin, Measures on Cantor space. Topology Proc,  {24} (1999) 1--34.

\Rf [Ak2] E Akin, Good measures on Cantor space. Trans Amer Math Soc,  {357} (2005) 2681--2722.

\Rf [ADMY] {E Akin, R Dougherty,  RD Mauldin, and A Yingst}, Which Bernoulli measures are good measures?. Colloq  Math,  {110} (2008) 243--291.

\Rf [BeH] S Bezuglyi \& D Handelman, Measures on Cantor set{\/\rm:} the good, the ugly, the bad.
Trans Amer Math Soc (to appear).

\Rf [BoH1] M Boyle \& D Handelman, Ordered equivalence, flow equivalence and ordered cohomology. Israel J Math 95 (1996) 169--210.

\Rf [BoH2] M Boyle \& D Handelman, unpublished drafts. 1993--2013.

\Rf [EHS] {EG Effros, David Handelman, \& Chao-Liang Shen}, Dimension groups and their affine representations. Amer J Math 102 (1980) 385--407.

\Rf [FO] SB Frick \& N Ormes, Dimension groups and invariant measures for polynomial odometers. Acta Appl Math (2013).

\Rf [GPS] {T Giordano,  IF Putnam, \& CF Skau},  Topological orbit equivalence and $C^*$--crossed products. J Reine Angew Math  469  (1995), 51--111.

\Rf [G] KR Goodearl, Partially ordered abelian groups with interpolation. Mathematical Surveys and Monographs, 20, American Mathematical Society, Providence RI, 1986.

\Rf [GH] KR Goodearl \& David Handelman, Metric completions of partially ordered abelian groups. Indiana Univ J Math 29 (1980) 861--895.

\Rf [GH2] KR Goodearl \& David Handelman,   Tensor products of dimension groups and $K_0$ of unit--regular rings. Canad J Math  38  (1986),  no\. 3, 633--658.

\Rf [H1] David Handelman,    Positive polynomials and product type actions of
compact groups. Mem Amer Math Soc  54  (1985),  320, xi+79 pp.

\Rf [H1A] David Handelman, Deciding eventual positivity of polynomials. Ergodic theory and dynamical systems 6 (1985) 57--79.

\Rf [H2] David Handelman,   Positive polynomials, convex integral polytopes, and a random walk problem. Lecture Notes in Mathematics, 1282, Springer--Verlag, Berlin, 1987, xii+136 pp.

\Rf [H3] D Handelman, Iterated multiplication of characters of compact connected Lie groups. J of Algebra 173 (1995) 67--96.

\Rf [H4] D Handelman, Free rank $n+1$ dense subgroups of $\text{\/\bf R}^{n}$ and their endomorphisms. J Funct Analysis  46  (1982), no\. 1, 1--27.

\Rf [H5] D Handelman, Matrices of positive polynomials. Electronic J Linear Algebra 19 (2009) 2--89.

\Rf [H6] D Handelman, Realizing dimension groups, good measures, and Toeplitz factors. submitted. (Formerly known as {\it Equal row and equal column sum realizations of dimension groups,} on Ar$\chi$iv.)

\Rf [HPS] {RH Herman, IF Putnam, \& CF Skau}, Ordered Bratteli diagrams, dimension groups and topological dynamics. Internat J Math  3  (1992),  no\. 6, 827--864.

\comment
\Rf [M] X Mla, A class of nonstationary adic transformations. Ann Inst Henri Poincar Prob \& Stat 42 (2006) 103--123.
\endcomment

\Rf [P] K Petersen, An adic dynamical system related to the Delannoy numbers. Ergodic theory \& dynamical systems 32 (2012) 809--823.

\Rf [Pu] IF Putnam, The C* algebras  associated with minimal  homeomorphisms of the Cantor set.
Pacific J Math 136 (1989) 329--353.

\Rf [R] J Renault, A groupoid approach to C*-algebras. Lecture Notes in Mathematics, 793, Springer--Verlag, Heidelberg, 1980.

\vskip 10pt

\noindent Mathematics Dept, University of Ottawa, Ottawa K1N 6N5 ON, Canada; dehsg\@uottawa.ca

\end

%% file: generic_macros






\font\rm=cmr10 \rm

\font\bf=cmb10
\font\Rm=cmr9 at 11pt
\rm
\font\it=cmsl9 at 10pt
 at 7pt

\font\Rrm=cmr17 at 16pt
   \font\Rm=cmr12 at 11.5pt

\long\def\Pf{\par\noindent {\it Proof.} }
\def\({\left(}
\def\){\right)}
\def\st{such that }
\def\qed{\hfill$\bullet$\vskip 4pt}
\def\quotes#1{{\lq\lq #1\rq\rq}}
\def\brcs#1{\left\{ #1\right\}}

\def\Log{\text{Log\,}}
\def\iso{\cong}
\def\wrt{with respect to }
\def\:{\,:}

\def\ker{\text{ker\,}}

\def\A{{\Cal A}}

\def\C{\text{\bf C}}
\def\T{\text{\bf T}}

\def\Im{{\text{Im}\,}}

\def\R{\text{\bf R}}
\def\N{\text{\bf N}}
\def\Z{\text{\bf Z}}
\def\Q{\text{\bf Q}}

\def\Arrow #1;#2.{#1\:#2 \to }

\def\Set#1#2{\brcs{#1 \left|\vphantom{#1 #2} \right.#2}}



\def\Rrr#1,#2{{\Cal J}_{#1,#2}}
\def\slfrac#1#2{{\raise -.07 ex\hbox{$^{#1}$}}\!/\raise .35 ex \hbox{${}_{#2}$}}
\def\ssf #1/#2{\slfrac {#1}{#2}}

\def\pd #1,#2.{\frac {\partial #1}{\partial #2}}

   \long\def\Lem
#1.#2\par{\vskip4pt{\baselineskip=13pt\font\it=cmsl12 at
11.5pt\Rm
   \noindent {\rm \uppercase{#1}} #2\vskip3pt

   }} 

\long\def\Proclaim #1.#2 \endproclaim{\vskip4pt{\baselineskip=13pt\font\it=cmsl12 at
11.5pt\Rm
   \noindent {\rm \uppercase{#1}} #2\vskip3pt

   }} 

\long\def\remark #1\endremark{\vskip 2pt \noindent {\it Remark\/} #1\par}

\long\def\Sectionhead #1.#2:\par #3{\vskip 4pt \noindent {\bf #1 #2}vskip 2pt\noindent\nospace #3}

\long\def\Title #1\par {\noindent{\Rrm #1}\vskip 9pt}

 \long\def\SubT #1.{\noindent {\it #1\/} } 
 
 \long\def\SecT
#1\par{\vskip 3pt \noindent {\bf #1}\vglue1pt
   \noindent}

\long\def\subtitle #1.{\vskip 2pt \noindent {\it #1}}

\long\def\Rmk#1\par{\vskip 1pt \noindent {\it
Remark.} #1\vskip2pt}

\long\def\Abstract #1\par{{\leftskip= 3 true cm \rightskip = 3 true cm \font\it=cmsl10 \font\rm=cmr10 \baselineskip = 10pt
\parindent=.35 true cm\rm\noindent 
{\it Abstract} #1\vskip 8pt

}}

\long\def\Author #1 \par{\noindent{\it #1}}

%% file: papermacros_entirefunctions
 
\scrollmode\NoBlackBoxes
\magnification=1100
\long\def\Abstract #1\par%
{\vskip .2 true cm{\leftskip 1 true in \rightskip 1 true in \font\rm=cmr8 \rm
\baselineskip=1pt \font\it=cmsl8 \font\bf=cmb10 at 8pt
\parindent=0em {\bf Abstract} #1

}}
\comment
\font\rm=Times at 10pt

\font\bf=TimesB
\font\Rm=Times at 11pt
\rm
\font\it=TimesI at 10pt
\endcomment

\long\def\Pf{\par\noindent {\it Proof.} }
\def\({\left(}
\def\){\right)}
\def\st{such that }
\def\qed{\hfill$\bullet$\vskip 4pt}
\def\quotes#1{{\lq\lq #1\rq\rq}}
\def\brcs#1{\left\{ #1\right\}}
\def\Set#1#2{\brcs{#1 \left|\vphantom{#1 #2} \right.#2}}

\def\C{\text{\bf C}}
\def\T{\text{\bf T}}

\def\Im{\text{Im\,}}

\def\Log{\text{Log\,}}
\def\iso{\cong}
\def\wrt{with respect to }
\def\:{\,:}
\def\Arrow #1;#2.{#1\:#2 \to }


\def\R{\text{\bf R}}
\def\N{\text{\bf N}}
\def\Z{\text{\bf Z}}
\def\Q{\text{\bf Q}}
 
\def\Rrr#1,#2{{\Cal J}_{#1,#2}}

\def\slfrac#1#2{{\raise -.07 ex\hbox{$^{#1}$}}\!/\raise .35 ex \hbox{${}_{#2}$}}
\def\ssf #1/#2{\slfrac {#1}{#2}}

\def\pd #1,#2.{\frac {\partial #1}{\partial #2}}


   \long\def\Title #1\par {\noindent{\Rrm #1}\vskip 9pt}
 \long\def\SubT #1.{\noindent {\it #1\/} }   \long\def\SecT
#1\par{\vskip 3pt \noindent {\bf #1}\vglue1pt
   \noindent}
\long\def\subtitle #1.{\vskip 2pt \noindent {\it #1}}

\long\def\Rmk#1\par{\vskip 1pt \noindent {\it
Remark.} #1\vskip2pt}


\def\oneone{\One.1}
\def\onetwo{\One.2}
\def\onethr{1.3}
\def\onefou{1.4}
\def\onefiv{1.5}
\def\onesix{1.6}
\def\onesev{1.7}

\def\twoone{3.1}
\def\twotwo{3.2}
\def\twothr{3.3}
\def\twofou{3.4}
\def\twofiv{3.5}

\def\threig{2.8}

\def\fouone{4.1}
\def\foutwo{4.2}
\def\fouthr{4.3}
\def\foufou{4.4}
\def\foufiv{4.5}
\def\fousix{4.6}
\def\fousev{4.7}

\def\fivone{5.1}
\def\fivtwo{5.2}
\def\fivthr{5.3}
\def\fivfou{5.4}
\def\fivfiv{5.5}
\def\fivsix{5.6}
\def\fivsev{5.7}
\def\fiveig{5.8}
\def\fivnin{5.9}
\def\fivten{5.10}

\def\sixone{6.1}
\def\sixtwo{6.2}
\def\sixthr{6.3}
\def\sixfou{6.4}
\def\sixfiv{6.5}
\def\sixsix{6.6}
\def\sixsev{6.7}
\def\sixeig{6.8}
\def\sixnin{6.9}

\def\Aone{A.1}
\def\Atwo{A.2}
\def\Athr{A.3}

\def\Bone{B.1}
\def\Btwo{B.2}
\def\Bthr{B.3}
